%Paper: q-alg/9505003
%From: Charilaos Aneziris <aneziris@ifh.de>
%Date: Tue, 2 May 1995 14:07:35 +0200 (METDST)

\magnification=\magstep1
\font\tenrm=cmr10
\font\teni=cmmi10
\font\tensy=cmsy10
\font\tenbf=cmbx10
\font\tentt=cmtt10
\font\tenit=cmti10
\font\tensl=cmsl10
\font\ninerm=cmr9
\font\ninei=cmmi9
\font\ninesy=cmsy9
\font\ninebf=cmbx9
\font\ninett=cmtt9
\font\nineit=cmti9
\font\ninesl=cmsl9
\font\eightrm=cmr7
\font\eighti=cmmi7
\font\eightsy=cmsy7
\font\eightbf=cmbx7
\font\eighttt=cmtt8
\font\eightit=cmti7
\font\eightsl=cmsl8
\skewchar\teni='177
\skewchar\tensy='60
\hyphenchar\tentt=-1
\skewchar\ninei='177
\skewchar\ninesy='60
\hyphenchar\ninett=-1
\skewchar\eighti='177
\skewchar\eightsy='60
\hyphenchar\eighttt=-1
\font\sixrm=cmr6
\font\sixi=cmmi6
\font\sixsy=cmsy6
\font\sixbf=cmbx6
\font\fiverm=cmr5
\font\fivei=cmmi5
\font\fivesy=cmsy5
\font\fivebf=cmbx5
\def\tenpoint{\def\rm{\fam0\tenrm}%
\textfont0=\tenrm \scriptfont0=\sevenrm \scriptscriptfont0=\fiverm
\textfont1=\teni \scriptfont1=\seveni \scriptscriptfont1=\fivei
\textfont2=\tensy \scriptfont2=\sevensy \scriptscriptfont2=\fivesy
\textfont3=\tenex \scriptfont3=\tenex \scriptscriptfont3=\tenex
\textfont\itfam=\tenit \def\it{\fam\itfam\tenit}%
\textfont\slfam=\tensl \def\sl{\fam\slfam\tensl}%
\textfont\ttfam=\tentt \def\tt{\fam\ttfam\tentt}%
\textfont\bffam=\tenbf \scriptfont\bffam=\sevenbf
\scriptscriptfont\bffam=\fivebf \def\bf{\fam\bffam\tenbf}%
\tt
\normalbaselineskip=16pt
\setbox\strutbox=\hbox{\vrule height8.5pt depth3.5pt width0pt}%
\let\sc=\eightrm \let\big=\tenbig \normalbaselines\rm}
\def\ninepoint{\def\rm{\fam0\ninerm}%
\textfont0=\ninerm \scriptfont0=\sixrm \scriptscriptfont0=\fiverm
\textfont1=\ninei \scriptfont1=\sixi \scriptscriptfont1=\fivei
\textfont2=\ninesy \scriptfont2=\sixsy \scriptscriptfont2=\fivesy
\textfont3=\tenex \scriptfont3=\tenex \scriptscriptfont3=\tenex
\textfont\itfam=\nineit \def\it{\fam\itfam\nineit}%
\textfont\slfam=\ninesl \def\sl{\fam\slfam\ninesl}%
\textfont\ttfam=\ninett \def\tt{\fam\ttfam\ninett}%
\textfont\bffam=\ninebf \scriptfont\bffam=\sixbf
\scriptscriptfont\bffam=\fivebf \def\bf{\fam\bffam\ninebf}%
\tt
\normalbaselineskip=9pt
\setbox\strutbox=\hbox{\vrule height8pt depth3pt width0pt}%
\let\sc=\sevenrm \let\big=\ninebig \normalbaselines\rm}
\def\eightpoint{\def\rm{\fam0\eightrm}%
\textfont0=\eightrm \scriptfont0=\sixrm \scriptscriptfont0=\fiverm
\textfont1=\eighti \scriptfont1=\sixi \scriptscriptfont1=\fivei
\textfont2=\eightsy \scriptfont2=\sixsy \scriptscriptfont2=\fivesy
\textfont3=\tenex \scriptfont3=\tenex \scriptscriptfont3=\tenex
\textfont\itfam=\eightit \def\it{\fam\itfam\eightit}%
\textfont\slfam=\eightsl \def\sl{\fam\slfam\eightsl}%
\textfont\ttfam=\eighttt \def\tt{\fam\ttfam\eighttt}%
\textfont\bffam=\eightbf \scriptfont\bffam=\sixbf
\scriptscriptfont\bffam=\fivebf \def\bf{\fam\bffam\eightbf}%
\tt
\normalbaselineskip=8pt
\setbox\strutbox=\hbox{\vrule height7pt depth2pt width0pt}%
\let\sc=\sixrm \let\big=\eightbig \normalbaselines\rm}
\parskip 0pt
\baselineskip 16pt
\null
\vskip 1 cm
\centerline {\bf A NOVEL APPROACH TO KNOT CLASSIFICATION}
\bigskip
\centerline {Charilaos Aneziris}
\centerline {\it Institut f\"ur Hochenergienphysik Zeuthen}
\centerline {\it DESY Deutsches Elektronen-Synchrotron}
\centerline {\it Platanenallee 6}
\centerline {\it 15738 Zeuthen}
\bigskip
\centerline {\bf Abstract}
\bigskip
In a previous paper (Ref. [1]) we suggested some algorithms that could be
useful in
solving the problem of knot classification. Here we continue this discussion
by answering questions raised in that paper and by commenting on practical
aspects related
to running these algorithms as computer programs. Further on, we reach a
solution by improving on
the efficiency of these algorithms and by introducing new ones.
\bigskip
\centerline {1. \qquad INTRODUCTION}
\bigskip
In Chern-Simons topological field theories the observables come out to be
the Wilson lines defined by $$W(R,C)=trP{\rm exp}\int _C A$$ where $R$ is an
irreducible representation of the gauge group $G$ and $C$ is a closed curve.
Due to the topological
nature of these theories, a continuous change on $C$ leaves $W$ invariant.
Therefore $W$ is defined not on the whole space of closed curves $C$, but
on a ``smaller" space of equivalence classes, two curves
being considered equivalent if they can be continuously deformed to each other.
In knot theory such curves are called {\bf ambient isotopic}, and thus the
observables of Chern-Simons topological field theories are related to
ambient isotopy classes which form the set of topologically distinct knots.
\par \medskip
When one considers now the correlation function associated with
$n$ such Wilson lines,
$<W(R_1,C_1)...W(R_n,C_n)>$, the result comes out to depend exclusively on the
$n$-component link which the $n$ Wilson lines form. In 1989 Witten showed that
when $R$ is the defining representation of $SU(2)$, the correlation function
comes out equal to the {\bf Jones polynomial} of the corresponding link
(Ref. [2]).
More recently it has been suggested in a number of papers that ultraviolet
divergencies of a renormalizable field theory are knot invariants (Ref. [3]).
In these
papers a relationship between Feynman diagrams and knot graphs has been
introduced, through which divergencies are related to properties of knots.
These, as well as a number of other interesting results published in between,
have established
the relevance of knot theory to High Energy Physics.
\par \medskip
It is also interesting to note that knot theory, in addition to High Energy,
has been found applicable to other branches of science such as Low
Dimensional Topology in mathematics (Ref. [4]), the
{\bf Potts model} of Condensed matter and the study of
critical phases in Statistical Physics (Ref. [5]), Quantum Gravity (Ref. [6]),
Chemistry (Ref. [7]), and even in the study of DNA in Biology (Ref. [8]).
In spite however of the extensive study of knot theory and its
applications during recent years, the problem of actually finding these
knots and of classifying them is far from being solved. While a number of
interesting properties of knots are still being obtained and one may find
lists with the simplest knots (sometimes contradicting each other) published in
relevant books (Ref. [9]), one has yet to establish a complete procedure or
method
which could yield all distinct classes of knots while avoiding undesirable
repetitions, in spite of a very good approach discussed by Thistlethwaite in
Ref. [10].
\par \medskip
In Ref. [1] the problem of knot classification was discussed by presenting
two ideas which could be useful in establishing algorithms that
could solve it. The first such idea was based on
providing knots with suitable ``names" dependent on two-dimensional
projections. One had then to check whether a particular ``name" actually
belonged to a knot, and if so, whether that knot had already been named by a
different such name. If those questions could be answered, the problem of
knot classification would have effectively been solved; while the first was
answered and the second had been extensively discussed, the lack of a complete
solution to the issue of avoiding repetitions under distinct names, was the
reason why a definite answer could not be provided. The second idea was based
on the fact that any closed curve may be continuously deformed to one
consisting of linear segments joining neighboring vertices of a cubic lattice.
By translating these vertices into a number sequence which could be then used
as a ``name", one would hope to establish a method to tackle knot
classification.
\par \medskip
Since the appearance of Ref. [1] we have further investigated the
approach related to two-dimensional projections; not only have we been able to
improve on identifying equivalent knots, but by generalising Reidemeister's
{\bf tricolorability test} we have been able to establish methods to show if
two knot projections are distinct. In Section 2 we provide some necessary
definitions and review the relevant material of Ref. [1], so that the reader
may
easily understand the rest of the paper. In Section 3 we present the main ideas
of developing a ``knot equivalence" algorithm, which is an improved version of
the algorithm of Ref. [1]. In Section 4 we discuss problems due to the
finiteness
of the allocated computer memory and running time, and suggest ways to
overcome such difficulties. In Section 5 we discuss how to develop a ``knot
inequivalence" algorithm, while in Sections 6 through 8 we focus on specific
examples of such algorithms. In particular, in Section 7 a new knot invariant
is presented, which emerged while applying a special class of inequivalence
algorithms. Finally in Section 9 we present and discuss the
results we got after running the computer programs.
\bigskip
\centerline {2. \qquad DEFINITIONS AND NOTATIONS}
\bigskip
A {\bf knot} is defined as a one-to-one homeomorphism from $S^1$ to some
manifold $\cal M$.
\par \medskip According to this definition, knots may exist on any manifold
$\cal M$, even though in most manifolds, including all $S^n$ for $n \ne 3$,
all knots are trivial. In this paper we are only considering knots in $S^3$,
the three-dimensional sphere, which may also be viewed as the
compactification of $R^3$. This compactification is understood as follows.
\par \medskip
Let $f(x)=2{\rm tan}^{-1}x$ be a function from the straight line $R^1$ to the
circle $S^1$. Such a function is one-to-one, since all points of $R^1$ are
mapped to some
point of $S^1$; no point of $R^1$ however is mapped to $\theta = \pm \pi$ of
$S^1$ (Figure 1a). By identifying the points $\pm \infty$ of $R^1$ and mapping
them to $\pi$ of $S^1$ we compactify $R^1$. In a similar manner one may
compactify $R^2$ and map the infinity circle of $R^2$ to one of the poles of
$S^2$ (Figure 1b).
\par \medskip
A compactification of $R^3$ is a process similar to the one described above,
where the infinity sphere of $R^3$ is mapped to some point of $S^3$, although
it is impossible to draw in three dimensions a figure similar to Figures 1a or
1b. By defining a knot as stated above (with ${\cal M} = S^3$), one may
distinguish two possibilities.
\par \medskip
a) The knot does not pass through the ``infinity" point of $S^3$. In such a
case the knot is a closed curve in $R^3$. (Figure 2a)
\par \medskip
b) The knot passes through the ``infinity" point. In such a case the knot is an
infinite open curve. (Figure 2b)
\par \medskip
One should notice that finite open curves do not belong to the space of knots,
and that the two cases mentioned above, even though they ``look" quite
different in $R^3$, are topologically equivalent, since by a continuous
displacement of the infinity point one may go from one case to the other.
When projected on $S^2$, knots become closed curves, although these curves may
now possess multiple points.
\par \medskip
Two knots are defined {\bf ambient isotopic} and considered equivalent if
there is a continuous path in the space of knots that connects them.
In other words, if $f_1(s)$ and $f_2(s)$ are the homeomorphisms
defining a pair of knots, these knots are considered equivalent if and only if
there is a function $g(s,t)$ defined on $S^1 \times [0,1]$ continuous with
respect to both $s$ and $t$, such that $g(s,0)=f_1(s)$ and $g(s,1)=f_2(s)$ for
all $s \in S^1$, and $\forall t \in [0,1]: s_1 \neq s_2 \Rightarrow g(s_1,t)
\neq g(s_2,t)$. In addition, knots satisfying $f_1(s)=f_2(-s)$ (inverse
orientation) or $f_1(s)=-f_2(s)$ (mirror symmetry) shall be considered
equivalent whether they belong to the same continuous component or not, in
order to avoid the problems of {\bf chirality} and {\bf orientation reversal}
which were briefly discussed in Ref. [1].
\par \medskip
While knots are studied after being projected to two dimensions, not all
possible
two-dimensional projections are suitable for our purpose. The projections we
are going to use are called {\bf regular} and satisfy the following conditions.
\par \medskip
a) While projections are not necessarily one-to-one, regular projections may
map at most two points of a knot to the same point of $S^2$.
\par \medskip
b) Let $f(s)$ be the function from $S^1$ to $S^3$ which defines the knot, and
let $f(s_1)$ and $f(s_2)$ be projected on the same point of $S^2$. The
difference $f'(s_2)-f'(s_1)$ should not be orthogonal to
$S^2$ at $p[f(s_1)]=p[f(s_2)]$ in order for
the projection to be regular.
\par \medskip
c) Finally a regular projection may not possess points $s$ such that the
derivative $f'(s)$ is orthogonal to $S^2$ at $p[f(s)]$.
\par \medskip
In Figures 3a, 3b and 3c we show segments of projections where conditions a),
b) and c) are respectively violated, and thus the corresponding projections
cannot be considered regular. Once a
projection fulfills these conditions and is therefore regular, one has to
distinguish between overcrossings and undercrossings at each double (crossing)
point, so that Figures 4a and 4b should be considered distinct.
\par \medskip
Knot projections that can be continuously deformed to each other and are thus
topologically equivalent, do belong to ambient isotopic knots as one would
expect. In addition, two knots may be ambient isotopic even if their
projections are inequivalent, so long as these projections are connected to
each other through a series of {\bf moves} named after Reidemeister. There are
three kinds of such moves and are shown in Figure 5. One may
notice that performing these moves involves passing through irregular
projections; the first Reidemeister move involves a projection violating
condition c), a second move involves a projection violating b), while a third
move a projection violating a).
\par \medskip
Let a knot $K$ possess a regular projection $P(K)$.
According to Ref. [1], a regular projection will be denoted by
$\{ (a_1,a_2),(a_3,a_4),...,(a_{2n-1},a_{2n}) \} $ such that
$a_i \in \{ 1,2,...,2n \} $ $\forall \quad i \in \{ 1,2,...,2n \}$ and
$i \neq j \Rightarrow a_i \neq a_j$.
This notation implies the existence of a ``starting" segment $s$ of $P(K)$
connecting two
crossing points; if one numbers the crossings encountered after $s$
by successive natural numbers from $1$ until $2n$, at the $k^{\rm th}$ double
point
the overcrossing will be assigned the number $a_{2k-1}$ and the undercrossing
the number $a_{2k}$ . An example is illustrated in
Figure 6.
\par \medskip
It is interesting to compare this notation to the {\bf Gauss word} of the
regular projection (Ref. [11]). Let $\alpha _1, \alpha _2,..., \alpha _{2n}$ be
the
Gauss word of a projection denoted by $\{ a_1,a_2),...,(a_{2n-1},a_{2n}) \}$.
Then $\alpha _{a_1} + \alpha _{a_2} = ... = \alpha _{a_{2n-1}} + \alpha
_{a_{2n}}=0$.
\par \medskip
By naming knot projections
this way, one is faced with the following questions.
\par \medskip
a) Do identical names imply equivalent knots? Yes, for prime knots at least,
although a chirality
ambiguity remains. By considering mirror symmetric knots equivalent, and by
focusing exclusively on prime knots, one may give an affirmative answer to this
question.
\par \medskip
b) Can all knots be assigned at least one such name? Yes, so long as the knots
are {\bf tame} and thus possess regular projections with finite number of
double points. We are not going to discuss {\bf wild} knots in this paper.
\par \medskip
c) Do all possible combinations of pairs of natural numbers denote some knot
projection? No, most of them in fact do not. A necessary (but not sufficient)
condition is that odd numbers must be paired to even numbers.
The sufficient condition is that all pairs of loops formed by a
projection are either ``tangent", meaning that they have at least one segment
in common, or that they intersect at an even number of points, so that they do
not violate the {\bf Jordan Curve Theorem} (Ref. [12]).
\par \medskip
d) How can one know whether two names
belong to equivalent projections? Each projection
may be shown to possess at most $4n$ distinct names, depending on the
starting position and the orientation. If orientation is invariant, a pair
$(i,j)$ is replaced by $(i+k,j+k)$ for some fixed
$k$, while if orientation is reversed $(i,j)$ is replaced by $(k-i,k-j)$.
(In case $i+k$ and/or $j+k$ come out larger than $2n$, one must subtract $2n$,
while if $k-i$ and/or $k-j$ come out non-positive, one must add $2n$).
Therefore one may easily check whether two names
belong to equivalent or distinct projections.
\par \medskip
e) Assuming now that two names belong to distinct projections
how can one know whether the projected knots are ambient
isotopic? While one may find the effect of a Reidemeister on a knot projection,
this does not answer the question. If of course one arrives from projection
$P_1$ to $P_2$ by performing a series of Reidemeister moves, equivalence has
been established; a failure to do so does not imply inequivalence, since one
may not have performed enough moves; as there is no upper limit in the
number of moves required to connect equivalent projections, one may not prove
knot inequivalence through such a procedure.
\par \medskip
Failing to adequately address the last of these questions was the ``Achille's
heel" of Ref. [1]. In spite of this weakness, we had been able to obtain the
first 35 non-trivial
knots, which are all the distinct knots that have at least one projection with
not more than 8 double points. Although we
had checked our results with those in the literature, we had been unable
at that point to prove by our own that these knots are in fact distinct. In
addition, when we attempted to get knots with 9 double points we had been
unable
to remove all repetitions, while the running time required had been enormous;
for
knots with 10 double points the running time would have definitely been
prohibitive
and thus no such attempt was made.
In the following Sections we are going to show how these
difficulties were later overcome by appropriately modifying the algorithms
and by introducing new ones.
\bigskip
\centerline {3. \qquad KNOT EQUIVALENCE ALGORITHMS}
\bigskip
According to the discussion up to this point, the set of names of knot
projections can be
divided into equivalence classes, where two names are considered equivalent
iff the corresponding knots are ambient isotopic.
One may easily check
that such a relation is indeed an equivalence relation, meaning that it
satisfies $N \sim N$, $N_1 \sim N_2 \Rightarrow N_2 \sim N_1$ and
$N_1 \sim N_2 \wedge N_2 \sim N_3 \Rightarrow N_1 \sim N_3$ for any names of
knot projections $N_1$, $N_2$, $N_3$ and $N$.
\par \medskip
We are now going to define a total ordering
relation $\leq$, where by ``total ordering" is meant that all names $N$,
$N_1$, $N_2$,
$N_3$ obey $$N \leq N, \quad N_1 \leq N_2 \wedge N_2 \leq
N_1 \Rightarrow N_1 = N_2, \quad N_1 \leq N_2 \wedge N_2 \leq N_3 \Rightarrow
N_1 \leq N_3$$ and for any two names $N_1$ and $N_2$, $$N_1 \quad \leq \quad
N_2 \qquad \vee \qquad N_2 \quad
\leq \quad N_1$$ This relation is going to be termed ``name preferrence",
signifying that if two equivalent knots $K_1$ and $K_2$ are projected to
$P_1$ and $P_2$ which bear the names $N_1$ and $N_2$,
where $N_1 \leq N_2$, then the
knot class containing $K_1$ and $K_2$ would ``prefer" to be ``addressed" by
the name $N_1$ than by the name $N_2$. Once such relations have been defined,
our task is to obtain the names belonging to each equivalence
class, and for each class to choose the ``most preferred" name to represent it.
When this procedure is completed, one may easily use these names to draw or
tie knots belonging to each class, and thus knot classification is solved.
\par \medskip
The ``preferrence relation" is defined as follows. Let $N_1$ and $N_2$ be two
names consisting of $n_1$ and $n_2$ pairs respectively. If $n_1 \neq n_2$,
then the name consisting of the smallest number of pairs is the preferred one.
Therefore, among all possible names the most preferred one is the empty set
which corresponds to the {\bf trivial} knot.
\par \medskip
Let $N_1$ and $N_2$ consist of the same number of pairs, $n$. $N_1$
associates $1$ with some even number $2a_1$, while $N_2$ associates $1$ with
$2b_1$, where
$a_1,b_1 \in \{ 1,2,...,n \} $. If $a_1 < b_1$, $N_1$ is preferred to $N_2$,
while if $a_1 > b_1$ then $N_2$ is preferred to $N_1$. If $a_1 = b_1$, let
$3$ be paired to $2a_2$ in $N_1$ and to $2b_2$ in $N_2$. We similarly compare
$a_2$ to $b_2$, and unless they are equal, one chooses as preferred name the
one pairing $3$ to the smallest possible even number. If $a_2 = b_2$, we focus
our attention to $5$ which is paired to $2a_3$ in $N_1$ and $2b_3$ in $N_2$
and proceed in a similar manner. For the time being, one does not care about
the ordering inside each pair.
\par \medskip
There is of course a chance that two names $N_1$ and $N_2$ pair all odd numbers
to the same even numbers, and thus one cannot determine which of them is
preferred. In such a case we proceed as follows. First we check whether $1$ is
at the left or at the right inside its pairs in $N_1$ and $N_2$. If $1$
appears at the left in $N_1$ and at the right in $N_2$, then $N_1$ is
preferred to $N_2$. If $1$ appears at the left in both pairs or at the right
in both pairs, we now check $3$, and, if necessary, $5$, $7$, ... , $2n-1$. If
the preferrence question is still unresolved, this can only mean that
$N_1 = N_2$, in which case both $N_1 \leq N_2$ and $N_2 \leq N_1$.
\par \medskip
One easily notices that the ``preferrence" relation satisfies the necessary
conditions in order to be total ordering. Equally important, one may
notice that each set of names, finite or infinite, does always possess a
most preferred element, since the number of pairs in a name cannot be less than
$0$; this would not have been the case had we defined
``preferrence" the other way round, and in such a case we would not have been
able to choose the most preferred name as the representative of its class.
Using the theory discussed
up to this point we construct an algorithm as follows.
\par \medskip
We first consider the set of non-negative integers $N_0 = \{ 0,1,2,... \} $.
For each of its elements, $n$, we consider the set of permutations $S_n$
defined as one-to-one functions from $\{ 1,2,...,n \} $ to itself. The number
of such permutations is $n!$, since $1$ may be mapped to any
element of $\{ 1,2,...,n \} $, while there are $n-1$ choices for $2$, $n-2$
for $3$ and so on. We order these $n!$ permutations as follows. A permutation
defined by a function $f$ is ahead of one defined by $g$ if $\exists i : \big(
j < i \Rightarrow f(j)=g(j) \big) \wedge f(i)<g(i)$. Having thus ordered all
such
permutations we proceed with each permutation $f$, by defining a function $p$
from $\{ 1,2,...,2n \} $ to itself as follows.
$$f(i)=j \qquad \Rightarrow \qquad p(2i-1)=2j \quad \wedge \quad p(2j)=2i-1$$
One may easily check that all knot ``shadows", defined as projections
where no distinction is made between overcrossings and undercrossings,
can be determined through some function $p$ obtained from a permutation $f$ as
shown above. The inverse however is not true. In order for $p$ to determine a
distinct shadow, two conditions must be satisfied.
\par \medskip
First, $p$ must not lead to a shadow already obtained. As stated
earlier, a shadow may be named through a total of at most $4n$ different
functions $p$. We keep only the one derived from the permutation ordered
ahead of all others.
\par \medskip
Second, $p$ must lead to an actual knot shadow. As a counterexample, the
permutation $p(1)=4$, $p(2)=9$, $p(3)=6$, $p(4)=1$, $p(5)=8$, $p(6)=3$,
$p(7)=10$,
$p(8)=5$, $p(9)=2$, $p(10)=7$ does not (Figure 7). In order to check the
``drawability" of $p$ one first obtains the loops formed by segments of the
shadow.
\par \medskip
We begin by assigning all elements of $\{ 1,2,...,n \} $ an element of
$\{ 0,1,-1 \} $. There is a number of $3^n$ such functions and they signify
all loops that can be drawn on a shadow. If a number $i$ is mapped to $0$,
this means that the crossing point $\bigl(2i-1,p(2i-1)\bigr)$ is not a vertex
of the loop. If $i$ is mapped to $1$, once the loop arrives at $p(2i-1)$ it
``turns" and proceeds from $2i-1$ to $2i$, $2i+1$, ... , until it reaches the
next vertex; finally $i$ mapped
to $-1$ denotes proceeding along the opposite direction. An example of such a
loop is shown in Figure 8. In this example, $1$ is assigned $0$, since $(1,4)$
does not belong to the loop, $2$ is assigned $-1$, since once the loop reaches
$(3,8)$ it turns towards $2$ rather than $4$, while $3$, $4$, $5$ and $6$ are
assigned $1$.
\par \medskip
For each such pair of loops one first checks
whether they are ``tangent" to each other, that is whether they share one or
more segments. If so, one may proceed with the next pair, while if not,
one counts the number of intersection points. For any actual shadow this
number is even due to the fact that a non-intersecting loop in $S^2$ divides
$S^2$ into two disjoint pieces (Jordan Curve Theorem) and thus a second loop
not tangent to the first, must ``enter" the first as many times as it ``exits"
and so the number of intersecting points is even (Figure 9). If
there is one or more pairs of non-tangent loops whose intersection number is
odd, then the corresponding shadow is undrawable. In the counterexample
provided above, the loops $1-2-3-4$ and $5-6-7-8$ are non-tangent and
intersect at one point, the point $(3,6)$, and thus the shadow is undrawable.
\par \medskip
Once a shadow passes all these tests, one is ready to proceed with the full
projection. In order to define such a projection once the shadow is given, one
has to choose at each crossing which is the overcrossing and
which the undercrossing. It may thus appear that there are $2^n$ possible
choices. By simultaneously switching all crossings however, one gets the mirror
symmetric
knot which is considered equivalent, and thus two such projections are
considered identical. This duplicity is overcome by the convention that the
double point assigned the number $1$ is always an overcrossing; therefore
instead of $2^n$ choices there are only $2^{n-1}$. Once a particular choice
has been made, one keeps track of it as follows. Up to this point, a pair
consisting of $i$ and $j$
was denoted by $p(i)=j$ and $p(j)=i$. If $i$ is
chosen as the overcrossing and $j$ as the undercrossing and thus $(i,j) \in N$,
the choice is denoted by a function $s$ from $\{ 1,2,...,2n \} $ to
$\{ \pm 1, \pm 2,..., \pm {2n} \} $ such that $s(i)=j$ and $s(j)=-i$.
\par \medskip
These $2^{n-1}$ choices are ordered as follows. We begin with the choice
where all odd numbers are overcrossings and all even are undercrossings. Such
a projection is called {\bf alternating}. This choice is denoted by a function
$b$ from $\{ 1,2,...,n \} $ to $\{ 0,1 \} $ such that $b(i)=0 \quad \forall
\quad i$. We
then proceed with a choice where
$1,3,5,...,2n-3$ are overcrossings and $2n-1$ is an undercrossing. Now
$b(1)=b(2)=...=b(n-1)=0$ and $b(n)=1$. We continue in a way similar to a
sequence of the digits of a binary number, where $b(i)=0$ denotes that
$s(2i-1)>0$ and $b(i)=1$ that $s(2i-1)<0$, until we end with
$b(1)=0$ and $b(2)=b(3)=...=b(n)=1$.
\par \medskip
Once a particular (full) projection is under consideration, the first thing
to check is whether this projection has been previously encountered under a
different name. While the possibility of a different permutation has already
been excluded, it may happen that symmetries of a shadow may lead to multiple
names of a projection; while keeping $f$ and $p$ the same, different values
of $b$ and $s$ may lead to the same projection. One such example is shown in
Figure 10. Here $p(1)=4$, $p(2)=5$, $p(3)=6$, $p(4)=1$, $p(5)=2$, $p(6)=3$. The
choices $b(1)=0$, $b(2)=0$, $b(3)=1$ and $b(1)=0$, $b(2)=1$, $b(3)=0$ lead to
identical projections.
\par \medskip
Once such cases have been excluded, one arrives to point where all
projections have been obtained and classified. Since the goal however is to
classify knots and not their projections, one still has to identify distinct
projections that belong to equivalent knots. This is checked by performing
Reidemeister moves. The whole procedure goes as follows.
\par \medskip
We begin with a particular projection and attempt to find moves that replace
it with a more ``preferrable" one according to the definition of this Section.
When dealing with first or second moves, one need only examine the ones
reducing the number of double points; third moves do not affect this number,
and thus one has to examine more carefully whether the new name is preferred.
If no such name can be found through any of the three kinds of moves, one
assigns a successive natural number to this projection. Conventionally we
assigned $0$ to the unknot, and went on with $1$ for the trefoil, $2$ for
{\bf Listing's knot} and so on (Figure 11).
\par \medskip
A name that can be directly
reduced to a preferred one, clearly cannot represent its class. On the other
hand, names that have been assigned these numbers and thus cannot be directly
reduced, are not guaranteed to represent their classes, since it may be
possible to find equivalences between these names. This is checked
as follows.
\par \medskip
Let $N$ be a name denied a number according to the procedure above because a
Reidemeister move
can reduce it to the preferred names $N_1$, $N_2$, ..., $N_f$. These
names have been assigned the numbers $n_1$, $n_2$, ..., $n_f$ either through
the procedure of the previous paragraph or retroactively through the procedure
we are now presenting. One assigns $N$ the number $$n_{min}=\min_{1\leq i
\leq f} n_i$$ Subsequently any name previously assigned a number
$n_i \in \{ n_1,n_2,...,n_f \} - \{ n_{min} \}$ is now assigned $n_{min}$. In
addition, we ``mark" all numbers $n_i \neq n_{min}$ in order to denote that
the names assigned to them have been found equivalent to a preferrable one. A
name may thus be initially be assigned some number $M$, but later this number
may change to some $m<M$.
\par \medskip
When the algorithm has ended running, we only keep the first name that has
been assigned to ``unmarked" numbers, since it is only these names that have
not been reduced directly or indirectly. Clearly, the further we run the
algorithm the better the results since there is a better chance to identify
possible equivalences; no matter how far the program is run, however, this
algorithm by itself cannot guarantee the absence of equivalencies; this
problem will be discussed from Section 5, after describing the practical
problems faced in running the algorithm and discussing ways to overcome them.
\bigskip
\centerline {4. \qquad COMPUTER MEMORY - RUNNING TIME}
\bigskip
While the theory described in Section 3 can lead to an algorithm correct in
principle, when one uses it to develop a computer program, two practical
difficulties arise.
\par \medskip
First, running a computer program requires a minimum memory, which means
a storage capacity where information is kept to be used later. When, for
example, one studies some shadow defined by a function $f$ taking values on
$\{ 1,2,...,n \} $, the computer needs $n$ memory units to store the values of
$f(1),f(2),...,f(n)$. In order for this to be achieved, one must request in
advance $n$ at least units for the function $f$.
\par \medskip
By looking at the algorithm suggested one notices that $f$ is just one of the
many functions to be used, each of them needing some space in the computer's
memory. If one neglects to request in advance the necessary memory, the
program will either stop running at some point if the account is VAX/VMS, or
even worse, come up with wrong results in a UNIX account. On the other hand, if
the memory requested exceeds the account's total capacity, the program cannot
be compiled or is not going to run. It is therefore essential to know or at
least estimate in advance the space each function requires, and minimize this
space as much as possible.
\par \medskip
Actually one-variable functions such as $f$ mentioned above,
require very limited space, and so long as one is ``generous" enough they do
not cause any problem. The real concern are the two-variable functions
$F(i,j)$, since the space they require is equal to the product of the space
required for each of the two variables.
\par \medskip
Fortunately most two-variable functions require at most $6n$ memory units
per variable, which gives a total of $36n^2$. For $n=10$ this gives $3600$
units which is insignificant compared to the total capacity which is equal to a
few million units. There are however two exceptions which require much more
memory units.
\par \medskip
The first case is related to functions recording the loops on a potential
shadow, necessary to check whether the shadow can be drawn. One of their two
variables takes values in $\{ 1,2,...,2n \} $ indicating whether a particular
double point is a vertex, belongs to a side or does not belong to the loop.
The other variable indicates the order of the loop, and can take values from
$1$ until $3^n$, since each
crossing can be mapped to $0$ or $\pm 1$.
For $n=10$ this is equal to $3^{10}=59049$ which when
multiplied by $2n=20$ exceeds one million. By avoiding self-intersecting loops
one can substantially reduce this number; instead of $3^n$ one gets a
number of about $2^n$. For $n=10$ the actual number of not self-intersecting
loops came out equal to $781$, and when this was multiplied by $2n=20$ the
total memory requested was merely $15620$ units.
\par \medskip
The second case is related to assigning numbers to the various knot
projections and is even more demanding. For $n=10$ for instance there is a
total of $26532$ different shadows, each giving rise to $2^{10-1}=512$
projections, which gives a total of more than thirteen million units. This
quantity is just above or below the limit of most accounts; even if one
succeeds in accomodating such a demand, one is completely hopeless
when running $n=11$.
\par \medskip
Therefore one needs to modify the algorithm presented, not because there is
some logical flaw in it, but in order to reduce the memory required. Since all
names appearing in the output are irreducible and have been
assigned successive integer numbers, we need only keep information on such
numbers and the corresponding names. Instead of performing just one
Reidemeister move on a reducible name $N$ and obtain the equivalent names
$N_1$, $N_2$, ..., $N_f$ which are not necessarily irreducible, we continue
reducing each $N_i$ until we reach the irreducible names ${\cal N}_1$,
${\cal N}_2$, ..., ${\cal N}_l$, where $l \geq f$. These names have been
assigned the numbers $n_1$, $n_2$, ..., $n_l$. We now ``identify" these numbers
which belong to equivalent knots by ``projecting" them to the smallest of
them, $n_{min}$.
\par \medskip
It is possible of course that some of them have already been
``projected" to some $m<n_{min}$. Due to this reason, we achieve better
results if instead of projecting $n_1$, ..., $n_l$ to $n_{min}$, we project
them instead to $m_{min}$, where $m_1$, ..., $m_l$ were the previous
projections of $n_i$. Subsequently, we project to $m_{min}$ any number that
had previously been projected to some other $m_i$.
\par \medskip
By doing so, after each knot projection has been studied, we only need keep
information regarding the irreducible knot names, their assigned numbers and
their projections. In the end one keeps only the irreducible names whose
assigned numbers are
projected to themselves. For $n=10$ the largest assigned number came out
$3977$. Since each such number requires $n$ memory units to record the $n$
pairs of the name, the memory required now is less than $40,000$ contrasted
to $13,000,000$ before. The memory problem has thus been solved up to
at least $n=13$.
\par \medskip
Having sufficient memory however is just one of the practical difficulties
encountered. If one could afford to wait long enough, a flawless algorithm
with sufficient memory would eventually produce an output. Unfortunately the
time demanded increases very fast, and unless one is extremely careful, it
may even exceed an average lifespan. One may reasonably estimate that for our
problem the time should be proportionate to $2^{n-1} \times n!$, so that each
input increase by $1$ will multiply the time needed by $2n$, and thus running
$n=11$ needs twenty times the time needed to run $n=10$. While using more
efficient computer machines (for example UNIX is much faster than VAX/VMS)
and optimizing can help up to a point, it is essential that one checks
carefully every step of the algorithm to see if it is really necessary and
whether it can be done faster.
\par \medskip
Comparing the two versions of the algorithm described before, one notices
that the one demanding less memory needs more time to run. This
is due to the fact that the first version requires performing just one move on
each knot projection,while the second requires to continue performing moves
until one reaches irreducible names. In addition, once an irreducible name
is reached, finding its assigned number can be very time consuming, since one
may have to check all natural numbers (up to $3977$ for $n=10$) to see which
one corresponds to the irreducible name. This occurs, because while before
the knot projection was mapped to the number, and thus the high memory
required, now the number is mapped to the projection. For each of these numbers
one also may  need to perform up to $n$ checks to compare the pairs of the knot
projection of the numbers with the projection whose number we want to find.
\par \medskip
The last difficulty can be overcome as follows. Instead of asking whether
the name to be checked is equal or not to the names defined through the
numbers 1,2,3,...,3977, one may ask whether this name is equal, more
preferrable or less preferrable to a ``middle" name such as the one defined for
instance by
the number 2000. If it is not equal, we have limited our scope among 2000
names, so we now compare the name to be checked to the names defined by 1000
or 3000, and thus limiting our scope by another factor of 2. By repeating this
process we finally obtain the number desired in at most 12 steps (since
$2^{12}>3977$) instead of an average of about 2000 steps. While the steps now
are slightly slower, since the question asked has three possible answers
instead of two, the advantage obtained by reducing the number of steps is
significant.
\par \medskip
Having overcome this hurdle, we now focus on the previous one. This can be
overcome by noticing that in the process described many of the equivalences
obtained are unnecessary. If for example one proves that knot class $A$ is
equivalent to knot class $B$ and $B$ is equivalent to $C$, one should not
waste any time performing a subroutine that will subsequently prove $A$ to be
equivalent to $C$. In general, showing that $n$ quantities are equivalent may
require
at least $n-1$ and at most $n^2$ comparisons depending on the algorithm's
efficiency. We had better make sure that this number is as close to $n-1$ as
possible.
\par \medskip
We now discuss how this principle works in practice. Let $N$ be a knot
name, which after one Reidemeister move can be reduced to $N_1$, $N_2$, ...,
$N_l$. In general $N_1$ can be similarly reduced to $N_{11}$, $N_{12}$, ...,
$N_{1f_1}$, while $N_2$ can be reduced to $N_{21}$, $N_{22}$, ..., $N_{2f_2}$
and so on, until we get a total of $f_1+f_2+...+f_l$ new names, and so on. One
may think that by checking if any names are equal to each other and avoiding
repetitions some time may be saved. That is indeed the case, although checking
for such equalities consumes additional time and thus the final saving is
marginal.
\par \medskip
One can save much more however by realising that the names $N_{1i}$, where
$1 \leq i \leq
f_1$, have already been shown equivalent to each other, and therefore there is
no reason to keep all $f_1$ names. We need only keep $N_{11}$ and go on with
$N_2$. Similarly for $N_2$ we only have to keep $N_{21}$ and so on. Once this
step is completed, we have been left with at most $l$ names $N_{11}$, $N_{21}$,
..., $N_{l1}$ instead of $l_1+l_2+...+l_f$, since the information we could
have got from the other names is redundant. If in addition, we are left at some
point with just one name, there is no reason to go on further, since no
additional equivalence can be obtained.
\par \medskip
One may further accelerate the program by noticing that first and second
Reidemeister moves ``commute" with each other in the sense that if a number of
such moves is possible, the total output is independent of the order with which
these moves are performed. On the contrary, third moves do not ``commute"
either with themselves or with first and second moves; performing a third
move may even destroy the ability to perform some other move (Figure 12).
Therefore if a name $N$ can be reduced to names $N_{\alpha _1}$, ...,
$N_{\alpha _i}$ through first moves, $N_{\beta _1}$, ..., $N_{\beta _j}$
through second moves and $N_{\gamma _1}$, ..., $N_{\gamma _k}$ through third
moves, while one must keep all $N_{\gamma }$ ones, one need only keep just one
of the $N_{\alpha }$ and $N_{\beta }$ names, preferrably the first one
encounters, and avoid checking for the existence of other such names.
\par \medskip
While however third moves may not ``commute", one notices the following.
Let $N$ be some name that can accept a first and a third move which do not
commute with each other. Let $N_\alpha $ be the name obtained after the first
move and $N_\gamma $ the name after the third move. One may see that
$N_\alpha $ is preferred to $N$, while $N$ is preferred
to $N_\gamma $ (Figure 13); therefore any information
provided by $N$ will be later provided again by $N_\gamma $. Therefore one
may ignore all names that can be reduced by first moves, which are the ones
whose shadows satisfy $f(1)=1$. This idea alone reduced the running time by a
factor of $15$.
\par \medskip
Finally one may save further by avoiding ``connected sums" (non-prime knots),
such as the knot shown in Figure 14.
These are denoted by names $N$ for which there are integers $k$ and $l$ other
than $1$ and $2n$ such that $(i,j) \in N, \qquad k \leq i \leq l
\Leftrightarrow k \leq j \leq l$. One does need to keep such names if they
contain $n-2$ pairs or fewer, since they may be equivalent to names with $n$
pairs or fewer that do not satisfy this property (Figure 15). By eliminating
however connected sums with $n-1$ and $n$ pairs one does reduce
substantially the running time.
\par \medskip
This is not the only reason for which connected sums are ignored. As stated
earlier, connected sums possess a name ambiguity, since it is possible for
two such knots to have the same name but be distinct. In addition, as we
shall see later, connected sums are not uniquely determined by their knot
groups. Most important, however, once prime knots are obtained, one may easily
get the connected sums by combining prime knots.
\par \medskip
Once these ideas are implemented, the algorithm is not merely faster than
before, but also faster than the one which required higher memory. While for
example $n=10$, the maximum value for which we were able to run the high
memory algorithm, required about six hours, now it can run in just over four
minutes. In addition, we were now able to run $n=11$ in about an hour and
twenty
minutes, and $n=12$ in about two days. The results we got can be
summarised as follows.
\par \medskip
\centerline {\vbox {\def\*{\hphantom{0}} \offinterlineskip
\halign {\strut#&\vrule#\quad&\hfil#\hfil&\quad\vrule#\quad&\hfil#\hfil&\quad
\vrule#\cr
\noalign {\hrule }
&&\omit \bf Crossing Number && \bf Maximum Number of Knots &\cr
\noalign {\hrule }
&&\*0&& \*\*\*\*1&\cr &&\*1&& \*\*\*\*0&\cr &&\*2&& \*\*\*\*0&\cr
&&\*3&& \*\*\*\*1&\cr &&\*4&& \*\*\*\*1&\cr &&\*5&& \*\*\*\*2&\cr
&&\*6&& \*\*\*\*3&\cr &&\*7&& \*\*\*\*7&\cr &&\*8&& \*\*\*21&\cr
&&\*9&& \*\*\*49&\cr &&10&& \*\*165&\cr &&11&& \*\*945&\cr
&&12&& 44322&\cr \noalign {\hrule}}}}
\bigskip
\centerline {5. \qquad SHOWING KNOT INEQUIVALENCE}
\bigskip
Even if one was able to solve all practical problems presented in the previous
Section, the presence of the input value of $n$ which acts as a cutoff
parameter, implies the possibility of repeated knot classes, since their
equivalence might only be shown for a higher value of $n$. Since $n$ cannot be
made infinite, there is no guarantee yet
that any two names appearing at the output might not be eventually shown
equivalent. It is therefore essential to establish some procedure which can
directly show knot inequivalences.
\par \medskip
Historically the first such procedure was presented during the 1920's by
Reidemeister and thus it was only then that the existence of non-trivial knots
was rigorously proved (Ref. [13]). His procedure was to subject knot
projections to a
suitable test, whose results were invariant under equivalence moves. Therefore
if some projection yielded a particular response, all equivalent projections
would yield identical responses; through contradiction one could thus show
that if two projections respond differently, they definitely belong to
inequivalent knots. The test he proposed for that purpose is an attempt to
``paint" a knot projections by using a choice of three colors and obeying the
following rules.
\par \medskip
1) All points of a particular ``strand", defined as the segment joining
consecutive undercrossings, must be colored the same.
\par \medskip
2) At each crossing, the three strands that meet each other must either be
colored the same, or three distinct colors must be used; this is equivalent to
disallowing for two strands to have one color and the third strand to have
another.
\par \medskip
Under these rules one can always paint any knot projection by using a single
color through the entire projection. But Reidemeister also imposed that all
three colors must be used in at least one strand. In Figures 16 a)-c) we show
that this test makes sense, meaning that when an equivalence move is performed
on some projection that fails the test, and thus covered entirely by the same
color $c$, it is impossible for the new projection
to succeed. In Figure 17 we show that the test is successful for the trefoil;
for the unknot, which can be represented by a circle, it is obviously not, and
therefore the trefoil is non-trivial.
\par \medskip
The obvious deficiency of this test is that it allows for only a ``YES" or
``NO" response, and thus can only divide knot projections into two classes.
Knots however responding affirmatively may still be distinct, and so can knots
responding negatively. Since the 1920's a number of far better tests have been
introduced, usually in the form of a ``polynomial", starting from the
{\bf Alexander} polynomial, already introduced in the 1920's (Ref. [14]), and
continuing until the {\bf HOMFLY} and the {\bf Akutsu Wadati}
polynomials which were introduced during the last ten years (Ref. [15]).
These polynomials are one or two variable functions which are defined on knot
and link projections through the {\bf skein relations}
$aP(L_+)+bP(L_-)+bP(L_0)=0$ (Figure 18). Care is taken so that the polynomials,
although originally defined on projections, are invariant under Reidemeister
moves (with the exception of the {\bf Kauffmann} polynomial that is not
invariant under the first move in order to distinguish among twisted bands,
Ref. [16]). If two projections are assigned distinct
polynomials, they are definitely inequivalent.
\par \medskip
One may wonder however how good is the definition based on the skein relations.
If $L_+$ possesses $n$ double points, $L_0$ possesses $n-1$ double points but
$L_-$ possesses $n$. It looks as if polynomials of $n$ crossing links can only
be obtained if they are already known. This problem is easily overcome by
noticing that a knot possessing a name $N$ such that $(i,j) \in N \Rightarrow
i < j$, where overcrossings always proceed undercrossings (or vice versa), is
trivial. Therefore through suitable crossing switches one arrives at a
projection which possesses $n$ double points but is trivial.
Subsequently, successive applications of the skein relation will lead to the
polynomial of the unknot, which is defined equal to $1$, and to polynomials
of disjoint links
defined $0$.
\par \medskip
In spite of the obvious advantages, these polynomials suffer from a few
setbacks. First, they are not one-to-one, and thus unable to fully
distinguish between inequivalent knots. In other words, inequivalent knots
may possess identical polynomials; one such example are the so called {\bf
mutant} knots (Ref. [17]). A more famous example where distinct knots have
identical Jones polynomials, had previously been shown by Birman in Ref. [18].
\par \medskip
Second, even if one were to tolerate a few such cases, where neither
equivalence
nor inequivalence could be demonstrated, the practical problems encountered
when trying to calculate these polynomials are substantial. First, one cannot
find the polynomial of a knot without considering multicomponent links as well.
In the skein relation, if $L_+$ is a knot, $L_-$ is also a knot, but $L_0$ is
a two-component link. In general, if $L_+$ is an $m$-component link, so is
$L_-$, but $L_0$ consists of $m+1$ components. Therefore, when calculating
retroactively the polynomial of some knot, one would encounter links with as
many as $n \over 3$ components. The knot notation thus introduced would not be
sufficient and one would have to generalise it, and simultaneously generalise
the equivalence algorithm so that it could identify equivalent links in
addition to equivalent knots.
\par \medskip
In addition, the memory required would be enormous; in order to calculate the
polynomial of an $n$ double point knot or link, the computer would have to
store the polynomials of all knot and link projections of up to $n-1$ double
points. Even worse, a $k$ degree polynomial demands not just $1$ but $k$
memory units to be stored. The running times would also be expected to become
prohibitive for fairly small values of $n$.
\par \medskip
One may expect to overcome these problems through the {\bf knot groups},
defined as the fundamental (first homotopic) groups of the knot complements.
This can also be understood as the group of loops that do not intersect the
knot. Knot groups possess the significant advantage that, when at least
only {\bf prime} knots are considered, two
inequivalent knots do always possess distinct groups (Ref. [19]).
In addition, finding the group of a knot does not involve any skein relations;
once a projection is given, one can easily obtain the {\bf Wirtinger}
presentation as follows (Ref. [20]). Every strand of the projection is mapped
to a generator $g_i$, while at every crossing, the generators of the meeting
strands satisfy either the relation $g_{i+1}=g_jg_ig_j^{-1}$ if the crossing
looks like Figure 19a, or $g_{i+1}=g_j^{-1}g_ig_j$ if it looks like Figure 19b.
\par \medskip
While obtaining presentations of knot groups is therefore trivial, the
presentations of two distinct projections are different, even if the
projections correspond to equivalent knots. Directly comparing two
presentations is just as hard as comparing the projections. On the other hand,
one may clearly understand now how Reidemeister's three color test works.
When one attempts to ``paint" a knot projection according to the rules of the
test, one essentially asks the following question.
\par \medskip
Let $G$ be a group generated by three generators $a$, $b$ and $c$ which
satisfy the relations $bab^{-1}=c$, $cac^{-1}=b$, $aba^{-1}=c$, $aca^{-1}=b$,
$cac^{-1}=b$, $cbc{-1}=a$. Is $G$ a {\bf coset} of the knot group $K$, or
equivalently, is there a group $H$ which is a subgroup of $K$ such that
$G=K/H$? If
there is such an $H$, all elements of $K$ can take the form $k=gh$ where
$k \in K$, $g \in G$ and $h \in H$, and are ``projected" to $g$. Each generator
of $K$ is projected to a generator of $G$, and this process is exactly what is
happenning when one ``paints" the knot. Clearly, if $G$ is a coset of a knot
group $K_1$ but not of some other knot group $K_2$, the groups are different
and thus the knots are inequivalent.
\par \medskip
Two questions arise at this point. First, if $G$ is actually a
group. Clearly, one cannot form a group by merely throwing some generators
with a set of constraints, without making sure that the constraints are
consistent. In the example above, one may notice that $G=P_3 \times Z$, where
$P_3$ is the $3$-permutation group. $G$ differs from $P_3$ since in $P_3$ one
must have $a^2=b^2=c^2=e$, while the constraints of $G$ merely imply that
$a^2=b^2=c^2$ and that they commute with all elements of $G$.
Second, is it possible to find additional ``color tests"? Now that we found a
method to distinguish between inequivalent knots, we would definitely want to
do more than merely separate knots into just two classes. In order to achieve
this goal, we first tried to generate a number of ``color tests" and then apply
them to knot projections in order to check their ``response".
\par \medskip
As mentioned above, finding suitable ``color tests" is not trivial; one
cannot just allow for a number of ``colors" $a_1$, $a_2$, ..., $a_r$ and
impose arbitrarily that at each crossing only some randomly selected triades
are permitted. While checking that the generators together with the constraints
do form a group may seem hard, that is not necessarily what is needed.
What we actually desire is a painting procedure such that if $K_1$ and $K_2$
are knot projections related through Reidemeister moves, for each permissible
``colorisation" of $K_1$ there is exactly one permissible ``colorisation" of
$K_2$. Once this is achieved one will be able to ask ``In how many distinct
ways
is it possible for a knot to be painted following these rules?", instead of
merely checking whether the painting is possible. Therefore the test one thus
obtains is even finer. This can be achieved as follows.
\par \medskip
Let $a_1$, $a_2$, ..., $a_r$ be the ``colors" to be used, and let
$a_k=f(a_i,a_j)$ imply that at each junction where $a_i$ meets $a_j$ as shown
in Figure 20, the color that will emerge is $a_k$. The invariance of the test
under Reidemeister moves is checked as follows.
\par \medskip
First move (Figure 21a): The strands before and after the crossing affected by
a first move must be colored the same. Therefore one requires $f(a_i,a_i)=a_i$.
This is equivalent to imposing $a_i^{\pm \epsilon}a_ia_i^{\mp \epsilon}=a_i$,
where $\epsilon \in \{ +1,-1 \} $.
\par \medskip
Second move (Figure 21b): The strands before and after the two crossings
affected by a second move must be colored the same, while the color of the
strand between them must be uniquely determined. Therefore one requires
$f(a_i,a_j)=f(a_k,a_j) \Leftrightarrow a_i=a_k$. This is equivalent to
imposing $a_j^{\pm \epsilon}a_ia_j^{\mp \epsilon} =
a_j^{\pm \epsilon}a_ka_j^{\mp \epsilon} \Leftrightarrow a_i=a_k$, where once
more $\epsilon \in \{ +1,-1 \} $.
\par \medskip
Third move (Figure 21c): While the colors of the top three strands
$a_i$, $a_j$ and $a_k$, are not
going to change, for the three strands belonging to the bottom segment one
requires that the colors at the two sides remain invariant too. The color at
the middle may change, but should be uniquely determined. This condition may be
written as follows. $$a_k=f(a_i,a_j) \vee a_m=f(a_l,a_j) \vee a_n=f(a_l,a_i)
\Rightarrow f(a_m,a_k)=f(a_n,a_j)$$ It implies that
$a_ja_ia_j^{-1}a_ja_la_j^{-1}a_ja_i^{-1}a_j^{-1}=a_ja_ia_la_i^{-1}a_j^{-1}$.
\par \medskip
While these conditions are necessary but not always sufficient to form a
group, they are sufficient to form a color test good enough for our purpose.
The first two conditions imply that all $n$ columns of the matrix $M_{ij}=
f(a_i,a_j)$ are $n$-permutations, under the constraint that the $i^{\rm th}$
permutation satisfies $p(i)=i$. Therefore the columns are effectively $n-1$
permutations. Consequently, one may start with $(n-1)!^n$ instead of $n^{n^2}$
possibilities, which is an enormous saving; for $n=4$ for example, $(n-1)!^n=
6^4=1296$, while $n^{n^2}=4^{16}=2^{32}>4 \times 10^9$. In addition, once the
remaining possibilities have been eliminated, one need only check whether the
third condition is satisfied.
\par \medskip
Even now however the possibilities increase enormously with $n$.
For $n=5$ the number of combinations
becomes $24^5>8 \times 10^6$, while for $n=6$ this number increases to
$120^6$ which is about $3 \times 10^{12}$. While the memory required for such a
program is insignificant, the running time increases very rapidly and thus
time-saving ideas are required.
\par \medskip
A first such idea takes advantage of the fact that each $n$ color test may be
written in up to $n!$ different ways by permuting the $n$ colors. If for
example $M_{ij}$ defines an acceptable test, so does
$p^{-1} \big( M_{p(i),p(j)} \big) $,
where $p$ is an $n$-permutation.
\par \medskip
In addition, once a particular combination fails, one may first want to know
the ``reason" of that failure instead of blindly proceeding to the next
combination. This means that one should find out the values of $M_{ij}$,
$M_{lj}$, $M_{la}$, $M_{mk}$ and $M_{nj}$ which do not satisfy the third
condition. Unless these are changed, there is no reason to test any other
combinations.
In particular, the third condition implies $M_{ij}=i \Leftrightarrow M_{ji}=j$.
By avoiding all possibilities violating this rule one may accelerate the
program enormously.
\par \medskip
Finally, one need not keep tests that can reduced into smaller ones. This
means that there should not exist any set $S \subset \{ 1,2,...,n \} $ other
than the empty set such that $i \in S \Rightarrow
M_{ij} \in S \quad \forall \quad j \in \{ 1,2,...,n \} $.
\par \medskip
Even after applying all the ideas above we were merely able to obtain color
tests up to eight colors. These are the following.
\par \vskip 1 cm
\def\*{\hphantom{0}}
\centerline {$\underline {THREE \quad COLORS}$}
\vskip 0.5 cm
\centerline {Test 1}
\vskip 0.5 cm \centerline
{\vbox
{\offinterlineskip
\halign{\strut#&\vrule#\quad&\hfil#\hfil&\quad
\vrule#\quad&\hfil#\hfil&\quad#&\hfil#\hfil&\quad#&\hfil#\hfil&\quad\vrule#\cr
\noalign{\hrule}
&&\omit \quad \hfil && \quad 1 \quad && \quad 2 \quad && \quad 3 \quad
& \cr \noalign {\hrule}
&&1&& 1&& 3&& 2&\cr
&&2&& 3&& 2&& 1&\cr
&&3&& 2&& 1&& 3&\cr
\noalign{\hrule}}}}
\vskip 1 cm
\centerline {$\underline {FOUR \quad COLORS}$}
\vskip 0.5 cm
\centerline {Test 2}
\vskip 0.5 cm \centerline
{\vbox
{\offinterlineskip
\halign{\strut#&\vrule#\quad&\hfil#\hfil&\quad
\vrule#\quad&\hfil#\hfil&\quad#&\hfil#\hfil&\quad#&\hfil#\hfil&\quad
#&\hfil#\hfil&\quad\vrule#\cr
\noalign{\hrule}
&&\omit \quad \hfil && \quad 1 \quad && \quad 2 \quad && \quad 3 \quad
&& \quad 4 \quad & \cr \noalign {\hrule}
&&1&& 1&& 4&& 2&& 3&\cr
&&2&& 3&& 2&& 4&& 1&\cr
&&3&& 4&& 1&& 3&& 2&\cr
&&4&& 2&& 3&& 1&& 4&\cr
\noalign{\hrule}}}}
\vfil \eject
\centerline {$\underline {FIVE \quad COLORS}$}
\vskip 0.5 cm
\centerline {Test 3}
\vskip 0.3 cm \centerline
{\vbox
{\offinterlineskip
\halign{\strut#&\vrule#\quad&\hfil#\hfil&\quad
\vrule#\quad&\hfil#\hfil&\quad#&\hfil#\hfil&\quad#&\hfil#\hfil&\quad
#&\hfil#\hfil&\quad#&\hfil#\hfil&\quad\vrule#\cr
\noalign{\hrule}
&&\omit \quad \hfil && \qquad 1 \qquad && \qquad 2 \qquad && \qquad 3 \qquad
&& \qquad 4 \qquad && \qquad 5 \qquad& \cr \noalign {\hrule}
&&1&& 1&& 4&& 5&& 3&& 2&\cr
&&2&& 3&& 2&& 4&& 5&& 1&\cr
&&3&& 2&& 5&& 3&& 1&& 4&\cr
&&4&& 5&& 1&& 2&& 4&& 3&\cr
&&5&& 4&& 3&& 1&& 2&& 5&\cr
\noalign{\hrule}}}}
\vskip 0.5 cm
\centerline {Test 4}
\vskip 0.3 cm \centerline
{\vbox
{\offinterlineskip
\halign{\strut#&\vrule#\quad&\hfil#\hfil&\quad
\vrule#\quad&\hfil#\hfil&\quad#&\hfil#\hfil&\quad#&\hfil#\hfil&\quad
#&\hfil#\hfil&\quad#&\hfil#\hfil&\quad\vrule#\cr
\noalign{\hrule}
&&\omit \quad \hfil && \qquad 1 \qquad && \qquad 2 \qquad && \qquad 3 \qquad
&& \qquad 4 \qquad && \qquad 5 \qquad& \cr \noalign {\hrule}
&&1&& 1&& 3&& 4&& 5&& 2&\cr
&&2&& 3&& 2&& 5&& 1&& 4&\cr
&&3&& 4&& 5&& 3&& 2&& 1&\cr
&&4&& 5&& 1&& 2&& 4&& 3&\cr
&&5&& 2&& 4&& 1&& 3&& 5&\cr
\noalign{\hrule}}}}
\vskip 1 cm
\centerline {$\underline {SIX \quad COLORS}$}
\vskip 0.5 cm
\centerline {Test 5}
\vskip 0.3 cm \centerline
{\vbox
{\offinterlineskip
\halign{\strut#&\vrule#\quad&\hfil#\hfil&\quad
\vrule#\quad&\hfil#\hfil&\quad#&\hfil#\hfil&\quad#&\hfil#\hfil&\quad
#&\hfil#\hfil&\quad#&\hfil#\hfil&\quad#&\hfil#\hfil&\quad\vrule#\cr
\noalign{\hrule}
&&\omit \quad \hfil && \qquad 1 \qquad && \qquad 2 \qquad && \qquad 3 \qquad
&& \qquad 4 \qquad && \qquad 5 \qquad && \qquad 6 \qquad& \cr \noalign {\hrule}
&&1&& 1&& 1&& 4&& 3&& 6&& 5&\cr
&&2&& 2&& 2&& 5&& 6&& 3&& 4&\cr
&&3&& 4&& 5&& 3&& 1&& 2&& 3&\cr
&&4&& 3&& 6&& 1&& 4&& 4&& 2&\cr
&&5&& 6&& 3&& 2&& 5&& 5&& 1&\cr
&&6&& 5&& 4&& 6&& 2&& 1&& 6&\cr
\noalign{\hrule}}}}
\vskip 0.5 cm
\centerline {Test 6}
\vskip 0.3 cm \centerline
{\vbox
{\offinterlineskip
\halign{\strut#&\vrule#\quad&\hfil#\hfil&\quad
\vrule#\quad&\hfil#\hfil&\quad#&\hfil#\hfil&\quad#&\hfil#\hfil&\quad
#&\hfil#\hfil&\quad#&\hfil#\hfil&\quad#&\hfil#\hfil&\quad\vrule#\cr
\noalign{\hrule}
&&\omit \quad \hfil && \qquad 1 \qquad && \qquad 2 \qquad && \qquad 3 \qquad
&& \qquad 4 \qquad && \qquad 5 \qquad && \qquad 6 \qquad& \cr \noalign {\hrule}
&&1&& 1&& 1&& 6&& 3&& 4&& 5&\cr
&&2&& 2&& 2&& 4&& 5&& 6&& 3&\cr
&&3&& 4&& 6&& 3&& 2&& 3&& 1&\cr
&&4&& 5&& 3&& 1&& 4&& 2&& 4&\cr
&&5&& 6&& 4&& 5&& 1&& 5&& 2&\cr
&&6&& 3&& 5&& 2&& 6&& 1&& 6&\cr
\noalign{\hrule}}}}
\vskip 13 cm
\centerline {$\underline {SEVEN \quad COLORS}$}
\vskip 1 cm
\centerline {Test 7}
\vskip 0.5 cm \centerline
{\vbox
{\offinterlineskip
\halign{\strut#&\vrule#\quad&\hfil#\hfil&\quad
\vrule#\quad&\hfil#\hfil&\quad#&\hfil#\hfil&\quad#&\hfil#\hfil&\quad
#&\hfil#\hfil&\quad#&\hfil#\hfil&\quad#&\hfil#\hfil&\quad
#&\hfil#\hfil&\quad\vrule#\cr
\noalign{\hrule}
&&\omit \quad \hfil && \qquad 1 \qquad && \qquad 2 \qquad && \qquad 3 \qquad
&& \qquad 4 \qquad && \qquad 5 \qquad && \qquad 6 \qquad&&
\qquad 7 \qquad& \cr \noalign {\hrule}
&&1&& 1&& 4&& 5&& 7&& 6&& 3&& 2&\cr
&&2&& 3&& 2&& 7&& 6&& 4&& 5&& 1&\cr
&&3&& 2&& 6&& 3&& 5&& 7&& 1&& 4&\cr
&&4&& 5&& 1&& 6&& 4&& 2&& 7&& 3&\cr
&&5&& 4&& 7&& 1&& 3&& 5&& 2&& 6&\cr
&&6&& 7&& 3&& 4&& 2&& 1&& 6&& 5&\cr
&&7&& 6&& 5&& 2&& 1&& 3&& 4&& 7&\cr
\noalign{\hrule}}}}
\vskip 1 cm
\centerline {Test 8}
\vskip 0.5 cm \centerline
{\vbox
{\offinterlineskip
\halign{\strut#&\vrule#\quad&\hfil#\hfil&\quad
\vrule#\quad&\hfil#\hfil&\quad#&\hfil#\hfil&\quad#&\hfil#\hfil&\quad
#&\hfil#\hfil&\quad#&\hfil#\hfil&\quad#&\hfil#\hfil&\quad
#&\hfil#\hfil&\quad\vrule#\cr
\noalign{\hrule}
&&\omit \quad \hfil && \qquad 1 \qquad && \qquad 2 \qquad && \qquad 3 \qquad
&& \qquad 4 \qquad && \qquad 5 \qquad && \qquad 6 \qquad&&
\qquad 7 \qquad& \cr \noalign {\hrule}
&&1&& 1&& 3&& 4&& 2&& 6&& 7&& 5&\cr
&&2&& 3&& 2&& 5&& 7&& 1&& 4&& 6&\cr
&&3&& 4&& 5&& 3&& 6&& 7&& 1&& 2&\cr
&&4&& 2&& 7&& 6&& 4&& 3&& 5&& 1&\cr
&&5&& 6&& 1&& 7&& 3&& 5&& 2&& 4&\cr
&&6&& 7&& 4&& 1&& 5&& 2&& 6&& 3&\cr
&&7&& 5&& 6&& 2&& 1&& 4&& 3&& 7&\cr
\noalign{\hrule}}}}
\vskip 1 cm
\centerline {Test 9}
\vskip 0.5 cm \centerline
{\vbox
{\offinterlineskip
\halign{\strut#&\vrule#\quad&\hfil#\hfil&\quad
\vrule#\quad&\hfil#\hfil&\quad#&\hfil#\hfil&\quad#&\hfil#\hfil&\quad
#&\hfil#\hfil&\quad#&\hfil#\hfil&\quad#&\hfil#\hfil&\quad
#&\hfil#\hfil&\quad\vrule#\cr
\noalign{\hrule}
&&\omit \quad \hfil && \qquad 1 \qquad && \qquad 2 \qquad && \qquad 3 \qquad
&& \qquad 4 \qquad && \qquad 5 \qquad && \qquad 6 \qquad&&
\qquad 7 \qquad& \cr \noalign {\hrule}
&&1&& 1&& 7&& 2&& 3&& 4&& 5&& 6&\cr
&&2&& 3&& 2&& 5&& 7&& 6&& 4&& 1&\cr
&&3&& 4&& 1&& 3&& 6&& 2&& 7&& 5&\cr
&&4&& 5&& 6&& 1&& 4&& 7&& 3&& 2&\cr
&&5&& 6&& 3&& 7&& 1&& 5&& 2&& 4&\cr
&&6&& 7&& 5&& 4&& 2&& 1&& 6&& 3&\cr
&&7&& 2&& 4&& 6&& 5&& 3&& 1&& 7&\cr
\noalign{\hrule}}}}
\vskip 13 cm
\centerline {$\underline {EIGHT \quad COLORS}$}
\vskip 0.3 cm
\centerline {Test 10}
\vskip 0.3 cm \centerline
{\vbox
{\offinterlineskip
\halign{\strut#&\vrule#\quad&\hfil#\hfil&\quad
\vrule#\quad&\hfil#\hfil&\quad#&\hfil#\hfil&\quad#&\hfil#\hfil&\quad
#&\hfil#\hfil&\quad#&\hfil#\hfil&\quad#&\hfil#\hfil&\quad
#&\hfil#\hfil&\quad#&\hfil#\hfil&\quad\vrule#\cr
\noalign{\hrule}
&&\omit \quad \hfil && \quad 1 \quad && \quad 2 \quad && \quad 3 \quad
&& \quad 4 \quad && \quad 5 \quad && \quad 6 \quad&&
\quad 7 \quad&& \quad 8 \quad & \cr \noalign {\hrule}
&&1&& 1&& 1&& 5&& 3&& 4&& 4&& 5&& 3&\cr
&&2&& 2&& 2&& 6&& 7&& 8&& 8&& 6&& 7&\cr
&&3&& 4&& 4&& 3&& 6&& 2&& 3&& 3&& 6&\cr
&&4&& 5&& 5&& 2&& 4&& 7&& 7&& 2&& 4&\cr
&&5&& 3&& 3&& 8&& 2&& 5&& 5&& 8&& 2&\cr
&&6&& 7&& 7&& 4&& 1&& 6&& 6&& 4&& 1&\cr
&&7&& 8&& 8&& 7&& 5&& 1&& 1&& 7&& 5&\cr
&&8&& 6&& 6&& 1&& 8&& 3&& 3&& 1&& 8&\cr
\noalign{\hrule}}}}
\vskip 0.3 cm
\centerline {Test 11}
\vskip 0.3 cm \centerline
{\vbox
{\offinterlineskip
\halign{\strut#&\vrule#\quad&\hfil#\hfil&\quad
\vrule#\quad&\hfil#\hfil&\quad#&\hfil#\hfil&\quad#&\hfil#\hfil&\quad
#&\hfil#\hfil&\quad#&\hfil#\hfil&\quad#&\hfil#\hfil&\quad
#&\hfil#\hfil&\quad#&\hfil#\hfil&\quad\vrule#\cr
\noalign{\hrule}
&&\omit \quad \hfil && \quad 1 \quad && \quad 2 \quad && \quad 3 \quad
&& \quad 4 \quad && \quad 5 \quad && \quad 6 \quad&&
\quad 7 \quad&& \quad 8 \quad& \cr \noalign {\hrule}
&&1&& 1&& 5&& 6&& 7&& 8&& 2&& 3&& 4&\cr
&&2&& 3&& 2&& 4&& 8&& 7&& 5&& 1&& 6&\cr
&&3&& 4&& 7&& 3&& 5&& 2&& 8&& 6&& 1&\cr
&&4&& 5&& 1&& 8&& 4&& 6&& 3&& 2&& 7&\cr
&&5&& 6&& 8&& 1&& 2&& 5&& 7&& 4&& 3&\cr
&&6&& 7&& 4&& 2&& 1&& 3&& 6&& 8&& 5&\cr
&&7&& 8&& 6&& 5&& 3&& 1&& 4&& 7&& 2&\cr
&&8&& 2&& 3&& 7&& 6&& 4&& 1&& 5&& 8&\cr
\noalign{\hrule}}}}
\vskip 0.3 cm
\par \medskip
Once these eleven tests were obtained, we applied them to the knots that
``survived" the equivalence algorithm. First we had to obtain the group of
each knot. The only difficulty here is being able to distinguish between
Figures 19a and 19b. We refer the reader to the end of paragraph 10 of Ref. [1]
about how to distinguish between these two cases.
\par \medskip
All knot projections, including the trivial ones, can yield a ``unicolor"
solution, that is a ``painting" where the entire projection is colored the
same. Such solutions are disregarded. In addition, one may notice that one
strand may be colored arbitrarily. If for example the color of some strand
changes from $a_i$ to $a_j$, all other colors $a_k$ may change to $f(a_k,b)$
where $f(a_i,b)=a_j$, and if the former coloring obeys the test, so will the
latter. While the existence of such $b$ is not always guaranteed, if the
test is irreducible, there is always going to be a sequence $b_1$, $b_2$,
..., $b_f$ such that $f \Big( ... \big( f(a_i,b_1),b_2 \big) ,...,b_f \Big)=
a_j$, and thus all other $a_k$ will be transformed accordingly. We therefore
imposed that only solutions ``coloring" the first strand with the first color
are going to count, thus reducing the running time by a factor equal to the
number of colors.
\par \medskip
After aplying these tests we got the following results. The trefoil
$\{ (1,4),(3,6),(5,2) \} $ proved its non-triviality from the very first test,
as was expected. The ``figure eight" knot, $\{ (1,4),(3,6),(5,8),(7,2) \} $
proved that
it is not the trefoil from failing the first test, and its non-triviality by
passing the second. The two five-crossing knots proved that they are
non-trivial, distinct from the trefoil and the ``figure eight" knot and
distinct from
each other as follows. They both failed the first two tests, but while
$\{ (1,6),(3,8),(5,10),(7,2),(9,4) \} $ passed the third,
$\{ (1,4),(3,8),(5,10),(7,2),(9,6) \} $ failed the first six tests and passed
the seventh. Similarly, two out of the three six-crossing knots proved that
they are distinct; the knot $\{ (1,4),(3,8),(5,10),(7,12),(9,2),(11,6) \} $,
on the other hand, failed all eleven tests.
\par \medskip
Consequently, after all the effort of developing these computer programs, and
converting them in order to fit the memory limitations and to reduce the
running times, all we were able to show was that there is one knot with
three crossings, one with four, and two with five. From six crossings and
higher the results would be ambiguous. Continuing along this route by
obtaining nine-color tests, while in
principle possible, required a practically ``infinite" amount of time. The
knot classification problem was thus merely solved for knots with up to five
crossings; this could hardly be called a success, since such a solution can
easily be obtained by hand in much less time than the time needed to develop
these algorithms.
\bigskip
\centerline {6. \qquad THE ``LINEAR" TESTS}
\bigskip
While obtaining all color tests may not be practical when the number of colors
is larger than eight, it is possible to get additional tests by generalising
the ones
already obtained for up to eight colors. While the generators of the three
color test, for instance, may be considered as the permutations
$\pmatrix {1&2&3\cr 2&1&3\cr}$, $\pmatrix {1&2&3\cr 3&2&1}$ and
$\pmatrix {1&2&3\cr 1&3&2\cr}$ satisfying $f(a_i,a_j)=a_ja_ia_j^{-1}$, they
may also be considered as the integer numbers $1$, $2$ and $3$ satisfying
$f(a_i,a_j)=2a_j-a_i \quad {\rm mod} 3$. Similarly the generators of the five
color tests may be considered as the numbers $1$, $2$, $3$, $4$ and $5$
satisfying $f(a_i,a_j)=4a_j-3a_i \quad {\rm mod} 5$ (test 3) and
$f(a_i,a_j)=2a_j-a_i
\quad {\rm mod} 5$ (test 4), although one would first have to rearrange them
by exchanging $2$ with $4$ in test 3 and $4$ with $5$ in test 4 to see this
relation. Similar relations are obeyed by the generators of the seven color
tests. One would thus desire to see whether one may obtain further such tests,
where the generators are related by $k_1a_i+k_2f(a_i,a_j)+k_3a_j=0 \quad
{\rm mod} k$; as discussed earlier, one would first have to find out whether
these tests are invariant under the Reidemeister moves.
\par \medskip
The condition imposed by the first move, $f(a_i,a_i)$, implies that
$k_1+k_2+k_3=0 \quad {\rm mod} k$. The condition due to the second move,
$f(a_i,a_j)=f(a_k,a_j)
\Leftrightarrow a_i=a_k$, implies that $k_1$ and $k_3$ are prime with
respect to $k$. Finally, the condition due to the third move may be written as
follows.
\par \medskip
For any natural numbers $a_i$, $a_j$, $a_l$, $a_m$, $a_n$, $a_p$ and $a_q$
satisfying
$$\alpha ) \qquad k_1a_i+k_2a_l+k_3a_j \quad = \quad 0 \quad {\rm mod} k$$
$$\beta ) \qquad k_1a_m+k_2a_n+k_3a_j \quad = \quad 0 \quad {\rm mod} k$$
$$\gamma ) \qquad k_1a_m+k_2a_p+k_3a_i \quad = \quad 0 \quad {\rm mod} k$$
$$\delta ) \qquad k_1a_n+k_2a_q+k_3a_l \quad = \quad 0 \quad {\rm mod} k$$
one must have $$k_1a_p+k_2a_q+k_3a_j \quad = \quad 0 \quad {\rm mod} k$$
By multiplying the
equation $\alpha )$ by $-k_3$, $\beta )$ by $-k_1$, $\gamma )$ by $k_1$,
$\delta )$ by $k_2$ and then adding together, one gets
$$-k_3^2a_j-k_1k_3a_j+k_1k_2a_p+k_2^2a_q=0 \quad {\rm mod} k$$
Since $k_1+k_2+k_3=0 \quad {\rm mod} k$, this may also be written as
$k_2(k_3a_j+k_1a_p+k_2a_q)=0 \quad {\rm mod} k$. Therefore, in order for the
condition to be satisfied, one must also impose $k_2 \neq 0 \quad {\rm mod}k$.
In fact, in order for the test to be irreducible, $k_2$ must be prime with
respect to $k$.
\par \medskip
Finally, if $k_1a_i+k_2a_j+k_3a_l=0 \quad {\rm mod} k$, obviously
$rk_1a_i+rk_2a_j+rk_3a_l=0 \quad {\rm mod}k$. In order to avoid such
repetitions, we impose $k_3$=1, and thus the ``linear" tests take the form
$f(a_i,a_j)=(s+1)a_j-sa_i \quad {\rm mod} k$. From now on a test defined by
such a relation is going to be called the {\it $[k,s]$ linear test}. In
addition,
if $s_1s_2=1 \quad {\rm mod} k$, then for each successful ``coloring" of a knot
under
the $[k,s_1]$ linear test, there is a successful ``coloring" of the mirror
symmetric knot under the $[k,s_2]$ linear test. Such tests are bound to
produce identical results, and therefore only one of them is performed,
chosen conventionally as the $[k,s]$ for $s=\min (s_1,s_2)$.
\par \medskip
As a corollary, since all $k_1$, $k_2$ and $k_3$ are prime with respect to
$k$, $k$ cannot be even, and thus the four, six and eight color tests obtained
before do not belong to this category of ``linear" tests. Therefore, even
though the method described here can yield an infinite number of tests, it is
unable to yield all of them. In Section 8 we are going to discuss more about
the even number color tests.
\par \medskip
Let $k=lm$, where $l$ and $m$ are prime with respect to each other. To
each successful coloring of $[k,s]$ corresponds exactly one coloring of
$[l,s {\rm mod} l]$ and one of $[m,s {\rm mod} m]$. Therefore if two knots
``agree" for all $l$ and $m$ color linear tests, they will also ``agree" for
all $k$ color linear tests. One need thus only apply tests whose number of
colors is either prime or a power of a prime number.
\par \medskip
Let now $i \rightarrow a_i$ be a successful coloring, where $i$ is a strand of
a knot and $a_i$ is its color. Then one may notice that $i \rightarrow a_i' =
c_1+c_2a_i$ is a successful coloring, too. This property reduces the running
time by a factor of up to $k(k-1)$; when $k$ is a higher than one power of a
prime
number however, one must make sure that $c_2$ is not a multiple of the prime
number, since then exchanging $a_i$ with $a_i'$ is not one-to-one.
\par \medskip
Using the theory described up to now we ran the corresponding
programmes by setting the cutoff parameters at $n=12$ for the equivalence
algorithm and at $49$ colors for the linear tests. Eventually $27$ colors
would prove enough to solve the knot classification problem up to $8$
crossings, since all knots that ``survived" the equivalence algorithm up to $8$
crossings ``reacted" differently from each other at some linear test of not
more than $27$ colors. That was not the case for the $9$ crossings knots; while
$46$ out of the $49$ knots proved to be distinct, there were three cases where
knots responded identically to knots that bore ``preferred" names.
\par \medskip
While the time needed to obtain new linear tests is infinitesimal, the
time needed to apply these tests becomes significant when the number of colors
is larger than 50, and in addition this amount of
time has to be multiplied by the number of tests, which is about half of
the number of colors, and by the number of the surviving knots. It might
be more productive instead of indefinitely running programmes with more and
more linear tests to find some method which may yield the result of a test in
advance. Since as we said earlier, the three color test decides whether
$P_3 \times Z$ is a coset of the knot group, the results of the tests may not
come randomly but be instead predetermined.
\vskip 2cm
\centerline {7. \qquad THE ``LINEAR" SYSTEM}
\bigskip
Running a computer program can yield very fast results when the crossing and
color numbers are small, but is not as efficient with large numbers. One
would therefore like to see if some solution can be found analytically in a
reasonable amount of time; one would also hope that such a solution would
provide more information about a knot than a computer output.
\par \medskip
Let $i \rightarrow a_i$ be a successful coloring of a knot projection under
some linear test $f(a_i,a_j)=(s+1)a_j-sa_i$. The constraints at the crossings
of the projection impose conditions $a_{i \pm 1}=(s+1)a_j-sa_i$, where the
strand $j$ separates the strands $i$ and $i \pm 1$. There are as many such
equations and unknowns as the crossing number, and thus one gets a homogeneous
system. Such systems always accept a solution $a_i=0 \quad \forall \quad
i$; in addition, this system accepts
non-zero solutions $a_i=c \neq 0 \quad \forall \quad i$, which implies
that the determinant is $0$. Both of these solutions correspond to a
``unicolor" painting, where the whole knot is painted the same. Since such
solutions exist for all knots, they are not helpful in proving inequivalencies.
\par \medskip
The situation is improved once the unknowns $a_i$ are replaced by
$x_i=a_{i+1}-a_i$ ($x_n=a_1-a_n$) and look for non-zero solutions for $x_i$.
This transformation is linear, and thus the new system remains linear. The
Jacobian of the transformation, $\rm {det} {\vartheta x_i \over \vartheta
a_j}=0$, and thus the transformation is singular. The acceptable solutions are
the ones for which the determinant of
the new system is $0$, and also the sum of the non-zero values of $x_i$ is $0$.
The latter is due to the fact that $x_i=a_{i+1}-a_i$, and therefore
$x_1+x_2+...+x_{n-1}+x_n=(a_2-a_1)+(a_3-a_2)+...+(a_n-a_{n-1})+(a_1-a_n)=0$.
\par \medskip
By applying the theory discussed above one finds the following. First, the
determinant of the system for the unknot is always a non-zero constant, which
is equal to
$1$ if one uses the $0$ crossings projection. This is hardly a surprise since
the unknot is supposed to fail all color tests, which is guaranteed by a
non-zero determinant. The unknot's failure is of course due to the fact that
its $0$ crossings projection cannnot be painted using more than one colors.
\par \medskip
For the trefoil one gets the systems
$$(s+1)x_1+x_2=0, \qquad (s+1)x_2+x_3=0,
\qquad (s+1)x_3+x_1=0$$ and
$$(s+1)x_1+sx_2=0, \qquad (s+1)x_2+sx_3=0, \qquad
(s+1)x_3+sx_1=0$$ depending on whether the the trefoil is right-handed or
left-handed. The two determinants are equal to $(s+1)^3+1$ and $(s+1)^3+s^3$
respectively. These are third degree polynomials whose roots are
$-1+\rm {exp} \big( {2\pi \nu \over 3} \big)$ and
$\Big( -1+\rm {exp} \big( {2\pi \nu \over 3} \big) \Big) ^{-1}$ respectively,
where $\nu \in
\{ 1,2,3 \} $.
\par \medskip
Such complex roots yield of course complex solutions for $x_i$,
which is hardly what one would expect; until now we were looking for integer
solutions, and were in fact able to find a number of them.
\par \medskip
This apparent paradox is resolved once we realise that all equations written up
to this point are modulo some unspecified integer number $k$. Therefore,
instead of trying to find the complex roots of the determinant, which may be
analytically impossible if the system consists of five or more equations, one
may simply set the determinant equal to $0$ by specifying $k$ to be a divisor
of the determinant.
\par \medskip
Not all such divisors are suitable; one also has to check that the sum of $x_i$
is also $0 {\rm mod} k$. Take for instance Listing's knot, whose system is
$(s+1)x_1+x_2=0, \quad (s+1)x_2+sx_3=0, \quad (s+1)x_3+x_4=0, \quad
(s+1)x_4+sx_1=0$. The determinant of this system is $(s+1)^4-s^2$, which is
equal to $15$ when $s=1$. Since the (prime) divisors of $15$ are $3$ and $5$,
the tests $[3,1]$ and $[5,1]$ stand a chance. Solving now the system
one gets $x_1=1$, $x_2=-2$, $x_3=4$, $x_4=-8$ and thus $x_1+x_2+x_3+x_4=-5$.
This becomes $0 {\rm mod} k$ for $k=5$ but not for $k=3$, and thus the test
$5,1$ is successful while $[3,1]$ is not. This result is in agreement
with the computer output.
One may similarly work for $s>1$ as well as with other knots and be able to
predetermine the test results.
\par \medskip
Time saving is not the only advantage of such an approach. One may have
noticed that while for the trefoil the roots of the determinant depend on
chirality, this is not the case for the unknot or for Listing's knot. Were
these roots to remain invariant under Reidemeister moves, one would establish
a criterion that could detect chiral knots.
When one performs Reidemeister moves that increases or decreases the crossing
number by $t$, the determinant is multiplied or divided by $(s+1)^t$.
This is going to add or subtract $t$ roots equal to $-1$. Therefore when
comparing the roots of the two determinants, one would first have to ignore
roots equal to $-1$.
\par \medskip
One does not need however to find independently the two determinants, obtain
their roots and ignore the $-1$ roots to detect chirality.
When a knot projection is replaced by its mirror image, the
equations $a_{i \pm 1}=(s+1)a_j-sa_i$ are replaced by
$a_i=(s+1)a_j-sa_{i \pm 1} \Leftrightarrow a_{i \pm 1}=(s^{-1}+1)a_j-s^{-1}a_i$
and thus $s$ is replaced by $s^{-1}$. Therefore if the knot is acheiral, for
each root $s$, $s^{-1}$ will also be a root. Performing Reidemeister moves and
adding or losing $-1$ roots does not pose a problem, since the inverse of $-1$
is also $-1$. Therefore the existence of a root $s$
when $s^{-1}$ is not a root, proves a knot's chirality.
\par \medskip
One cannot of course use this method to prove mirror symmetry; the roots may
come in pairs $(s,s^{-1})$ even if a knot is chiral. In order to prove that a
knot is acheiral one would have to modify the equivalence algorithm so that
mirror symmetric knots are not automatically considered equivalent. (See
Section 11, Ref. 1)
\par \medskip
One would desire at this point to know whether it is possible to obtain a new
knot characteristic that would play a role similar to the knot polynomials.
Since the determinant is multiplied (divided) by $(s+1)^t$ when the knot is
replaced by an equivalent whose crossing number, $c$, is larger (smaller) by
$t$, it is not the determinant itself, but the ratio of the determinant divided
by $(s+1)^c$, that can play the role of a knot
characteristic.
Such a characteristic would carry in fact some important advantages. First, it
can be obtained without using any skein relations. It thus requires much less
memory and running time than obtaining any of the known knot polynomials. The
most involved part of the calculation is just the calculation of the
determinant of a matrix whose number of rows and columns is equal to the
crossing number, and whose elements can take one of the values $0$, $1$, $s$
and $s+1$.
Second, it can be useful in detecting chirality; indeed, unless $c_i=c_{r-i}
\forall i \in \{ 0,1,2,...,r \} $, the knot is chiral. Until now the only
polynomial known to detect chirality is the HOMFLY polynomial, which depends
on two variables, while the one obtained from the determinant depends on just
one variable. Finally, it provides a lower value for a knot's crossing number;
let $f(s)={P(s) \over (s+1)^{\nu}}$ be a knot's characteristic such that
$P(s)$ is a polynomial satisfying $P(-1) \neq 0$. The knot's crossing number is
greater or equal to $\nu$.
\par \medskip
A careful reader may detect a similarity between this characteristic and the
Alexander polynomial, since the latter may also be defined from the determinant
of a linear system. The similarity however ends here, as the system
corresponding to the Alexander polynomial involves the unknowns $a_i$ instead
of $x_i$. Since the transformation from $a_i$ to $x_i$ is singular, the two
systems and their determinants are distinct.
\par \medskip
In spite of these advantages the characteristic defined above is not ``sharp"
enough;
as one may notice from the tables listed in Section 9, it fails to distinguish
among one pair of 8 crossing knots and two pairs of 9 crossing knots. Since
these pairs do not coincide with the pairs that were ``tied" after the linear
tests, the knot classification problem is now solved up to 9 crossings; all
84 knots surviving the equivalence algorithm either possess different
characteristics, or react differently to linear tests. In fact one need not go
further than the 5 color tests to complete the classification.
\par \medskip
Having solved the problem up to 9 crossings, we proceeded with the 165 knots
with 10 crossings which were obtained from the equivalence algorithm. These
knots possessed 156 distinct polynomials. Of the remaining 9 knots, 8 were
able to show their distinctess through linear tests of up to 7 colors, while
the last knot was shown distinct through the four-color non-linear test.
Therefore the problem is now solved up to 10 crossings.
\par \medskip
Finally we would like to comment on ``multiple successes". Most knots will
either fail some test, or if they do pass it, there is just one way to be
``painted", apart from the symmetry $a_i \rightarrow c_1+c_2a_i$. For a few
knots
however this is not the case. The knot $\{ (1,6),(3,8),(5,10),(7,12),(9,14),
(11,16),(13,2),(15,4) \} $ is such an example. It not only passes the test
$[3,1]$, but it does so in $4$ distinct ways (Figures 22 and 23).
It also passes the test $[7,2]$ in $8$ ways, the test $[13,3]$ in
$14$ ways and the test $[19,7]$ in $20$. One may notice that it passes
$[k,s]$ tests in $k+1$ ways. This occurs beacuse its characteristic possesses
multiple roots. As one may actually check from the tables of Section 9, the
characteristic of this knot, which is ordered $29^{\rm th}$, is
$8(s^4+2s^3+3s^2+2s+1)(s+1)^{-6}$, which is equal to $8(s^2+s+1)^2(s+1)^{-6}$,
and thus possesses two double roots, ${-1 \pm i \sqrt3 \over 2}$.
Consequently, its solutions lie on a two-dimensional manifold; even if one
$x_i$ is fixed, the rest of them are not uniquely
determined, unless a $x_i \neq x_j$ is also fixed. Such cases are
obviously rare, and this example is the only case for up to eight crossings.
\bigskip
\centerline {8. \qquad THE ``POLYHEDRAL" TESTS}
\bigskip
While the ``linear" tests discussed in the two previous Sections can account
for the odd color tests listed in Section 5, we have yet to study the
even color tests. A careful look at the four color test reveals a striking
resemblance to the symmetry group of a regular tetrahedron. Indeed, each
element of the matrix $M_{ij}$ coincides with a $2\pi \over 3$ rotation of
vertex $i$ with respect to the axis passing through vertex $j$. One could thus
generalise this argument and introduce $6$, $8$, $12$ and $20$ color tests
related to symmetry groups of the regular cube, octahedron, dodekahedron
and eikosahedron.
In each such case $M_{ij}$ is the
vertex obtained from $i$ after a fixed angle rotation with respect to the axis
going from the center to vertex $j$. Therefore, in order for the number of
vertices to be finite, they should form a regular polyhedron.
\par \medskip
Such tests can easily be shown to obey all conditions imposed by Reidemeister
moves which were discussed in Section 5. They may not however be irreducible,
as we are now going to show.
Take for example the octahedron (Figure 24). If one defines as $M_{ij}$ the
position of vertex $i$ after it is rotated by $\pi \over 2$ with respect to the
axis passing through $j$, one indeed gets Test 6, which is expected. A
$\pi $ rotation however does not yield Test 5, but a reducible test, since
each vertex could be mapped either to itself or to its opposite. Therefore
such a test would be reduced to three two-color tests, each of them being
further reduced to one-color tests, and one would merely get trivial tests. A
similar problem is encountered if one works now on the cube (Figure 25), by
defining $M_{ij}$ as a $2\pi \over 3$ rotation of $i$ with respect to the
axis passing through $j$; this time the test would be reduced to two four-color
tests, one involving the vertices 1,3,5,7, and the other the vertices 2,4,6,8.
Therefore just two out of the five even color tests, namely Test 2 and Test 6,
can be accounted as
regular polyhedra tests.
\par \medskip
If one generalises such tests to include reflections in addition to rotations,
one is able to obtain more tests. One may define for instance
$M_{ij}$ to be not a rotation but a reflection on an octahedron with respect
to a plane passing through its center, vertex $j$ and the middle of a side
perpendicular to the axis passing through $j$. Similarly, if $j'$ is the
vertex opposite to $j$, $M_{ij'}$ is defined as a reflection with respect to a
plane passing through the center and $j'$, but perpendicular to the previous
one (Figure 26). Under these requirements one may notice that the test obtained
is identical to Test 5. Similarly, through a more involved definition that
combines rotations and reflections on a cube, one may come up with Test 10.
Test 11 however remains unaccounted, which shows that there is more to the
theory of ``color tests" than the ``linear" tests discussed before and the
``polyhedral" tests introduced here.
\bigskip
\centerline {9. \qquad RESULTS AND CONCLUSION}
\bigskip
After running the equivalence algorithm for $n=12$, a total of 45517 knots
``survived", as one may see from the table of Section 4. This
means that no two of them were shown to be equivalent to each
other. This of course does not prove their inequivalence, since some knots in
order to be shown equivalent, require a series of Reidemeister moves that
involves knots of $13$ crossings or more. Starting from the next page, we
shall list all knots
with up to ten crossings, since these knots have been proved inequivalent.
Then we shall list their characteristics as defined in Section 7, and record
their ``responses" to a few linear tests so that their distinctness can be
shown unambiguously.
Here the word ``response" means the number of distinct
ways in which a knot can be ``painted" according to some linear test; symmetric
``colorings" related by $a'_i=c_1+c_2a_i$ are counted just once. (There is
however an exception, as mentioned in Section 7. The knots ordered
$73^{\rm rd}$ and $165^{\rm th}$ possess identical polynomials and respond
identically to linear tests. Their distinctness can be demonstrated through
the four-color tetrahedral test; while there are five distinct ways to
``paint" the $165^{\rm th}$ knot, there is just one for the $73^{\rm rd}$.)
\par \medskip
To summarise now the work discussed in this paper, our goal has been to fully
classify all distinct knots, and this has been carried out in two stages.
First, knots shown to be equivalent were eliminated, with the
exception of one knot from each equivalence class, the one bearing the most
``preferrable" name. Since however this is an infinite process, we had to set
a cutoff parameter; this parameter was set at $n=12$, implying that knot
projections with up to $12$ crossings would be used to identify possible
equivalences, and thus any series of equivalence moves involving projections
with $13$ or more double points would not be considered. Second, the knots
that were not previously eliminated, would be tested in order to show their
distinctness. The most efficient test we were able to obtain is related to a
new characteristic introduced in this paper. A few cases that were not resolved
through this characteristic would be resolved through the {\it color
tests} which are a generalisation of Reidemeister's tricolorability.
\par \medskip
Using these techniques we were able to solve the knot classification problem
up to $10$ crossings in about two days. It is our hope and belief that
given sufficient computer memory and running time, the theory and the
algorithms discussed in this paper can help solve the problem up to any
crossing number.
\vfil \eject
\eightpoint
{\nopagenumbers
\def\*{\hphantom{0}}
\centerline
{\vbox
{\offinterlineskip
\halign{\strut#&\vrule#&\hfil#\hfil&\vrule#&\hfil#\hfil&\
\vrule#&\hfil#\hfil&\vrule#\cr
\noalign{\hrule}
&&\omit {\bf Order} && {\bf Crossings}&&
{\bf Knot's Name} &\cr \noalign {\hrule}
&&\* \* 0&& \* 0&& $\{ \}$&\cr
&&\* \* 1&& \* 3&& $\{ (1,4),(3,6),(5,2) \}$&\cr
&&\* \* 2&& \* 4&& $\{ (1,4),(3,6),(5,8),(7,2) \}$&\cr
&&\* \* 3&& \* 5&& $\{ (1,4),(3,8),(5,10),(7,2),(9,6) \}$&\cr
&&\* \* 4&& \* 5&& $\{ (1,6),(3,8),(5,10),(7,2),(9,4) \}$&\cr
&&\* \* 5&& \* 6&& $\{ (1,4),(3,8),(5,10),(7,2),(9,12),(11,6) \}$&\cr
&&\* \* 6&& \* 6&& $\{ (1,4),(3,8),(5,10),(7,12),(9,2),(11,6) \}$&\cr
&&\* \* 7&& \* 6&& $\{ (1,4),(3,8),(5,12),(7,10),(9,2),(11,6) \}$&\cr
&&\* \* 8&& \* 7&& $\{ (1,4),(3,8),(5,10),(7,12),(9,2),(11,14),(13,6) \}$&\cr
&&\* \* 9&& \*
7&& $\{ (1,4),(3,8),(5,12),(7,2),(9,14),(11,6),(13,10) \}$&\cr
&&\* 10&& \* 7&& $\{ (1,4),(3,10),(5,12),(7,14),(9,2),(11,8),(13,6) \}$&\cr
&&\* 11&& \* 7&& $\{ (1,4),(3,10),(5,14),(7,12),(9,2),(11,8),(13,6) \}$&\cr
&&\* 12&& \* 7&& $\{ (1,6),(3,10),(5,12),(7,14),(9,2),(11,4),(13,8) \}$&\cr
&&\* 13&& \* 7&& $\{ (1,6),(3,10),(5,12),(7,14),(9,4),(11,2),(13,8) \}$&\cr
&&\* 14&& \* 7&& $\{ (1,8),(3,10),(5,12),(7,14),(9,2),(11,4),(13,6) \}$&\cr
&&\* 15&& \* 8&& $\{ (1,4),(3,8),(5,10),(7,14),(9,2),(11,16),(13,6),(15,12)
\}$&\cr
&&\* 16&& \* 8&& $\{ (1,4),(3,8),(5,12),(7,2),(9,14),(11,6),(13,16),(15,10)
\}$&\cr
&&\* 17&& \* 8&& $\{ (1,4),(3,8),(5,12),(7,2),(14,9),(11,6),(16,13),(10,15)
\}$&\cr
&&\* 18&& \* 8&& $\{ (1,4),(3,8),(12,5),(7,2),(9,14),(6,11),(13,16),(15,10)
\}$&\cr
&&\* 19&& \* 8&& $\{ (1,4),(3,8),(5,12),(7,2),(9,14),(11,16),(13,6),(15,10)
\}$&\cr
&&\* 20&& \* 8&& $\{ (1,4),(3,8),(12,5),(7,2),(14,9),(16,11),(6,13),(10,15)
\}$&\cr
&&\* 21&& \* 8&& $\{ (1,4),(3,8),(5,12),(7,2),(9,16),(11,14),(13,6),(15,10)
\}$&\cr
&&\* 22&& \* 8&& $\{ (1,4),(3,8),(5,14),(7,10),(9,2),(11,16),(13,6),(15,12)
\}$&\cr
&&\* 23&& \* 8&& $\{ (1,4),(3,10),(5,12),(7,14),(9,2),(11,16),(13,6),(15,8)
\}$&\cr
&&\* 24&& \* 8&& $\{ (1,4),(3,10),(5,12),(7,14),(9,2),(11,16),(13,8),(15,6)
\}$&\cr
&&\* 25&& \* 8&& $\{ (1,4),(3,10),(5,12),(7,14),(9,16),(11,2),(13,6),(15,8)
\}$&\cr
&&\* 26&& \* 8&& $\{ (1,4),(3,10),(5,12),(7,14),(9,16),(11,2),(13,8),(15,6)
\}$&\cr
&&\* 27&& \* 8&& $\{ (1,4),(3,10),(5,14),(7,16),(9,12),(11,2),(13,8),(15,6)
\}$&\cr
&&\* 28&& \* 8&& $\{ (1,4),(3,10),(5,16),(7,14),(9,12),(11,2),(13,8),(15,6)
\}$&\cr
&&\* 29&& \* 8&& $\{ (1,6),(3,8),(5,10),(7,12),(9,14),(11,16),(13,2),(15,4)
\}$&\cr
&&\* 30&& \* 8&& $\{ (1,6),(3,8),(5,12),(7,2),(9,14),(11,16),(13,4),(15,10)
\}$&\cr
&&\* 31&& \* 8&& $\{ (1,6),(3,8),(5,12),(7,14),(9,4),(11,16),(13,2),(15,10)
\}$&\cr
&&\* 32&& \* 8&& $\{ (1,6),(3,8),(5,14),(7,12),(9,4),(11,16),(13,2),(15,10)
\}$&\cr
&&\* 33&& \* 8&& $\{ (1,6),(3,10),(5,12),(7,14),(9,16),(11,4),(13,2),(15,8)
\}$&\cr
&&\* 34&& \* 8&& $\{ (1,6),(3,10),(5,12),(7,16),(9,14),(11,4),(13,2),(15,8)
\}$&\cr
&&\* 35&& \* 8&& $\{ (1,6),(3,12),(5,10),(7,16),(9,14),(11,4),(13,2),(15,8)
\}$&\cr
&&\* 36&&\* 9&& $\{ (1,4),(3,8),(5,10),(7,14),(9,2),(11,16),(13,6),(15,18),
(17,12) \}$&\cr
&&\* 37&&\* 9&& $\{ (1,4),(3,8),(5,10),(7,14),(9,2),(16,11),(13,6),(18,15),
(12,17) \}$&\cr
&&\* 38&&\* 9&& $\{ (1,4),(3,8),(5,10),(14,7),(9,2),(11,16),(6,13),(15,18),
(17,12) \}$&\cr
&&\* 39&&\* 9&& $\{ (1,4),(3,8),(5,10),(14,7),(9,2),(16,11),(6,13),(18,15),
(12,17) \}$&\cr
&&\* 40&&\* 9&& $\{ (1,4),(3,8),(5,10),(7,14),(9,2),(11,16),(13,18),(15,6),
(17,12) \}$&\cr
&&\* 41&&\* 9&& $\{ (1,4),(3,8),(5,10),(14,7),(9,2),(16,11),(18,13),(6,15),
(12,17) \}$&\cr
&&\* 42&&\* 9&& $\{ (1,4),(3,8),(5,10),(7,14),(9,2),(11,18),(13,16),(15,6),
(17,12) \}$&\cr
&&\* 43&&\* 9&& $\{ (1,4),(3,8),(5,12),(7,2),(9,16),(11,6),(13,18),(15,10),
(17,14) \}$&\cr
&&\* 44&&\* 9&& $\{ (1,4),(3,8),(5,12),(7,2),(9,16),(11,14),(13,6),(15,18),
(17,10) \}$&\cr
&&\* 45&&\* 9&& $\{ (1,4),(3,8),(5,12),(7,14),(9,2),(11,16),(13,18),(15,10),
(17,6) \}$&\cr
&&\* 46&&\* 9&& $\{ (1,4),(3,8),(5,14),(7,2),(9,16),(11,18),(13,6),(15,12),
(17,10) \}$&\cr
&&\* 47&&\* 9&& $\{ (1,4),(3,8),(5,14),(7,2),(9,18),(11,16),(13,6),(15,12),
(17,10) \}$&\cr
&&\* 48&&\* 9&& $\{ (1,4),(3,8),(5,14),(7,10),(9,2),(11,16),(13,18),(15,6),
(17,12) \}$&\cr
&&\* 49&&\* 9&& $ \hphantom{0} \hphantom{0} \hphantom{0} \hphantom {0}
\{ (1,4),(3,8),(5,14),(7,10),(9,2),(11,18),(13,16),(15,6),(17,12)
\} \hphantom{0} \hphantom{0} \hphantom{0} \hphantom{0} $&\cr
&&\* 50&&\* 9&& $\{ (1,4),(3,8),(5,14),(7,12),(9,2),(11,16),(13,18),(15,10),
(17,6) \}$&\cr
&&\* 51&&\* 9&& $\{ (1,4),(3,10),(5,12),(7,14),(9,2),(11,18),(13,16),(15,6),
(17,8) \}$&\cr
&&\* 52&&\* 9&& $\{ (1,4),(3,10),(5,12),(7,14),(9,2),(11,18),(13,16),(15,8),
(17,6) \}$&\cr
&&\* 53&&\* 9&& $\{ (1,4),(3,10),(5,12),(7,14),(9,16),(11,2),(13,6),(15,18),
(17,8) \}$&\cr
&&\* 54&&\* 9&& $\{ (1,4),(3,10),(5,12),(7,14),(9,16),(11,2),(13,18),(15,8),
(17,6) \}$&\cr
&&\* 55&&\* 9&& $\{ (1,4),(3,10),(5,12),(7,16),(9,2),(11,8),(13,18),(15,6),
(17,14) \}$&\cr
&&\* 56&&\* 9&& $\{ (1,4),(3,10),(5,12),(7,16),(9,14),(11,2),(13,18),(15,8),
(17,6) \}$&\cr
&&\* 57&&\* 9&& $\{ (1,4),(3,10),(5,14),(7,12),(9,16),(11,2),(13,6),(15,18),
(17,8) \}$&\cr
&&\* 58&&\* 9&& $\{ (1,4),(3,10),(5,14),(7,12),(16,9),(11,2),(13,6),(18,15),
(8,17) \}$&\cr
&&\* 59&&\* 9&& $\{ (1,4),(3,10),(14,5),(12,7),(9,16),(11,2),(6,13),(15,18),
(17,8) \}$&\cr
&&\* 60&&\* 9&& $\{ (1,4),(3,10),(5,14),(7,16),(9,2),(11,18),(13,8),(15,6),
(17,12) \}$&\cr
&&\* 61&&\* 9&& $\{ (1,4),(3,10),(5,14),(7,16),(9,12),(11,2),(13,18),(15,6),
(17,8) \}$&\cr
&&\* 62&&\* 9&& $\{ (1,4),(3,10),(5,14),(7,16),(9,12),(11,2),(13,18),(15,8),
(17,6) \}$&\cr
\noalign{\hrule}}}}
\vskip 1 cm \centerline
{\vbox
{\offinterlineskip
\halign{\strut#&\vrule#&\hfil#\hfil&\vrule#&\hfil#\hfil&\
\vrule#&\hfil#\hfil&\vrule#\cr
\noalign{\hrule}
&&\omit {\bf Order} && {\bf Crossings}&&
{\bf Knot's Name} &\cr \noalign {\hrule}
&&\* 63&&\* 9&& $\{ (1,4),(3,10),(5,16),(7,14),(9,2),(11,18),(13,8),(15,6),
(17,12)
\}$&\cr
&&\* 64&&\* 9&& $\{ (1,4),(3,12),(5,14),(7,16),(9,18),(11,2),(13,10),(15,6),
(17,8)
\}$&\cr
&&\* 65&&\* 9&& $\{ (1,4),(3,12),(5,14),(7,16),(9,18),(11,2),(13,10),(15,8),
(17,6)
\}$&\cr
&&\* 66&&\* 9&& $\{ (1,4),(3,12),(5,16),(7,18),(9,14),(11,2),(13,8),(15,10),
(17,6)
\}$&\cr
&&\* 67&&\* 9&& $\{ (1,4),(3,12),(5,16),(7,18),(9,14),(11,2),(13,10),(15,8),
(17,6)
\}$&\cr
&&\* 68&&\* 9&& $\{ (1,4),(3,12),(5,18),(7,16),(9,14),(11,2),(13,10),(15,8),
(17,6)
\}$&\cr
&&\* 69&&\* 9&& $\{ (1,6),(3,8),(5,10),(7,16),(9,14),(11,18),(13,4),(15,2),
(17,12)
\}$&\cr
&&\* 70&&\* 9&& $\{ (1,6),(3,8),(5,10),(7,16),(9,14),(18,11),(13,4),(15,2),
(12,17)
\}$&\cr
&&\* 71&&\* 9&& $\{ (1,6),(3,10),(5,14),(7,12),(9,16),(11,2),(13,18),(15,4),
(17,8)
\}$&\cr
&&\* 72&&\* 9&& $\{ (1,6),(10,3),(14,5),(7,12),(16,9),(2,11),(13,18),(4,15),
(8,17)
\}$&\cr
&&\* 73&&\* 9&& $\{ (1,6),(3,10),(5,14),(7,18),(9,4),(11,16),(13,2),(15,8),
(17,12)
\}$&\cr
&&\* 74&&\* 9&& $\{ (1,6),(3,10),(5,14),(7,18),(9,4),(11,16),(13,8),(15,2),
(17,12)
\}$&\cr
&&\* 75&&\* 9&& $\{ (1,6),(3,10),(5,14),(7,18),(9,16),(11,2),(13,8),(15,4),
(17,12)
\}$&\cr
&&\* 76&&\* 9&& $\{ (1,6),(3,12),(5,14),(7,16),(9,18),(11,2),(13,4),(15,10),
(17,8)
\}$&\cr
&&\* 77&&\* 9&& $\{ (1,6),(3,12),(5,14),(7,16),(9,18),(11,4),(13,2),(15,10),
(17,8)
\}$&\cr
&&\* 78&&\* 9&& $\{ (1,6),(3,12),(5,14),(7,18),(9,16),(11,2),(13,4),(15,10),
(17,8)
\}$&\cr
&&\* 79&&\* 9&& $\{ (1,6),(3,12),(5,14),(7,18),(9,16),(11,4),(13,2),(15,10),
(17,8)
\}$&\cr
&&\* 80&&\* 9&& $\{ (1,6),(3,16),(5,14),(7,12),(9,4),(11,2),(13,18),(15,10),
(17,8)
\}$&\cr
&&\* 81&&\* 9&& $\{ (1,8),(3,12),(5,14),(7,16),(9,18),(11,2),(13,4),(15,6),
(17,10)
\}$&\cr
&&\* 82&&\* 9&& $\{ (1,8),(3,12),(5,14),(7,16),(9,18),(11,2),(13,6),(15,4),
(17,10)
\}$&\cr
&&\* 83&&\* 9&& $\{ (1,8),(3,12),(5,16),(7,14),(9,18),(11,4),(13,2),(15,6),
(17,10)
\}$&\cr
&&\* 84&&\* 9&& $\{ (1,10),(3,12),(5,14),(7,16),(9,18),(11,2),(13,4),(15,6),
(17,8)
\}$&\cr
&&\* 85&&10&& $\{ (1,4),(3,8),(5,10),(7,14),(9,2),(11,16),(13,18),(15,6),
(17,20),(19,12) \}$&\cr
&&\* 86&&10&& $\{ (1,4),(3,8),(5,10),(7,14),(9,2),(16,11),(18,13),(15,6),
(20,17),(12,19) \}$&\cr
&&\* 87&&10&& $\{ (1,4),(3,8),(5,10),(14,7),(9,2),(11,16),(13,18),(6,15),
(17,20),(19,12) \}$&\cr
&&\* 88&&10&& $\{ (1,4),(3,8),(5,10),(7,14),(9,2),(11,18),(13,6),(15,20),
(17,12),(19,16) \}$&\cr
&&\* 89&&10&& $\{ (1,4),(3,8),(5,10),(14,7),(9,2),(18,11),(6,13),(20,15),
(12,17),(16,19) \}$&\cr
&&\* 90&&10&& $\{ (1,4),(3,8),(5,10),(7,14),(9,2),(11,18),(13,16),(15,6),
(17,20),(19,12) \}$&\cr
&&\* 91&&10&& $\{ (1,4),(3,8),(5,10),(7,16),(9,2),(11,18),(13,20),(15,6),
(17,14),(19,12) \}$&\cr
&&\* 92&&10&& $\{ (1,4),(3,8),
(5,10),(7,16),(9,2),(11,20),(13,18),(15,6),(17,14),(19,12)
\}$&\cr
&&\* 93&&10&& $\{ (1,4),(3,8),
(5,12),(7,2),(9,14),(11,18),(13,6),(15,20),(17,10),(19,16)
\}$&\cr
&&\* 94&&10&& $\{ (1,4),(3,8),
(5,12),(7,2),(14,9),(18,11),(13,6),(20,15),(10,17),(16,19)
\}$&\cr
&&\* 95&&10&& $\{ (1,4),(3,8),
(12,5),(7,2),(9,14),(11,18),(6,13),(15,20),(17,10),(19,16)
\}$&\cr
&&\* 96&&10&& $\{ (1,4),(3,8),
(12,5),(7,2),(14,9),(18,11),(6,13),(20,15),(10,17),(16,19)
\}$&\cr
&&\* 97&&10&& $\{ (1,4),(3,8),
(5,12),(7,2),(9,16),(11,6),(13,18),(15,10),(17,20),(19,14)
\}$&\cr
&&\* 98&&10&& $\{ (1,4),(3,8),
(5,12),(7,2),(9,16),(11,6),(18,13),(15,10),(20,17),(14,19)
\}$&\cr
&&\* 99&&10&& $\{ (1,4),(3,8),
(5,12),(7,2),(16,9),(11,6),(13,18),(10,15),(17,20),(19,14)
\}$&\cr
&&100&&10&& $\{ (1,4),(3,8),
(5,12),(7,2),(16,9),(11,6),(18,13),(10,15),(20,17),(14,19)
\}$&\cr
&&101&&10&& $\{ (1,4),(3,8),
(5,12),(7,2),(9,16),(11,6),(13,18),(15,20),(17,10),(19,14)
\}$&\cr
&&102&&10&& $\{ (1,4),(3,8),
(5,12),(7,2),(16,9),(11,6),(18,13),(20,15),(10,17),(14,19)
\}$&\cr
&&103&&10&& $\{ (1,4),(3,8),
(12,5),(7,2),(9,16),(6,11),(13,18),(15,20),(17,10),(19,14)
\}$&\cr
&&104&&10&& $\{ (1,4),(3,8),
(12,5),(7,2),(16,9),(6,11),(18,13),(20,15),(10,17),(14,19)
\}$&\cr
&&105&&10&& $\{ (1,4),(3,8),
(5,12),(7,2),(9,16),(11,6),(13,20),(15,18),(17,10),(19,14)
\}$&\cr
&&106&&10&& $\{ (1,4),(3,8),
(12,5),(7,2),(16,9),(6,11),(20,13),(18,15),(10,17),(14,19)
\}$&\cr
&&107&&10&& $\{ (1,4),(3,8),
(5,12),(7,2),(9,18),(11,14),(13,6),(15,20),(17,10),(19,16)
\}$&\cr
&&108&&10&& $\{ (1,4),(3,8),
(5,12),(7,14),(9,2),(11,16),(13,20),(15,18),(17,10),(19,6)
\}$&\cr
&&109&&10&& $\{ (1,4),(3,8),
(5,12),(7,18),(9,2),(11,16),(13,20),(15,6),(17,10),(19,14)
\}$&\cr
&&110&&10&& $\{ (1,4),(3,8),
(12,5),(18,7),(9,2),(16,11),(20,13),(6,15),(10,17),(14,19)
\}$&\cr
&&111&&10&& $\{ (1,4),(3,8),
(5,14),(7,2),(9,16),(11,18),(13,6),(15,20),(17,10),(19,12)
\}$&\cr
&&112&&10&& $\{ (1,4),(3,8),
(5,14),(7,2),(16,9),(18,11),(13,6),(20,15),(10,17),(12,19)
\}$&\cr
&&113&&10&& $\{ (1,4),(3,8),
(14,5),(7,2),(9,16),(11,18),(6,13),(15,20),(17,10),(19,12)
\}$&\cr
&&114&&10&& $\{ (1,4),(3,8),
(14,5),(7,2),(16,9),(18,11),(6,13),(20,15),(10,17),(12,19)
\}$&\cr
&&115&&10&& $\{ (1,4),(3,8),
(5,14),(7,2),(9,16),(11,18),(13,6),(15,20),(17,12),(19,10)
\}$&\cr
&&116&&10&& $\{ (1,4),(3,8),
(5,14),(7,2),(16,9),(18,11),(13,6),(20,15),(12,17),(10,19)
\}$&\cr
&&117&&10&& $\{ (1,4),(3,8),
(14,5),(7,2),(9,16),(11,18),(6,13),(15,20),(17,12),(19,10)
\}$&\cr
&&118&&10&& $\{ (1,4),(3,8),
(14,5),(7,2),(16,9),(18,11),(6,13),(20,15),(12,17),(10,19)
\}$&\cr
&&119&&10&& $\{ (1,4),(3,8),
(5,14),(7,2),(9,16),(11,18),(13,20),(15,6),(17,10),(19,12)
\}$&\cr
&&120&&10&& $\{ (1,4),(3,8),
(14,5),(7,2),(16,9),(18,11),(20,13),(6,15),(10,17),(12,19)
\}$&\cr
&&121&&10&& $\{ (1,4),(3,8),
(5,14),(7,2),(9,16),(11,18),(13,20),(15,6),(17,12),(19,10)
\}$&\cr
&&122&&10&& $\{ (1,4),(3,8),
(14,5),(7,2),(16,9),(18,11),(20,13),(6,15),(10,17),(12,19)
\}$&\cr
&&123&&10&& $\{ (1,4),(3,8),
(5,14),(7,2),(9,18),(11,16),(13,6),(15,12),(17,20),(19,10)
\}$&\cr
&&124&&10&& $\{ (1,4),(3,8),
(5,14),(7,2),(9,18),(11,20),(13,16),(15,6),(17,12),(19,10)
\}$&\cr
&&125&&10&& $\{ (1,4),(3,8),
(5,14),(7,2),(9,20),(11,18),(13,16),(15,6),(17,12),(19,10)
\}$&\cr
\noalign{\hrule}}}}
\vskip 1 cm \centerline
{\vbox
{\offinterlineskip
\halign{\strut#&\vrule#&\hfil#\hfil&\vrule#&\hfil#\hfil&\
\vrule#&\hfil#\hfil&\vrule#\cr
\noalign{\hrule}
&&\omit {\bf Order} && {\bf Crossings}&&
{\bf Knot's Name} &\cr \noalign {\hrule}
&&126&&10&& $\{ (1,4),(3,8),(5,14),(7,10),(9,2),(11,18),(13,6),(15,20),
(17,12),(19,16)
\}$&\cr
&&127&&10&& $\{ (1,4),(3,8),(5,14),(7,12),(9,2),(11,16),(13,20),(15,18),
(17,10),(19,6)
\}$&\cr
&&128&&10&& $\{ (1,4),(3,8),(5,16),(7,10),(9,2),(11,18),(13,20),(15,6),
(17,14),(19,12)
\}$&\cr
&&129&&10&& $\{ (1,4),(3,8),(5,16),(7,10),(9,2),(11,20),(13,18),(15,6),
(17,14),(19,12)
\}$&\cr
&&130&&10&& $\{ (1,4),(3,8),(5,18),(7,12),(9,2),(11,16),(13,20),(15,6),
(17,10),(19,14)
\}$&\cr
&&131&&10&& $\{ (1,4),(3,8),(18,5),(12,7),(9,2),(16,11),(20,13),(6,15),
(10,17),(14,19)
\}$&\cr
&&132&&10&& $\{ (1,4),(3,10),(5,12),(7,14),(9,16),(11,2),(13,20),(15,18),
(17,8),(19,6)
\}$&\cr
&&133&&10&& $\{ (1,4),(3,10),(5,12),(7,14),(9,18),(11,2),(13,6),(15,20),
(17,8),(19,16)
\}$&\cr
&&134&&10&& $\{ (1,4),(3,10),(5,12),(7,14),(9,18),(11,2),(13,16),(15,6),
(17,20),(19,8)
\}$&\cr
&&135&&10&& $\{ (1,4),(3,10),(5,12),(7,16),(9,2),(11,8),(13,18),(15,20),
(17,6),(19,14)
\}$&\cr
&&136&&10&& $\{ (1,4),(3,10),(5,12),(7,16),(9,2),(11,8),(13,20),(15,18),
(17,6),(19,14)
\}$&\cr
&&137&&10&& $\{ (1,4),(3,10),(5,12),(7,16),(9,2),(11,20),(13,6),(15,18),
(17,8),(19,14)
\}$&\cr
&&138&&10&& $\{ (1,4),(3,10),(5,12),(7,16),(9,2),(11,20),(13,8),(15,18),
(17,6),(19,14)
\}$&\cr
&&139&&10&& $\{ (1,4),(3,10),(5,12),(7,16),(9,14),(11,2),(13,20),(15,18),
(17,8),(19,6)
\}$&\cr
&&140&&10&& $\{ (1,4),(3,10),(5,12),(7,16),(9,18),(11,2),(13,20),(15,6),
(17,8),(19,14)
\}$&\cr
&&141&&10&& $\{ (1,4),(3,10),(5,12),(7,16),(9,18),(11,2),(13,20),(15,8),
(17,6),(19,14)
\}$&\cr
&&142&&10&& $\{ (1,4),(3,10),(5,12),(7,16),(9,20),(11,2),(13,8),(15,18),
(17,6),(19,14)
\}$&\cr
&&143&&10&& $\{ (1,4),(3,10),(5,14),(7,12),(9,2),(11,16),(13,18),(15,20),
(17,8),(19,6)
\}$&\cr
&&144&&10&& $\{ (1,4),(3,10),(5,14),(12,7),(9,2),(16,11),(13,18),(15,20),
(8,17),(19,6)
\}$&\cr
&&145&&10&& $\{ (1,4),(3,10),(5,14),(7,12),(9,18),(11,2),(13,6),(15,20),
(17,8),(19,16)
\}$&\cr
&&146&&10&& $\{ (1,4),(3,10),(5,14),(7,12),(9,18),(11,2),(13,16),(15,6),
(17,20),(19,8)
\}$&\cr
&&147&&10&& $\{ (1,4),(3,10),(5,14),(7,16),(9,2),(11,8),(13,18),(15,20),
(17,12),(19,6)
\}$&\cr
&&148&&10&& $\{ (1,4),(3,10),(5,14),(7,16),(9,2),(11,18),(13,8),(15,6),
(17,20),(19,12)
\}$&\cr
&&149&&10&& $\{ (1,4),(3,10),(5,14),(7,16),(9,2),(18,11),(13,8),(15,6),
(20,17),(12,19)
\}$&\cr
&&150&&10&& $\{ (1,4),(3,10),(14,5),(16,7),(9,2),(11,18),(8,13),(6,15),
(17,20),(19,12)
\}$&\cr
&&151&&10&& $\{ (1,4),(3,10),(5,14),(7,16),(9,2),(11,18),(13,20),(15,6),
(17,8),(19,12)
\}$&\cr
&&152&&10&& $\{ (1,4),(3,10),(14,5),(16,7),(9,2),(18,11),(20,13),(6,15),
(8,17),(12,19)
\}$&\cr
&&153&&10&& $\{ (1,4),(3,10),(5,14),(7,16),(9,2),(11,18),(13,20),(15,8),
(17,6),(19,12)
\}$&\cr
&&154&&10&& $\{ (1,4),(3,10),(5,14),(7,16),(9,2),(18,11),(20,13),(15,8),
(17,6),(12,19)
\}$&\cr
&&155&&10&& $\{ (1,4),(3,10),(14,5),(16,7),(9,2),(11,18),(13,20),(8,15),
(6,17),(19,12)
\}$&\cr
&&156&&10&& $\{ (1,4),(3,10),(14,5),(16,7),(9,2),(18,11),(20,13),(8,15),
(6,17),(12,19)
\}$&\cr
&&157&&10&& $\{ (1,4),(3,10),(5,14),(7,16),(9,2),(11,20),(13,18),(15,6),
(17,8),(19,12)
\}$&\cr
&&158&&10&& $\{ (1,4),(3,10),(5,14),(7,16),(9,2),(11,20),(13,18),(15,8),
(17,6),(19,12)
\}$&\cr
&&159&&10&& $\{ (1,4),(3,10),(5,14),(7,18),(9,2),(11,16),(13,6),(15,20),
(17,8),(19,12)
\}$&\cr
&&160&&10&& $\{ (1,4),(3,10),(5,14),(7,18),(9,2),(11,16),(13,8),(15,20),
(17,12),(19,6)
\}$&\cr
&&161&&10&& $\{ (1,4),(3,10),(5,14),(7,18),(9,2),(11,16),(13,20),(15,8),
(17,12),(19,6)
\}$&\cr
&&162&&10&& $\{ (1,4),(3,10),(5,16),(7,12),(9,2),(11,8),(13,18),(15,20),
(17,6),(19,14)
\}$&\cr
&&163&&10&& $\{ (1,4),(3,10),(5,16),(7,12),(9,2),(11,8),(13,20),(15,18),
(17,6),(19,14)
\}$&\cr
&&164&&10&& $\{ (1,4),(3,10),(5,16),(7,14),(9,2),(11,8),(13,18),(15,20),
(17,12),(19,6)
\}$&\cr
&&165&&10&& $\{ (1,4),(3,10),(5,16),(7,14),(9,2),(11,18),(13,8),(15,6),
(17,20),(19,12)
\}$&\cr
&&166&&10&& $\{ (1,4),(3,10),(5,16),(7,14),(9,2),(11,20),(13,8),(15,18),
(17,6),(19,12)
\}$&\cr
&&167&&10&& $\{ (1,4),(3,10),(5,16),(7,18),(9,12),(11,2),(13,20),(15,8),
(17,6),(19,14)
\}$&\cr
&&168&&10&& $\{ (1,4),(3,10),(5,18),(7,16),(9,12),(11,2),(13,20),(15,8),
(17,6),(19,14)
\}$&\cr
&&169&&10&& $\{ (1,4),(3,12),(5,14),(7,16),(9,18),(11,2),(13,10),(15,20),
(17,8),(19,6)
\}$&\cr
&&170&&10&& $\{ (1,4),(3,12),(5,14),(7,16),(9,18),(11,2),(13,20),(15,6),
(17,8),(19,10)
\}$&\cr
&&171&&10&& $\{ (1,4),(3,12),(5,14),(7,16),(9,18),(11,2),(13,20),(15,6),
(17,10),(19,8)
\}$&\cr
&&172&&10&& $\{ (1,4),(3,12),(5,14),(7,16),(9,18),(11,2),(13,20),(15,10),
(17,8),(19,6)
\}$&\cr
&&173&&10&& $\{ (1,4),(3,12),(5,14),(7,16),(9,18),(11,20),(13,2),(15,6),
(17,8),(19,10)
\}$&\cr
&&174&&10&& $\{ (1,4),(3,12),(5,14),(7,16),(9,18),(11,20),(13,2),(15,6),
(17,10),(19,8)
\}$&\cr
&&175&&10&& $\{ (1,4),(3,12),(5,14),(7,16),(9,18),(11,20),(13,2),(15,10),
(17,8),(19,6)
\}$&\cr
&&176&&10&& $\{ (1,4),(3,12),(5,14),(7,16),(9,20),(11,18),(13,2),(15,8),
(17,6),(19,10)
\}$&\cr
&&177&&10&& $\{ (1,4),(3,12),(5,14),(7,18),(9,16),(11,2),(13,10),(15,20),
(17,8),(19,6)
\}$&\cr
&&178&&10&& $\{ (1,4),(3,12),(5,14),(7,18),(9,16),(11,2),(13,20),(15,10),
(17,8),(19,6)
\}$&\cr
&&179&&10&& $\{ (1,4),(3,12),(5,14),(7,18),(9,16),(11,20),(13,2),(15,10),
(17,8),(19,6)
\}$&\cr
&&180&&10&& $\{ (1,4),(3,12),(5,16),(7,14),(9,18),(11,2),(13,8),(15,20),
(17,10),(19,6)
\}$&\cr
&&181&&10&& $\{ (1,4),(3,12),(5,16),(14,7),(9,18),(11,2),(8,13),(15,20),
(17,10),(19,6)
\}$&\cr
&&182&&10&& $\{ (1,4),(3,12),(16,5),(7,14),(18,9),(11,2),(13,8),(20,15),
(10,17),(6,19)
\}$&\cr
&&183&&10&& $\{ (1,4),(3,12),(16,5),(14,7),(18,9),(11,2),(8,13),(20,15),
(10,17),(6,19)
\}$&\cr
&&184&&10&& $\{ (1,4),(3,12),(5,16),(7,14),(9,18),(11,2),(13,20),(15,6),
(17,10),(19,8)
\}$&\cr
&&185&&10&& $\{ (1,4),(3,12),(5,16),(7,18),(9,14),(11,2),(13,10),(15,20),
(17,6),(19,8)
\}$&\cr
&&186&&10&& $\{ (1,4),(3,12),(5,16),(7,18),(9,14),(11,2),(13,10),(15,20),
(17,8),(19,6)
\}$&\cr
&&187&&10&& $\{ (1,4),(3,12),(5,16),(7,18),(9,20),(11,14),(13,2),(15,10),
(17,6),(19,8)
\}$&\cr
&&188&&10&& $\{ (1,4),(3,12),(5,16),(7,18),(9,20),(11,14),(13,2),(15,10),
(17,8),(19,6)
\}$&\cr
\noalign{\hrule}}}}
\vskip 1 cm \centerline
{\vbox
{\offinterlineskip
\halign{\strut#&\vrule#&\hfil#\hfil&\vrule#&\hfil#\hfil&\
\vrule#&\hfil#\hfil&\vrule#\cr
\noalign{\hrule}
&&\omit {\bf Order} && {\bf Crossings}&&
{\bf Knot's Name} &\cr \noalign {\hrule}
&&189&&10&& $\{ (1,4),(3,12),(5,16),(7,20),(9,18),(11,2),(13,8),(15,6),
(17,10),(19,14)
\}$&\cr
&&190&&10&& $\{ (1,4),(3,12),(5,18),(7,20),(9,14),(11,16),(13,2),(15,10),
(17,8),(19,6)
\}$&\cr
&&191&&10&& $\{ (1,4),(3,12),(5,18),(7,20),(9,16),(11,14),(13,2),(15,10),
(17,8),(19,6)
\}$&\cr
&&192&&10&& $\{ (1,4),(3,12),(5,20),(7,18),(9,16),(11,14),(13,2),(15,10),
(17,8),(19,6)
\}$&\cr
&&193&&10&& $\{ (1,6),(3,8),(5,10),(7,14),(9,16),(11,18),(13,20),(15,2),
(17,4),(19,12)
\}$&\cr
&&194&&10&& $\{ (1,6),(3,8),(5,10),(7,14),(9,16),(18,11),(20,13),(15,2),
(17,4),(12,19)
\}$&\cr
&&195&&10&& $\{ (1,6),(3,8),(5,10),(7,14),(9,16),(11,20),(13,18),(15,2),
(17,4),(19,12)
\}$&\cr
&&196&&10&& $\{ (1,6),(3,8),(5,10),(7,14),(9,16),(20,11),(18,13),(15,2),
(17,4),(12,19)
\}$&\cr
&&197&&10&& $\{ (1,6),(3,8),(5,12),(7,2),(9,16),(11,4),(13,18),(15,20),
(17,10),(19,14)
\}$&\cr
&&198&&10&& $\{ (1,6),(3,8),(5,12),(7,2),(16,9),(11,4),(18,13),(20,15),
(10,17),(14,19)
\}$&\cr
&&199&&10&& $\{ (1,6),(3,8),(5,14),(7,2),(9,16),(11,18),(13,4),(15,20),
(17,10),(19,12)
\}$&\cr
&&200&&10&& $\{ (1,6),(3,8),(5,14),(7,2),(9,16),(11,18),(13,4),(15,20),
(17,12),(19,10)
\}$&\cr
&&201&&10&& $\{ (1,6),(3,8),(5,14),(7,2),(9,16),(11,18),(13,20),(15,4),
(17,10),(19,12)
\}$&\cr
&&202&&10&& $\{ (1,6),(3,8),(5,14),(7,2),(9,16),(11,18),(13,20),(15,4),
(17,12),(19,10)
\}$&\cr
&&203&&10&& $\{ (1,6),(3,8),(5,14),(7,16),(9,4),(11,18),(13,20),(15,2),
(17,10),(19,12)
\}$&\cr
&&204&&10&& $\{ (1,6),(3,8),(5,14),(7,16),(9,4),(11,18),(13,20),(15,2),
(17,12),(19,10)
\}$&\cr
&&205&&10&& $\{ (1,6),(3,8),(5,14),(7,18),(9,16),(11,4),(13,20),(15,10),
(17,2),(19,12)
\}$&\cr
&&206&&10&& $\{ (1,6),(3,8),(5,14),(7,18),(9,16),(11,4),(20,13),(15,10),
(17,2),(12,19)
\}$&\cr
&&207&&10&& $\{ (1,6),(3,8),(5,16),(7,14),(9,4),(11,18),(13,20),(15,2),
(17,10),(19,12)
\}$&\cr
&&208&&10&& $\{ (1,6),(3,8),(5,16),(7,14),(9,4),(11,18),(13,20),(15,2),
(17,12),(19,10)
\}$&\cr
&&209&&10&& $\{ (1,6),(3,8),(5,18),(7,14),(9,16),(11,4),(13,20),(15,2),
(17,10),(19,12)
\}$&\cr
&&210&&10&& $\{ (1,6),(3,10),(5,12),(7,14),(9,18),(11,16),(13,20),(15,2),
(17,4),(19,8)
\}$&\cr
&&211&&10&& $\{ (1,6),(10,3),(5,12),(7,14),(9,18),(16,11),(13,20),(15,2),
(4,17),(19,8)
\}$&\cr
&&212&&10&& $\{ (1,6),(3,10),(5,12),(7,20),(9,18),(11,16),(13,8),(15,2),
(17,4),(19,14)
\}$&\cr
&&213&&10&& $\{ (1,6),(3,10),(5,14),(7,2),(9,16),(11,18),(13,20),(15,8),
(17,4),(19,12)
\}$&\cr
&&214&&10&& $\{ (1,6),(3,10),(5,14),(7,2),(9,16),(11,20),(13,18),(15,8),
(17,4),(19,12)
\}$&\cr
&&215&&10&& $\{ (1,6),(3,10),(5,14),(7,16),(9,2),(11,18),(13,4),(15,20),
(17,8),(19,12)
\}$&\cr
&&216&&10&& $\{ (1,6),(3,10),(5,14),(7,16),(9,4),(11,18),(13,2),(15,20),
(17,12),(19,8)
\}$&\cr
&&217&&10&& $\{ (1,6),(3,10),(5,14),(7,16),(9,18),(11,4),(13,20),(15,2),
(17,8),(19,12)
\}$&\cr
&&218&&10&& $\{ (1,6),(3,10),(5,14),(7,16),(9,18),(11,4),(20,13),(15,2),
(17,8),(12,19)
\}$&\cr
&&219&&10&& $\{ (1,6),(3,10),(5,14),(7,18),(9,2),(11,16),(13,20),(15,4),
(17,8),(19,12)
\}$&\cr
&&220&&10&& $\{ (1,6),(3,10),(5,14),(7,18),(9,16),(11,4),(13,20),(15,2),
(17,8),(19,12)
\}$&\cr
&&221&&10&& $\{ (1,6),(3,10),(5,14),(7,18),(9,16),(11,4),(20,13),(15,2),
(17,8),(12,19)
\}$&\cr
&&222&&10&& $\{ (1,6),(3,10),(5,16),(7,14),(9,2),(11,18),(13,8),(15,20),
(17,4),(19,12)
\}$&\cr
&&223&&10&& $\{ (1,6),(10,3),(16,5),(7,14),(2,9),(18,11),(13,8),(15,20),
(4,17),(12,19)
\}$&\cr
&&224&&10&& $\{ (1,6),(3,10),(5,16),(7,14),(9,18),(11,4),(13,20),(15,2),
(17,12),(19,8)
\}$&\cr
&&225&&10&& $\{ (1,6),(3,10),(5,16),(7,20),(9,14),(11,2),(13,18),(15,4),
(17,8),(19,12)
\}$&\cr
&&226&&10&& $\{ (1,6),(3,10),(5,16),(7,20),(9,14),(11,4),(13,18),(15,2),
(17,12),(19,8)
\}$&\cr
&&227&&10&& $\{ (1,6),(3,10),(5,18),(7,12),(9,4),(11,16),(13,20),(15,8),
(17,2),(19,14)
\}$&\cr
&&228&&10&& $\{ (1,6),(3,10),(5,18),(7,14),(9,2),(11,16),(13,20),(15,8),
(17,4),(19,12)
\}$&\cr
&&229&&10&& $\{ (1,6),(10,3),(18,5),(7,14),(2,9),(16,11),(13,20),(15,8),
(4,17),(19,12)
\}$&\cr
&&230&&10&& $\{ (1,6),(3,10),(5,18),(7,14),(9,16),(11,4),(13,20),(15,8),
(17,2),(19,12)
\}$&\cr
&&231&&10&& $\{ (1,6),(3,10),(5,18),(7,16),(9,14),(11,4),(13,20),(15,8),
(17,2),(19,12)
\}$&\cr
&&232&&10&& $\{ (1,6),(3,10),(5,20),(7,14),(9,16),(11,18),(13,4),(15,8),
(17,2),(19,12)
\}$&\cr
&&233&&10&& $\{ (1,6),(3,12),(5,14),(7,16),(9,18),(11,2),(13,4),(15,20),
(17,8),(19,10)
\}$&\cr
&&234&&10&& $\{ (1,6),(3,12),(5,14),(7,16),(9,18),(11,2),(13,4),(15,20),
(17,10),(19,8)
\}$&\cr
&&235&&10&& $\{ (1,6),(3,12),(5,14),(7,16),(9,18),(11,4),(13,2),(15,20),
(17,10),(19,8)
\}$&\cr
&&236&&10&& $\{ (1,6),(3,12),(5,14),(7,16),(9,18),(11,20),(13,4),(15,2),
(17,8),(19,10)
\}$&\cr
&&237&&10&& $\{ (1,6),(3,12),(5,14),(7,16),(9,18),(11,20),(13,4),(15,2),
(17,10),(19,8)
\}$&\cr
&&238&&10&& $\{ (1,6),(3,12),(5,14),(7,18),(9,20),(11,16),(13,4),(15,2),
(17,10),(19,8)
\}$&\cr
&&239&&10&& $\{ (1,6),(3,12),(5,14),(7,20),(9,18),(11,16),(13,4),(15,2),
(17,10),(19,8)
\}$&\cr
&&240&&10&& $\{ (1,6),(3,14),(5,12),(7,16),(9,18),(11,20),(13,4),(15,2),
(17,8),(19,10)
\}$&\cr
&&241&&10&& $\{ (1,6),(3,14),(5,12),(7,16),(9,18),(11,20),(13,4),(15,2),
(17,10),(19,8)
\}$&\cr
&&242&&10&& $\{ (1,6),(3,14),(5,12),(7,18),(9,20),(11,16),(13,4),(15,2),
(17,10),(19,8)
\}$&\cr
&&243&&10&& $\{ (1,6),(3,14),(5,12),(7,20),(9,18),(11,16),(13,4),(15,2),
(17,10),(19,8)
\}$&\cr
&&244&&10&& $\{ (1,6),(3,16),(5,12),(7,14),(9,18),(11,4),(13,20),(15,2),
(17,8),(19,10)
\}$&\cr
&&245&&10&& $\{ (1,6),(3,16),(5,12),(7,14),(9,18),(11,4),(13,20),(15,2),
(17,10),(19,8)
\}$&\cr
&&246&&10&& $\{ (1,6),(3,16),(5,18),(7,14),(9,2),(11,4),(13,20),(15,8),
(17,10),(19,12)
\}$&\cr
&&247&&10&& $\{ (1,8),(3,10),(5,12),(7,14),(9,16),(11,18),(13,20),(15,2),
(17,4),(19,6)
\}$&\cr
&&248&&10&& $\{ (1,8),(3,10),(5,14),(7,16),(9,2),(11,18),(13,20),(15,6),
(17,4),(19,12)
\}$&\cr
&&249&&10&& $\{ (1,8),(3,10),(5,16),(7,14),(9,2),(11,18),(13,20),(15,6),
(17,4),(19,12)
\}$&\cr
\noalign{\hrule}}}}
\par\medskip
\centerline
{\vbox
{\offinterlineskip
\halign{\strut#&\vrule#&\hfil#\hfil&\vrule#&\hfil#\hfil&\vrule#&\hfil#\hfil&
\vrule#&\hfil#\hfil&\vrule#&\hfil#\hfil&\vrule#&\hfil#\hfil&
\vrule#&\hfil#\hfil&\vrule#&\hfil#\hfil&\vrule#&\hfil#\hfil&
\vrule#\cr
\noalign{\hrule}
&&{\bf Order}&& {\bf Knot Characteristic}
&&[3,1]&&[5,1]&&[5,2]&&[7,1]&&[7,2]&&[7,3]&\cr
\noalign {\hrule}
&&\*\*0&&$1$&&0&&0&&0&&0&&0&&0&\cr
&&\*\*1&&$(s^3+3s^2+3s+2)(s+1)^{-3}$&&1&&0&&0&&0&&1&&0&\cr
&&\*\*2&&$(s^4+4s^3+5s^2+4s+1)(s+1)^{-4}$&&0&&1&&0&&0&&0&&0&\cr
&&\*\*3&&$(s^4+4s^3+5s^2+6s)(s+1)^{-4}$&&0&&0&&0&&1&&0&&0&\cr
&&\*\*4&&$(s^5+5s^4+10s^3+15s^2+10s+3)(s+1)^{-5}$&&0&&1&&0&&0&&0&&0&\cr
&&\*\*5&&$(2s^4+4s^3+8s^2+4s+2)(s+1)^{-4}$&&0&&0&&0&&0&&0&&1&\cr
&&\*\*6&&$(s^6+6s^5+13s^4+20s^3+17s^2+9s+2)(s+1)^{-6}$&&0&&0&&0&&0&&0&&0&\cr
&&\*\*7&&$(s^5+5s^4+7s^3+9s^2+2s)(s+1)^{-5}$&&1&&0&&1&&0&&0&&1&\cr
&&\*\*8&&$(4s^6+16s^5+31s^4+30s^3+19s^2+7s+1)(s+1)^{-7}$&&1&&0&&0&&1&&0&&0&\cr
&&\*\*9&&$(4s^4+10s^3+12s^2+9s)(s+1)^{-5}$&&0&&0&&1&&0&&0&&0&\cr
&&\*10&&$(s^6+6s^5+13s^4+24s^3+18s^2+8s)(s+1)^{-6}$&&0&&0&&0&&0&&0&&1&\cr
&&\*11&&$(s^5+5s^4+7s^3+12s^2)(s+1)^{-5}$&&0&&0&&1&&0&&0&&0&\cr
&&\*12&&$(s^6+6s^5+14s^4+29s^3+22s^2+8s)(s+1)^{-6}$&&0&&0&&0&&0&&0&&0&\cr
&&\*13&&$(s^5+5s^4+8s^3+16s^2)(s+1)^{-5}$&&1&&1&&0&&0&&0&&0&\cr
&&\*14&&$(s^7+7s^6+21s^5+49s^4+63s^3+49s^2+21s+4)(s+1)^{-7}$
&&0&&0&&0&&1&&0&&0&\cr
&&\*15&&$(s^7+6s^6+15s^5+29s^4+30s^3+23s^2+6s)(s+1)^{-7}$&&0&&0&&0&&0&&0&&0&\cr
&&\*16&&$(s^6+6s^5+13s^4+24s^3+19s^2+12s)(s+1)^{-6}$&&1&&0&&1&&0&&1&&0&\cr
&&\*17&&$(2s^7+9s^6+24s^5+37s^4+39s^3+24s^2+9s)(s+1)^{-7}$
&&1&&0&&0&&0&&1&&0&\cr
&&\*18&&$(s^7+7s^6+20s^5+36s^4+38s^3+28s^2+11s+3)(s+1)^{-7}$
&&1&&1&&0&&0&&1&&0&\cr
&&\*19&&$(3s^7+14s^6+35s^5+48s^4+44s^3+25s^2+9s+2)(s+1)^{-7}$
&&1&&0&&0&&0&&1&&0&\cr
&&\*20&&$(s^8+8s^7+26s^6+60s^5+89s^4+88s^3+56s^2+20s+3)(s+1)^{-8}$&&1&&0&&0&&0
&&1&&0&\cr
&&\*21&&$(6s^5+13s^4+21s^3+16s^2+7s+2)(s+1)^{-6}$&&0&&1&&1&&0&&0&&0&\cr
&&\*22&&$(2s^5+11s^4+20s^3+18s^2+4s)(s+1)^{-6}$&&0&&0&&0&&0&&0&&0&\cr
&&\*23&&$(3s^7+14s^6+35s^5+51s^4+49s^3+27s^2+9s+2)(s+1)^{-7}$
&&0&&0&&0&&0&&0&&0&\cr
&&\*24&&$(6s^5+13s^4+24s^3+18s^2+7s+2)(s+1)^{-6}$&&0&&0&&0&&0&&0&&0&\cr
&&\*25&&$(s^8+8s^7+25s^6+56s^5+80s^4+81s^3+52s^2+19s+3)(s+1)^{-8}$
&&0&&0&&0&&0&&0&&0&\cr
&&\*26&&$(s^7+7s^6+17s^5+33s^4+33s^3+23s^2+6s)(s+1)^{-7}$&&1&&0&&1&&0&&1&&1&\cr
&&\*27&&$(s^7+7s^6+16s^5+30s^4+27s^3+16s^2+3s)(s+1)^{-7}$&&0&&0&&0&&0&&0&&0&\cr
&&\*28&&$(s^6+6s^5+9s^4+16s^3+3s^2)(s+1)^{-6}$&&0&&0&&0&&0&&0&&0&\cr
&&\*29&&$(8s^5+16s^4+24s^3+16s^2+8s)(s+1)^{-6}$&&4&&1&&0&&0&&8&&0&\cr
&&\*30&&$(s^8+8s^7+24s^6+52s^5+71s^4+68s^3+42s^2+16s+3)(s+1)^{-8}$
&&1&&0&&0&&1&&1&&0&\cr
&&\*31&&$(2s^8+13s^7+36s^6+64s^5+75s^4+64s^3+36s^2+13s+2)(s+1)^{-8}$
&&0&&0&&0&&0&&0&&0&\cr
&&\*32&&$(3s^7+11s^6+30s^5+46s^4+47s^3+27s^2+9s+2)(s+1)^{-7}$
&&0&&1&&0&&1&&0&&0&\cr
&&\*33&&$(2s^8+13s^7+36s^6+64s^5+75s^4+64s^3+36s^2+13s+2)(s+1)^{-8}$
&&0&&1&&0&&0&&0&&1&\cr
&&\*34&&$(2s^7+11s^6+23s^5+32s^4+25s^3+15s^2+2s)(s+1)^{-7}$
&&0&&0&&0&&0&&0&&0&\cr
&&\*35&&$(2s^5+13s^4+10s^3+13s^2+2s)(s+1)^{-6}$&&0&&0&&0&&0&&0&&0&\cr
&&\*36&&$(s^8+8s^7+27s^6+56s^5+74s^4+65s^3+35s^2+12s+2)(s+1)^{-8}$
&&0&&0&&0&&0&&0&&0&\cr
&&\*37&&$(s^9+7s^8+23s^7+58s^6+98s^5+123s^4+105s^3+55s^2+16s+2)
(s+1)^{-9}$&&0&&0&&1&&0&&0&&0&\cr
&&\*38&&$\*\*\*\*(s^9+11s^8+46s^7+104s^6+146s^5+134s^4+85s^3+36s^2+9s+1)
(s+1)^{-9}\*\*\*\*$&&0&&0&&0&&0&&0&&1&\cr
&&\*39&&$(s^9+9s^8+31s^7+67s^6+98s^5+110s^4+88s^3+47s^2+15s+2)(s+1)^{-9}$
&&0&&0&&1&&0&&0&&0&\cr
&&\*40&&$(4s^8+22s^7+62s^6+100s^5+114s^4+89s^3+47s^2+15s+2)(s+1)^{-9}$
&&0&&0&&0&&0&&0&&1&\cr
&&\*41&&$(4s^8+29s^7+85s^6+144s^5+151s^4+99s^3+39s^2+9s+1)(s+1)^{-9}$
&&0&&0&&0&&1&&0&&0&\cr
&&\*42&&$(4s^7+18s^6+40s^5+44s^4+37s^3+15s^2+2s)(s+1)^{-8}$
&&0&&0&&0&&0&&0&&0&\cr
&&\*43&&$(6s^5+18s^4+22s^3+20s^2)(s+1)^{-6}$&&0&&0&&0&&0&&0&&0&\cr
&&\*44&&$(2s^7+9s^6+25s^5+42s^4+43s^3+30s^2+11s+3)(s+1)^{-7}$
&&1&&0&&1&&0&&1&&0&\cr
&&\*45&&$(6s^8+32s^7+81s^6+120s^5+118s^4+75s^3+33s^2+9s+1)(s+1)^{-9}$
&&0&&0&&1&&0&&0&&0&\cr
&&\*46&&$(6s^5+15s^4+27s^3+21s^2+12s)(s+1)^{-6}$&&1&&1&&0&&0&&1&&0&\cr
&&\*47&&$(9s^3+9s^2+12s)(s+1)^{-5}$&&0&&0&&0&&0&&0&&0&\cr
&&\*48&&$(s^8+10s^7+37s^6+70s^5+80s^4+56s^3+27s^2+8s+1)(s+1)^{-8}$
&&0&&0&&0&&0&&0&&0&\cr
&&\*49&&$(2s^6+14s^5+32s^4+30s^3+19s^2+7s+1)(s+1)^{-7}$&&1&&0&&0&&0&&0&&0&\cr
&&\*50&&$(s^8+9s^7+29s^6+57s^5+71s^4+61s^3+33s^2+12s+2)(s+1)^{-8}$
&&0&&0&&0&&0&&0&&0&\cr
&&\*51&&$(2s^8+12s^7+35s^6+65s^5+76s^4+60s^3+30s^2+9s+1)(s+1)^{-8}$
&&0&&0&&1&&1&&0&&0&\cr
&&\*52&&$(2s^7+9s^6+25s^5+42s^4+45s^3+33s^2+11s+3)(s+1)^{-7}$
&&0&&1&&0&&0&&0&&0&\cr
&&\*53&&$(4s^8+24s^7+69s^6+111s^5+125s^4+96s^3+49s^2+15s+2)(s+1)^{-9}$
&&1&&0&&0&&0&&0&&0&\cr
&&\*54&&$(6s^8+32s^7+81s^6+120s^5+118s^4+75s^3+33s^2+9s+1)(s+1)^{-9}$
&&0&&0&&0&&0&&0&&0&\cr
&&\*55&&$(s^7+7s^6+16s^5+35s^4+29s^3+20s^2)(s+1)^{-7}$&&1&&1&&0&&0&&1&&0&\cr
&&\*56&&$(9s^6+30s^5+54s^4+44s^3+24s^2+8s+1)(s+1)^{-8}$&&0&&0&&0&&0&&0&&0&\cr
&&\*57&&$(4s^7+20s^6+45s^5+50s^4+42s^3+17s^2+2s)(s+1)^{-8}$
&&4&&1&&1&&0&&0&&1&\cr
&&\*58&&$(3s^7+20s^6+47s^5+66s^4+53s^3+27s^2+8s+1)(s+1)^{-8}$
&&4&&0&&1&&0&&0&&1&\cr
&&\*59&&$(10s^7+51s^6+99s^5+106s^4+76s^3+35s^2+9s+1)(s+1)^{-9}$
&&4&&0&&0&&0&&0&&0&\cr
&&\*60&&$(6s^6+24s^5+44s^4+54s^3+34s^2+12s)(s+1)^{-7}$&&0&&0&&0&&0&&0&&0&\cr
&&\*61&&$(s^8+10s^7+38s^6+74s^5+84s^4+57s^3+27s^2+8s+1)(s+1)^{-8}$
&&1&&0&&0&&0&&0&&1&\cr
&&\*62&&$(2s^6+15s^5+35s^4+31s^3+19s^2+7s+1)(s+1)^{-7}$&&0&&0&&0&&0&&0&&0&\cr
\noalign{\hrule}}}}
\par
\centerline
{\vbox
{\offinterlineskip
\halign{\strut#&\vrule#&\hfil#\hfil&\vrule#&\hfil#\hfil&\vrule#&\hfil#\hfil&
\vrule#&\hfil#\hfil&\vrule#&\hfil#\hfil&\vrule#&\hfil#\hfil&
\vrule#&\hfil#\hfil&\vrule#&\hfil#\hfil&\vrule#&\hfil#\hfil&
\vrule#\cr
\noalign{\hrule}
&&{\bf Order}&& {\bf Knot Characteristic}
&&[3,1]&&[5,1]&&[5,2]&&[7,1]&&[7,2]&&[7,3]&\cr
\noalign {\hrule}
&&\*63&&$(6s^5+18s^4+22s^3+20s^2)(s+1)^{-6}$&&0&&1&&0&&1&&0&&0&\cr
&&\*64&&$(s^8+8s^7+25s^6+65s^5+90s^4+83s^3+42s^2+10s)(s+1)^{-8}$
&&1&&0&&0&&0&&1&&0&\cr
&&\*65&&$(s^7+7s^6+17s^5+41s^4+34s^3+20s^2)(s+1)^{-7}$&&0&&0&&1&&0&&0&&0&\cr
&&\*66&&$(s^8+8s^7+24s^6+58s^5+78s^4+72s^3+37s^2+10s)(s+1)^{-8}$
&&1&&0&&0&&0&&1&&0&\cr
&&\*67&&$(s^7+7s^6+16s^5+35s^4+28s^3+15s^2)(s+1)^{-7}$&&0&&0&&0&&0&&0&&0&\cr
&&\*68&&$(s^6+6s^5+9s^4+20s^3)(s+1)^{-6}$&&1&&1&&0&&0&&0&&0&\cr
&&\*69&&$(s^8+8s^7+28s^6+59s^5+77s^4+63s^3+28s^2+6s)(s+1)^{-8}$
&&1&&0&&0&&0&&0&&0&\cr
&&\*70&&$(s^9+9s^8+33s^7+73s^6+117s^5+123s^4+84s^3+36s^2+9s+1)(s+1)^{-9}$
&&4&&0&&0&&0&&0&&0&\cr
&&\*71&&$(12s^6+34s^5+59s^4+46s^3+25s^2+8s+1)(s+1)^{-8}$&&0&&0&&0&&8&&0&&1&\cr
&&\*72&&$(s^7+7s^6+18s^5+43s^4+40s^3+21s^2+3s)(s+1)^{-7}$&&0&&6&&0&&0&&0&&0&\cr
&&\*73&&$(s^7+7s^6+17s^5+41s^4+35s^3+25s^2)(s+1)^{-7}$&&1&&0&&0&&0&&1&&1&\cr
&&\*74&&$(9s^7+36s^6+62s^5+69s^4+47s^3+22s^2+6s+1)(s+1)^{-8}$
&&1&&0&&1&&0&&1&&0&\cr
&&\*75&&$(16s^6+52s^5+52s^4+41s^3+18s^2+6s+1)(s+1)^{-8}$&&0&&1&&1&&0&&0&&0&\cr
&&\*76&&$(s^8+8s^7+26s^6+72s^5+101s^4+88s^3+42s^2+10s)(s+1)^{-8}$
&&0&&0&&0&&0&&0&&0&\cr
&&\*77&&$(s^7+7s^6+18s^5+47s^4+39s^3+20s^2)(s+1)^{-7}$&&0&&0&&0&&0&&0&&0&\cr
&&\*78&&$(s^7+7s^6+18s^5+47s^4+38s^3+15s^2)(s+1)^{-7}$&&1&&0&&0&&1&&0&&1&\cr
&&\*79&&$(s^6+6s^5+11s^4+30s^3)(s+1)^{-6}$&&0&&0&&0&&0&&0&&0&\cr
&&\*80&&$(2s^7+11s^6+33s^5+48s^4+45s^3+24s^2+7s+1)(s+1)^{-7}$
&&1&&6&&0&&0&&1&&0&\cr
&&\*81&&$(s^8+8s^7+27s^6+79s^5+113s^4+99s^3+47s^2+10s)(s+1)^{-8}$
&&0&&0&&0&&0&&0&&0&\cr
&&\*82&&$(s^7+7s^6+19s^5+53s^4+44s^3+20s^2)(s+1)^{-7}$&&1&&0&&0&&0&&0&&0&\cr
&&\*83&&$(s^6+6s^5+12s^4+35s^3)(s+1)^{-6}$&&4&&0&&0&&0&&0&&0&\cr
&&\*84&&$(s^9+9s^8+36s^7+114s^6+216s^5+270s^4+222s^3+117s^2+36s+5)(s+1)^{-9}$
&&1&&0&&0&&0&&1&&0&\cr
&&\*85&&$(9s^{8}+54s^{7}+134s^{6}+196s^{5}+175s^{4}+104s^{3}+42s^{2}+10s+1)
(s+1)^{-10}$
&&0&&1&&1&&0&&0&&0&\cr
&&\*86&&$(6s^{8}+43s^{7}+129s^{6}+216s^{5}+213s^{4}+124s^{3}+39s^{2}+5s)
(s+1)^{-10}$
&&0&&1&&0&&0&&0&&0&\cr
&&\*87&&$(2s^{9}+21s^{8}+82s^{7}+168s^{6}+220s^{5}+193s^{4}+114s^{3}+44s^{2}
+10s+1)(s+1)^{-10}$
&&0&&1&&0&&1&&0&&0&\cr
&&\*88&&$(2s^{7}+14s^{6}+39s^{5}+67s^{4}+64s^{3}+39s^{2}+9s)(s+1)^{-8}$
&&1&&1&&0&&0&&1&&0&\cr
&&\*89&&$(3s^{9}+26s^{8}+86s^{7}+152s^{6}+179s^{5}+145s^{4}+84s^{3}+35s^{2}+9s
+1)(s+1)^{-10}$
&&1&&1&&0&&0&&1&&0&\cr
&&\*90&&$(4s^{8}+24s^{7}+68s^{6}+106s^{5}+109s^{4}+75s^{3}+34s^{2}+9s+1)
(s+1)^{-9}$
&&0&&0&&0&&0&&0&&0&\cr
&&\*91&&$(s^{9}+7s^{8}+24s^{7}+61s^{6}+99s^{5}+122s^{4}+92s^{3}+42s^{2}+8s)
(s+1)^{-9}$
&&0&&0&&0&&0&&0&&0&\cr
&&\*92&&$(s^{8}+6s^{7}+17s^{6}+38s^{5}+45s^{4}+43s^{3}+12s^{2})(s+1)^{-8}$
&&1&&0&&0&&0&&0&&0&\cr
&&\*93&&$(8s^{7}+28s^{6}+64s^{5}+74s^{4}+60s^{3}+30s^{2}+10s+2)(s+1)^{-8}$
&&0&&0&&0&&0&&0&&0&\cr
&&\*94&&$(s^{9}+10s^{8}+39s^{7}+90s^{6}+131s^{5}+131s^{4}+83s^{3}+32s^{2}+6s)
(s+1)^{-9}$
&&0&&0&&1&&0&&0&&0&\cr
&&\*95&&$(6s^{8}+29s^{7}+75s^{6}+122s^{5}+132s^{4}+95s^{3}+45s^{2}+12s+1)
(s+1)^{-9}$
&&0&&0&&0&&0&&0&&0&\cr
&&\*96&&$(s^{9}+9s^{8}+31s^{7}+79s^{6}+128s^{5}+144s^{4}+106s^{3}+44s^{2}+8s)
(s+1)^{-9}$
&&0&&0&&1&&0&&0&&0&\cr
&&\*97&&$(9s^{5}+15s^{4}+27s^{3}+15s^{2}+9s)(s+1)^{-6}$
&&0&&1&&0&&0&&0&&0&\cr
&&\*98&&$(6s^{7}+27s^{6}+56s^{5}+76s^{4}+61s^{3}+29s^{2}+6s)(s+1)^{-8}$
&&0&&0&&0&&0&&0&&1&\cr
&&\*99&&$(2s^{8}+10s^{7}+29s^{6}+55s^{5}+68s^{4}+58s^{3}+28s^{2}+8s+1)
(s+1)^{-8}$
&&0&&0&&1&&0&&0&&0&\cr
&&100&&$(s^9+9s^8+29s^7+68s^6+101s^5+105s^4+70s^3+24s^2+3s)(s+1)^{-9}$
&&0&&0&&0&&0&&0&&0&\cr
&&101&&$(s^8+8s^7+24s^6+58s^5+80s^4+84s^3+48s^2+15s)(s+1)^{-8}$
&&0&&0&&0&&0&&0&&0&\cr
&&102&&$(3s^9+19s^8+59s^7+114s^6+152s^5+139s^4+89s^3+36s^2+9s)(s+1)^{-9}$
&&0&&0&&1&&0&&0&&1&\cr
&&103&&$(s^9+8s^8+30s^7+75s^6+124s^5+150s^4+121s^3+67s^2+22s+4)(s+1)^{-9}$
&&0&&0&&0&&0&&0&&0&\cr
&&104&&$(8s^9+46s^8+129s^7+222s^6+254s^5+205s^4+118s^3+48s^2+13s+2)(s+1)^{-10}$
&&0&&0&&0&&0&&0&&0&\cr
&&105&&$(s^7+7s^6+16s^5+35s^4+30s^3+25s^2)(s+1)^{-7}$
&&0&&0&&0&&0&&0&&0&\cr
&&106&&$(12s^7+51s^6+91s^5+100s^4+72s^3+35s^2+11s+2)(s+1)^{-9}$
&&0&&1&&0&&0&&0&&0&\cr
&&107&&$(9s^5+15s^4+27s^3+12s^2+6s)(s+1)^{-6}$
&&0&&0&&0&&1&&0&&0&\cr
&&108&&$(4s^9+28s^8+87s^7+159s^6+189s^5+159s^4+87s^3+28s^2+4s)(s+1)^{-10}$
&&0&&0&&0&&0&&0&&1&\cr
&&109&&$(s^8+6s^7+18s^6+43s^5+54s^4+60s^3+16s^2)(s+1)^{-8}$
&&1&&0&&0&&0&&0&&0&\cr
&&110&&$(10s^7+52s^6+109s^5+122s^4+76s^3+32s^2+9s+1)(s+1)^{-9}$
&&1&&0&&1&&0&&0&&0&\cr
&&111&&$(s^8+8s^7+25s^6+65s^5+92s^4+95s^3+53s^2+15s)(s+1)^{-8}$
&&0&&0&&0&&0&&0&&0&\cr
&&112&&$(3s^9+19s^8+62s^7+127s^6+179s^5+173s^4+116s^3+49s^2+12s)(s+1)^{-9}$
&&0&&0&&0&&0&&0&&0&\cr
&&113&&$(s^9+9s^8+34s^7+85s^6+137s^5+160s^4+125s^3+67s^2+22s+4)(s+1)^{-9}$
&&0&&0&&0&&0&&0&&1&\cr
&&114&&$(8s^9+50s^8+146s^7+255s^6+292s^5+232s^4+129s^3+50s^2+13s+2)(s+1)^{-10}$
&&0&&0&&0&&0&&0&&0&\cr
&&115&&$(s^7+7s^6+17s^5+41s^4+36s^3+30s^2)(s+1)^{-7}$
&&0&&0&&0&&0&&0&&0&\cr
&&116&&$(6s^7+23s^6+56s^5+76s^4+70s^3+37s^2+12s)(s+1)^{-8}$
&&0&&1&&1&&0&&0&&0&\cr
&&117&&$(s^8+8s^7+25s^6+53s^5+63s^4+56s^3+24s^2+8s)(s+1)^{-8}$
&&0&&0&&0&&0&&0&&0&\cr
&&118&&$(16s^7+64s^6+111s^5+118s^4+81s^3+37s^2+11s+2)(s+1)^{-9}$
&&0&&0&&0&&0&&0&&1&\cr
&&119&&$(4s^9+26s^8+84s^7+164s^6+212s^5+182s^4+110s^3+44s^2+12s+2)(s+1)^{-9}$
&&0&&0&&0&&0&&0&&0&\cr
&&120&&$(s^{10}+10s^9+42s^8+128s^7+263s^6+382s^5+393s^4+277s^3+127s^2+34s+4)
(s+1)^{-10}$&&0&&0&&0&&0&&0&&0&\cr
&&121&&$(8s^7+32s^6+74s^5+86s^4+68s^3+32s^2+10s+2)(s+1)^{-8}$
&&0&&0&&1&&0&&0&&0&\cr
&&122&&$(s^9+9s^8+32s^7+87s^6+146s^5+165s^4+120s^3+48s^2+8s)(s+1)^{-9}$
&&0&&0&&0&&0&&0&&1&\cr
&&123&&$(4s^6+16s^5+32s^4+44s^3+32s^2+16s)(s+1)^{-7}$
&&1&&0&&0&&0&&1&&1&\cr
&&124&&$(8s^7+28s^6+60s^5+72s^4+60s^3+30s^2+10s+2)(s+1)^{-8}$
&&1&&0&&0&&1&&1&&0&\cr
&&125&&$(12s^5+22s^4+32s^3+20s^2+8s+2)(s+1)^{-7}$
&&0&&0&&0&&0&&0&&0&\cr
\noalign{\hrule}}}}
\centerline
{\vbox
{\offinterlineskip
\halign{\strut#&\vrule#&\hfil#\hfil&\vrule#&\hfil#\hfil&\vrule#&\hfil#\hfil&
\vrule#&\hfil#\hfil&\vrule#&\hfil#\hfil&\vrule#&\hfil#\hfil&
\vrule#&\hfil#\hfil&\vrule#&\hfil#\hfil&\vrule#&\hfil#\hfil&
\vrule#\cr
\noalign{\hrule}
&&{\bf Order}&& {\bf Characteristic}
&&[3,1]&&[5,1]&&[5,2]&&[7,1]&&[7,2]&&[7,3]&\cr
\noalign {\hrule}
&&126&&$(3s^6+17s^5+34s^4+39s^3+9s^2)(s+1)^{-7}$&&0&&1&&0&&0&&0&&0&\cr
&&127&&$(4s^8+24s^7+70s^6+108s^5+110s^4+75s^3+34s^2+9s+1)(s+1)^{-9}$
&&1&&0&&0&&0&&0&&0&\cr
&&128&&$(2s^6+14s^5+35s^4+43s^3+29s^2+6s)(s+1)^{-7}$&&0&&0&&1&&0&&0&&0&\cr
&&129&&$(3s^4+16s^3+23s^2+6s)(s+1)^{-6}$&&0&&0&&0&&1&&1&&0&\cr
&&130&&$(2s^8+16s^7+52s^6+100s^5+112s^4+80s^3+30s^2+4s)(s+1)^{-9}$
&&1&&0&&0&&0&&0&&0&\cr
&&131&&$(4s^8+26s^7+70s^6+108s^5+103s^4+65s^3+27s^2+7s+1)(s+1)^{-9}$
&&1&&0&&0&&0&&0&&0&\cr
&&132&&$(4s^9+28s^8+87s^7+159s^6+189s^5+159s^4+87s^3+28s^2+4s)(s+1)^{-10}$
&&0&&0&&0&&0&&0&&0&\cr
&&133&&$(s^9+8s^8+27s^7+69s^6+110s^5+131s^4+96s^3+42s^2+8s)(s+1)^{-9}$
&&0&&0&&0&&0&&0&&0&\cr
&&134&&$(12s^8+64s^7+147s^6+205s^5+179s^4+105s^3+42s^2+10s+1)(s+1)^{-10}$
&&4&&0&&0&&0&&0&&0&\cr
&&135&&$(s^9+8s^8+26s^7+65s^6+100s^5+115s^4+83s^3+38s^2+8s)(s+1)^{-9}$
&&0&&1&&0&&0&&0&&0&\cr
&&136&&$(s^8+7s^7+18s^6+41s^5+43s^4+38s^3+8s^2)(s+1)^{-8}$
&&0&&0&&0&&0&&0&&0&\cr
&&137&&$(8s^7+28s^6+64s^5+78s^4+66s^3+32s^2+10s+2)(s+1)^{-8}$
&&1&&1&&1&&0&&1&&0&\cr
&&138&&$(6s^8+28s^7+66s^6+103s^5+97s^4+62s^3+26s^2+7s+1)(s+1)^{-9}$
&&1&&0&&0&&0&&0&&0&\cr
&&139&&$(4s^8+22s^7+62s^6+103s^5+115s^4+84s^3+39s^2+10s+1)(s+1)^{-9}$
&&0&&0&&0&&0&&0&&1&\cr
&&140&&$(s^9+8s^8+28s^7+73s^6+120s^5+147s^4+109s^3+46s^2+8s)(s+1)^{-9}$
&&1&&0&&0&&0&&0&&0&\cr
&&141&&$(s^8+7s^7+20s^6+47s^5+57s^4+56s^3+16s^2)(s+1)^{-}$
&&0&&0&&0&&0&&0&&0&\cr
&&142&&$(2s^9+19s^8+64s^7+123s^6+169s^5+157s^4+100s^3+38s^2+6s)(s+1)^{-10}$
&&0&&0&&0&&0&&0&&0&\cr
&&143&&$(s^9+9s^8+34s^7+79s^6+129s^5+137s^4+100s^3+42s^2+9s)(s+1)^{-9}$
&&1&&0&&0&&0&&1&&0&\cr
&&144&&$(s^{10}+9s^9+40s^8+113s^7+210s^6+269s^5+236s^4+138s^3+48s^2+7s)
(s+1)^{-10}$&&1&&0&&1&&0&&1&&1&\cr
&&145&&$(s^8+7s^7+19s^6+44s^5+50s^4+47s^3+12s^2)(s+1)^{-8}$
&&1&&0&&1&&1&&0&&1&\cr
&&146&&$(4s^8+26s^7+74s^6+113s^5+113s^4+76s^3+34s^2+9s+1)(s+1)^{-9}$
&&1&&0&&1&&0&&0&&0&\cr
&&147&&$(2s^9+13s^8+39s^7+79s^6+108s^5+111s^4+75s^3+35s^2+6s)(s+1)^{-9}$
&&1&&0&&1&&0&&1&&1&\cr
&&148&&$(s^8+8s^7+25s^6+65s^5+91s^4+89s^3+48s^2+15s)(s+1)^{-8}$
&&1&&1&&0&&0&&1&&0&\cr
&&149&&$(2s^9+13s^8+47s^7+102s^6+154s^5+157s^4+112s^3+49s^2+12s)(s+1)^{-9}$
&&1&&0&&0&&0&&1&&0&\cr
&&150&&$(2s^9+15s^8+51s^7+109s^6+156s^5+158s^4+110s^3+53s^2+15s+3)(s+1)^{-9}$
&&1&&0&&0&&1&&1&&0&\cr
&&151&&$(4s^9+26s^8+84s^7+164s^6+216s^5+192s^4+118s^3+46s^2+12s+2)(s+1)^{-9}$
&&1&&1&&0&&0&&1&&0&\cr
&&152&&$(s^{10}+10s^9+43s^8+137s^7+288s^6+413s^5+411s^4+281s^3+127s^2+34s+4)
(s+1)^{-10}$&&1&&0&&0&&0&&1&&0&\cr
&&153&&$(8s^7+32s^6+74s^5+90s^4+74s^3+34s^2+102+2)(s+1)^{-8}$
&&1&&0&&0&&1&&1&&0&\cr
&&154&&$(3s^8+26s^7+83s^6+150s^5+171s^4+128s^3+60s^2+16s+2)(s+1)^{-9}$
&&1&&0&&0&&0&&1&&0&\cr
&&155&&$(2s^9+14s^8+51s^7+109s^6+154s^5+141s^4+84s^3+27s^2+3s)(s+1)^{-9}$
&&1&&0&&0&&0&&1&&0&\cr
&&156&&$(s^9+9s^8+33s^7+95s^6+163s^5+179s^4+124s^3+48s^2+8s)(s+1)^{-9}$
&&1&&1&&1&&0&&1&&0&\cr
&&157&&$(8s^7+32s^6+70s^5+88s^4+74s^3+34s^2+10s+2)(s+1)^{-8}$
&&0&&0&&0&&0&&0&&0&\cr
&&158&&$(16s^5+28s^4+42s^3+24s^2+8s+2)(s+1)^{-7}$&&0&&0&&0&&0&&0&&0&\cr
&&159&&$(s^7+7s^6+18s^5+47s^4+42s^3+24s^2+8s+2)(s+1)^{-7}$
&&0&&1&&1&&0&&0&&0&\cr
&&160&&$(s^9+7s^8+24s^7+61s^6+100s^5+125s^4+99s^3+51s^2+12s)(s+1)^{-9}$
&&0&&0&&1&&0&&0&&0&\cr
&&161&&$(8s^7+28s^6+68s^5+80s^4+66s^3+32s^2+10s+2)(s+1)^{-8}$
&&0&&0&&0&&1&&0&&1&\cr
&&162&&$(3s^7+12s^6+30s^5+39s^4+36s^3+15s^2+6s)(s+1)^{-7}$
&&0&&0&&0&&0&&0&&0&\cr
&&163&&$(6s^5+9s^4+21s^3+9s^2+6s)(s+1)^{-6}$&&0&&0&&0&&0&&0&&1&\cr
&&164&&$(2s^8+10s^7+32s^6+64s^5+72s^4+58s^3+24s^2+8s)(s+1)^{-8}$
&&1&&0&&0&&0&&1&&0&\cr
&&165&&$(s^7+7s^6+17s^5+41s^4+35s^3+25s^2)(s+1)^{-7}$&&1&&0&&0&&0&&1&&1&\cr
&&166&&$(9s^5+15s^4+27s^3+15s^2+9s)(s+1)^{-6}$&&0&&0&&0&&0&&0&&0&\cr
&&167&&$(3s^7+20s^6+53s^5+84s^4+70s^3+34s^2+6s)(s+1)^{-8}$
&&1&&0&&0&&1&&0&&0&\cr
&&168&&$(3s^6+17s^5+34s^4+39s^3+9s^2)(s+1)^{-7}$&&0&&0&&0&&0&&0&&0&\cr
&&169&&$(2s^9+13s^8+39s^7+79s^6+108s^5+111s^4+75s^3+35s^2+6s)(s+1)^{-9}$
&&1&&0&&0&&0&&1&&0&\cr
&&170&&$(4s^9+26s^8+84s^7+168s^6+226s^5+204s^4+126s^3+48s^2+12s+2)(s+1)^{-9}$
&&1&&0&&0&&0&&1&&0&\cr
&&171&&$(8s^7+32s^6+74s^5+94s^4+80s^3+36s^2+10s+2)(s+1)^{-8}$
&&0&&0&&0&&0&&0&&0&\cr
&&172&&$(8s^7+28s^6+64s^5+82s^4+72s^3+34s^2+10s+2)(s+1)^{-8}$
&&0&&0&&0&&0&&0&&0&\cr
&&173&&$(s^{10}+10s^9+41s^8+120s^7+238s^6+350s^5+364s^4+260s^3+121s^2+33s+4)
(s+1)^{-10}$&&0&&0&&0&&0&&0&&0&\cr
&&174&&$(s^9+9s^8+31s^7+81s^6+131s^5+156s^4+113s^3+46s^2+8s)(s+1)^{-9}$
&&1&&1&&1&&0&&0&&1&\cr
&&175&&$(s^9+9s^8+30s^7+76s^6+119s^5+135s^4+96s^3+42s^2+8s)(s+1)^{-9}$
&&0&&1&&0&&0&&0&&0&\cr
&&176&&$(s^8+8s^7+22s^6+52s^5+61s^4+56s^3+16s^2)(s+1)^{-8}$
&&4&&0&&1&&1&&0&&1&\cr
&&177&&$(2s^7+15s^6+28s^5+48s^4+40s^3+29s^2+6s)(s+1)^{-8}$
&&0&&1&&0&&0&&0&&0&\cr
&&178&&$(s^5+22s^4+40s^3+24s^2+8s+2)(s+1)^{-7}$&&1&&1&&0&&0&&0&&0&\cr
&&179&&$(s^8+8s^7+21s^6+48s^5+53s^4+43s^3+12s^2)(s+1)^{-8}$
&&0&&0&&0&&0&&0&&0&\cr
&&180&&$(9s^5+15s^4+33s^3+15s^2+9s)
(s+1)^{-6}$&&1&&0&&1&&0&&0&&0&\cr
&&181&&$(s^8+8s^7+28s^6+60s^5+78s^4+67s^3+33s^2+10s+2)(s+1)^{-8}$
&&0&&1&&0&&1&&0&&0&\cr
&&182&&$(s^9+9s^8+31s^7+84s^6+133s^5+129s^4+76s^3+24s^2+3s)(s+1)^{-9}$
&&0&&1&&0&&0&&0&&1&\cr
&&183&&$(s^{10}+8s^9+32s^8+88s^7+165s^6+224s^5+207s^4+120s^3+40s^2+6s)
(s+1)^{-10}$&&1&&0&&1&&1&&0&&0&\cr
&&184&&$(16s^5+28s^4+46s^3+26s^2+8s+2)(s+1)^{-7}$&&1&&0&&0&&0&&0&&0&\cr
&&185&&$(3s^7+12s^6+30s^5+42s^4+39s^3+15s^2+6s)(s+1)^{-7}$
&&0&&0&&0&&0&&0&&0&\cr
&&186&&$(6s^5+9s^4+24s^3+9s^2+6s)(s+1)^{-6}$&&1&&0&&0&&0&&0&&0&\cr
&&187&&$(s^9+9s^8+29s^7+73s^6+110s^5+121s^4+80s^3+29s^2+4s)(s+1)^{-9}$
&&0&&0&&1&&0&&0&&0&\cr
&&188&&$(s^8+8s^7+20s^6+46s^5+47s^4+38s^3+8s^2)(s+1)^{-8}$
&&0&&1&&0&&0&&0&&0&\cr
\noalign{\hrule}}}}
\centerline
{\vbox
{\offinterlineskip
\halign{\strut#&\vrule#&\hfil#\hfil&\vrule#&\hfil#\hfil&\vrule#&\hfil#\hfil&
\vrule#&\hfil#\hfil&\vrule#&\hfil#\hfil&\vrule#&\hfil#\hfil&
\vrule#&\hfil#\hfil&\vrule#&\hfil#\hfil&\vrule#&\hfil#\hfil&
\vrule#\cr
\noalign{\hrule}
&&{\bf Order}&& {\bf Characteristic}
&&[3,1]&&[5,1]&&[5,2]&&[7,1]&&[7,2]&&[7,3]&\cr
\noalign {\hrule}
&&189&&$(3s^7+10s^6+32s^5+53s^4+49s^3+28s^2)(s+1)^{-7}$&&0&&0&&0&&1&&0&&0&\cr
&&190&&$(s^9+9s^8+28s^7+68s^6+98s^5+100s^4+63s^3+25s^2+4s)(s+1)^{-9}$
&&1&&0&&0&&0&&1&&0&\cr
&&191&&$(s^8+8s^7+19s^6+42s^5+39s^4+25s^3+4s^2)(s+1)^{-8}$
&&0&&1&&1&&1&&0&&0&\cr
&&192&&$(s^7+7s^6+11s^5+25s^4+4s^3)(s+1)^{-7}$&&0&&0&&0&&0&&0&&0&\cr
&&193&&$(12s^8+53s^7+123s^6+173s^5+167s^4+107s^3+45s^2+10s)(s+1)^{-9}$
&&1&&0&&0&&0&&1&&1&\cr
&&194&&$(s^{10}+9s^9+45s^8+131s^7+249s^6+324s^5+298s^4+191s^3+82s^2+21s+2)
(s+1)^{-10}$&&1&&0&&0&&0&&1&&1&\cr
&&195&&$(21s^6+50s^5+75s^4+64s^3+35s^2+10s)(s+1)^{-8}$&&1&&0&&0&&0&&1&&0&\cr
&&196&&$(s^9+8s^8+36s^7+87s^6+132s^5+128s^4+83s^3+32s^2+8s+1)(s+1)^{-9}$
&&1&&0&&1&&0&&1&&0&\cr
&&197&&$(3s^8+15s^7+45s^6+81s^5+102s^4+81s^3+45s^2+15s+3)(s+1)^{-8}$
&&0&&0&&0&&0&&0&&0&\cr
&&198&&$(s^{10}+10s^9+40s^8+115s^7+224s^6+318s^5+323s^4+226s^3+103s^2+27s+3)
(s+1)^{-10}$&&0&&0&&0&&0&&0&&0&\cr
&&199&&$(3s^8+15s^7+45s^6+84s^5+108s^4+87s^3+48s^2+15s+3)(s+1)^{-8}$
&&0&&0&&0&&1&&0&&1&\cr
&&200&&$(6s^6+15s^5+39s^4+42s^3+33s^2+12s+3)(s+1)^{-7}$&&0&&0&&0&&0&&0&&0&\cr
&&201&&$(s^{10}+10s^9+39s^8+110s^7+209s^6+294s^5+296s^4+210s^3+100s^2+29s+4)
(s+1)^{-10}$&&0&&0&&0&&0&&0&&0&\cr
&&202&&$(s^9+9s^8+29s^7+73s^6+111s^5+124s^4+87s^3+38s^2+8s)(s+1)^{-9}$
&&0&&0&&0&&0&&0&&0&\cr
&&203&&$(2s^{10}+17s^9+63s^8+152s^7+255s^6+320s^5+295s^4+199s^3+91s^2+25s+3)
(s+1)^{-10}$&&1&&0&&0&&1&&1&&0&\cr
&&204&&$(6s^8+38s^7+91s^6+142s^5+141s^4+103s^3+50s^2+15s+2)(s+1)^{-9}$
&&0&&1&&0&&0&&0&&0&\cr
&&205&&$(6s^8+38s^7+91s^6+142s^5+141s^4+103s^3+50s^2+15s+2)(s+1)^{-9}$
&&0&&0&&0&&0&&0&&0&\cr
&&206&&$(8s^8+45s^7+115s^6+162s^5+139s^4+82s^3+31s^2+8s+1)(s+1)^{-9}$
&&1&&0&&0&&0&&0&&0&\cr
&&207&&$(4s^9+22s^8+70s^7+142s^6+200s^5+188s^4+122s^3+48s^2+12s+2)(s+1)^{-9}$
&&1&&0&&0&&0&&1&&0&\cr
&&208&&$(8s^7+24s^6+62s^5+82s^4+76s^3+36s^2+10s+2)(s+1)^{-8}$
&&0&&0&&1&&0&&0&&0&\cr
&&209&&$(3s^8+12s^7+45s^6+87s^5+114s^4+87s^3+45s^2+12s+3)(s+1)^{-8}$
&&0&&0&&0&&0&&0&&0&\cr
&&210&&$(24s^6+52s^5+72s^4+56s^3+28s^2+8s)(s+1)^{-8}$&&1&&1&&0&&1&&1&&0&\cr
&&211&&$(7s^8+29s^7+77s^6+127s^5+139s^4+101s^3+49s^2+15s+2)(s+1)^{-9}$
&&1&&1&&0&&0&&0&&0&\cr
&&212&&$(2s^8+10s^7+34s^6+68s^5+74s^4+60s^3+20s^2+8s)(s+1)^{-8}$
&&0&&1&&0&&0&&0&&0&\cr
&&213&&$(2s^{10}+17s^9+63s^8+150s^7+246s^6+301s^5+269s^4+175s^3+78s^2+22s+3)
(s+1)^{-10}$&&0&&0&&0&&0&&0&&0&\cr
&&214&&$(2s^9+15s^8+46s^7+91s^6+118s^5+113s^4+67s^3+25s^2+3s)(s+1)^{-9}$
&&0&&0&&0&&1&&0&&1&\cr
&&215&&$(3s^8+15s^7+48s^6+87s^5+111s^4+87s^3+48s^2+15s+3)(s+1)^{-8}$
&&0&&1&&0&&0&&0&&0&\cr
&&216&&$(9s^5+12s^4+33s^3+12s^2+9s)(s+1)^{-6}$&&0&&0&&0&&0&&0&&0&\cr
&&217&&$(2s^{10}+17s^9+63s^8+152s^7+255s^6+320s^5+295s^4+199s^3+91s^2+25s+3)
(s+1)^{-10}$&&1&&1&&0&&0&&1&&1&\cr
&&218&&$(2s^{10}+17s^9+70s^8+172s^7+284s^6+332s^5+277s^4+166s^3+67s^2+17s+2)
(s+1)^{-10}$&&0&&6&&0&&0&&0&&1&\cr
&&219&&$(s^9+9s^8+30s^7+76s^6+120s^5+138s^4+103s^3+51s^2+12s)(s+1)^{-9}$
&&4&&0&&1&&0&&8&&1&\cr
&&220&&$(2s^9+15s^8+46s^7+93s^6+125s^5+125s^4+81s^3+35s^2+6s)(s+1)^{-9}$
&&0&&0&&1&&0&&0&&0&\cr
&&221&&$(2s^{10}+17s^9+68s^8+159s^7+244s^6+260s^5+191s^4+95s^3+28s^2+3s)
(s+1)^{-10}$&&0&&1&&1&&1&&0&&0&\cr
&&222&&$(s^9+8s^8+27s^7+69s^6+111s^5+134s^4+103s^3+51s^2+12s)(s+1)^{-9}$
&&0&&0&&0&&1&&0&&0&\cr
&&223&&$(s^9+12s^8+49s^7+117s^6+167s^5+150s^4+82s^3+27s^2+4s)(s+1)^{-9}$
&&1&&1&&0&&0&&0&&0&\cr
&&224&&$(8s^7+24s^6+66s^5+80s^4+70s^3+34s^2+10s+2)(s+1)^{-8}$
&&0&&0&&0&&0&&0&&0&\cr
&&225&&$(2s^9+21s^8+71s^7+137s^6+188s^5+174s^4+110s^3+41s^2+6s)(s+1)^{-10}$
&&1&&0&&0&&1&&0&&0&\cr
&&226&&$(3s^7+9s^6+30s^5+42s^4+42s^3+15s^2+6s)(s+1)^{-7}$
&&0&&0&&0&&0&&0&&0&\cr
&&227&&$(s^7+7s^6+17s^5+41s^4+38s^3+40s^2)(s+1)^{-7}$&&1&&1&&0&&1&&0&&0&\cr
&&228&&$(3s^8+15s^7+48s^6+84s^5+105s^4+84s^3+48s^2+15s+3)(s+1)^{-8}$
&&4&&0&&0&&0&&8&&0&\cr
&&229&&$(3s^{10}+24s^9+95s^8+222s^7+336s^6+343s^5+240s^4+119s^3+40s^2+9s+1)
(s+1)^{-10}$&&0&&0&&0&&8&&0&&0&\cr
&&230&&$(4s^9+22s^8+74s^7+152s^6+208s^5+186s^4+116s^3+46s^2+12s+2)(s+1)^{-9}$
&&0&&1&&0&&0&&0&&1&\cr
&&231&&$(8s^7+28s^6+68s^5+84s^4+72s^3+34s^2+10s+2)(s+1)^{-8}$
&&1&&6&&1&&0&&1&&0&\cr
&&232&&$(3s^8+15s^7+48s^6+87s^5+108s^4+81s^3+42s^2+12s+3)(s+1)^{-8}$
&&0&&0&&0&&0&&0&&0&\cr
&&233&&$(3s^8+15s^7+48s^6+93s^5+120s^4+93s^3+48s^2+15s+3)(s+1)^{-8}$
&&0&&0&&0&&0&&0&&0&\cr
&&234&&$(6s^6+15s^5+45s^4+48s^3+33s^2+12s+3)(s+1)^{-7}$&&1&&0&&0&&0&&0&&0&\cr
&&235&&$(6s^5+9s^4+30s^3+9s^2+6s)(s+1)^{-6}$&&0&&1&&0&&0&&0&&0&\cr
&&236&&$(2s^{10}+17s^9+63s^8+152s^7+255s^6+320s^5+295s^4+199s^3+91s^2+25s+3)
(s+1)^{-10}$&&1&&0&&0&&0&&1&&0&\cr
&&237&&$(2s^9+15s^8+46s^7+93s^6+125s^5+125s^4+81s^3+35s^2+6s)(s+1)^{-9}$
&&0&&0&&0&&0&&0&&0&\cr
&&238&&$(2s^9+15s^8+44s^7+86s^6+111s^5+106s^4+64s^3+22s^2+3s)(s+1)^{-9}$
&&0&&0&&0&&1&&0&&0&\cr
&&239&&$(2s^8+13s^7+29s^6+46s^5+41s^4+28s^3+3s^2)(s+1)^{-8}$
&&1&&0&&0&&0&&0&&0&\cr
&&240&&$(2s^8+21s^7+55s^6+105s^5+128s^4+114s^3+66s^2+22s+3)(s+1)^{-9}$
&&0&&0&&1&&0&&0&&0&\cr
&&241&&$(2s^7+19s^6+34s^5+56s^4+46s^3+29s^2+6s)(s+1)^{-8}$
&&0&&0&&0&&0&&0&&0&\cr
&&242&&$(2s^7+19s^6+32s^5+51s^4+39s^3+22s^2+3s)(s+1)^{-8}$
&&0&&0&&1&&0&&0&&0&\cr
&&243&&$(2s^6+17s^5+13s^4+25s^3+3s^2)(s+1)^{-7}$&&0&&1&&0&&0&&1&&0&\cr
&&244&&$(3s^8+15s^7+48s^6+90s^5+114s^4+84s^3+42s^2+12s+3)(s+1)^{-8}$
&&0&&0&&0&&1&&0&&0&\cr
&&245&&$(6s^6+15s^5+45s^4+45s^3+30s^2+9s+3)(s+1)^{-7}$&&1&&0&&0&&1&&0&&0&\cr
&&246&&$(16s^8+68s^7+158s^6+227s^5+221s^4+139s^3+55s^2+12s)(s+1)^{-9}$
&&0&&1&&0&&0&&0&&0&\cr
&&247&&$(2s^{10}+15s^9+55s^8+125s^7+200s^6+231s^5+200s^4+125s^3+55s^2+15s+2)
(s+1)^{-10}$&&0&&0&&0&&0&&0&&0&\cr
&&248&&$(2s^{10}+17s^9+61s^8+143s^7+234s^6+287s^5+259s^4+172s^3+78s^2+22s+3)
(s+1)^{-10}$&&1&&0&&0&&0&&1&&0&\cr
&&249&&$(2s^8+21s^7+53s^6+98s^5+116s^4+100s^3+56s^2+19s+3)(s+1)^{-9}$
&&1&&0&&0&&0&&1&&0&\cr
\noalign{\hrule}}}}
\vfil \eject}
\tenpoint
\centerline {ACKNOWLEDGMENTS}
\bigskip
I would like to thank the Brookhaven National Laboratory for their
hospitality during a part of this work, and in particular William Marciano
and Frank Paige for their help, especially in making possible the
use of the computer facilities of their group. I am also pleased to
acknowledge interesting discussions with Konstantinos Anagnostopoulos about
computer applications. Finally I would like to thank Wei Chen for carefully
reading through the manuscript, and for a number of helpful remarks and
suggestions.
\bigskip
\centerline {REFERENCES}
\bigskip
1) C. Aneziris, DESY preprint, {\it DESY 94-230}.
\par \smallskip
2) E. Witten, {\it Comm. Math. Phys.} {\bf 121} (1989) 351-399; V.F.R. Jones,
{\it Notices of AMS} {\bf 33} (1986) 219-225.
\par \smallskip
3) D. Kreimer, preprints {\it UTAS-PHYS-94-25} and {\it UTAS-PHYS-95-10}.
\par \smallskip
4) J.S. Birman {\it ``Braids, Links and Mapping Class Groups"}, Annals of
Mathematics Sudies No. 82, Princeton University Press, 1974.
\par \smallskip
5) F.Y. Wu,
{\it Rev. Mod. Phys.} {\bf 64} 1099-1132;
H. Temperley and E. Lieb,
{\it Proc. Roy. Soc. London} A {\bf 322}
(1971) 251-280;
R. Baxter,
{\it J. Phys. A} {\bf 20} (1987) 5241-5261 and
{\it Exactly solved models in Statistical Mechanics}, Academic
Press, 1982.
\par \smallskip
6) L. Smolin, {\it Cont. Math.} {\bf 71} (1988) 55-98; C. Rovelli and
L. Smolin, {\it Phys. Rev. Lett.} {\bf 61} (1988) 155.
\par \smallskip
7) D.W. Sumners, {\it Studies in Physical and Theoretical Chemistry}, {\bf 51}
(1987) 3-22.
\par \smallskip
8) G. Kolata, {\it Science} {\bf 231} (1986) 1506-1508; W.F. Pohl,
{\it Math.
Intelligencer} {\bf 12} (1980) 20-27; W.R. Bauer, F.H.C. Crick and J.H. White,
{\it Scientific American} {\bf 243} (1980) No. 1, 100-113; N.R. Cozzarelli,
S.J. Spengler and A. Stasiak, {\it Cell} {\bf 42} (1985) 325-334; D.W. Sumners,
{\it Math. Intelligencer} {\bf 12} (1990) 71-80; C. Ernst and D.W. Sumners,
{\it Math. Proc. Camb. Phil. Soc.} {\bf 108} (1990) 489-515.
\par \smallskip
9) D. Rolfsen, {\it ``Knots and Links"}, Berkeley, CA,: Publish or Perish,
Inc., 1976; G. Burde, H. Zieschang, {\it ``Knots"}, Berlin: de Gruyter, 1985;
L.H. Kauffman, {\it ``Knots and Physics"}, World Scientific Publishing Co.
Pte. Ltd., 1991.
\par \smallskip
10) M.B. Thistlethwaite, L.M.S. Lecture Notes no 93 pp. 1-76, Cambridge
University Press, 1985.
\par \smallskip
11) C.F. Gauss, {\it Werke VIII}, pp. 271-286; J.S. Carter, ``How Surfaces
Intersect in Space", World Scientific, 1993.
\par \smallskip
12) C. Jordan, {\it Course d' Analyse}, Paris, 1893; O. Veblen, {\it Trans.
Amer. Math. Soc.} {\bf 6} (1905) 83.
\par \smallskip
13) K. Reidemeister, {\it Abh. Math. Sem. Univ. Hamburg} {\bf 5} (1927) 7-23.
\par \smallskip
14) J.W. Alexander, {\it Trans. AMS} {\bf 30} (1928) 275-306.
\par \smallskip
15) P. Freyd, D. Yetter, J. Hoste, W.B.R. Lickorish, K.C. Millet, A. Ocneanu
{\it Bull. AMS} {\bf 12} (1985) 239-246; Y. Akutsu and M. Wadati, {\it J. Phys.
Soc. Jpn.}, {\bf 56} (1987) 839-842 and 3039-3051.
\par \smallskip
16) L.H. Kauffman, {\it Trans. AMS} {\bf 318} (1990) 417-471.
\par \smallskip
17) P. Ramadevi, T.R. Govindarajan and R.K. Kaul, preprint
{\it hep-th 9412084}.
\par \smallskip
18) J.S. Birman, {\it Invent. Math.} {\bf 81} (1985) 287-294.
\par \smallskip
19) C.M. Gordon and J. Luecke, University of Texas preprint, 1988;
B.A. Cipra, {\it Science} {\bf 241} (1988) 1291-1292.
\par \smallskip
20) W. Wirtinger, {\it Jahresber. DMV} {\bf 14} (1905) 517.
\vfil \eject
\centerline {$\underline {\rm FIGURE \qquad CAPTIONS}$}
\par \bigskip
Figure 1: a) Compactification of $S^1$, b) Compactification of $S^2$.
\par \medskip
Figure 2: a) Finite closed knot, b) Infinite open knot.
\par \medskip
Figure 3: Irregular knot projections.
\par \medskip
Figure 4: Identical shadows, distinct projections.
\par \medskip
Figure 5: Reidemeister moves.
\par \medskip
Figure 6: Naming a knot projection starting from $S$.
\par \medskip
Figure 7: An undrawable name.
\par \medskip
Figure 8: A loop drawn on a knot shadow.
\par \medskip
Figure 9: Non-tangent loops have an even number of intersecting points.
\par \medskip
Figure 10: Identical projections with distinct names.
\par \medskip
Figure 11: Irreducible knots are assigned successive natural numbers.
\par \medskip
Figure 12: Third moves that destroy the ability to perform first or second
moves.
\par \medskip
Figure 13: Third move $N_{\gamma } \rightarrow N$, followed by first move
$N \rightarrow N_{\alpha }$. $N_{\alpha }$ is preferred to $N$ which is
preferred to $N_{\gamma }$.
\par \medskip
Figure 14: A connected sum.
\par \medskip
Figure 15: By performing a second move a connected sum acquires a prime shadow.
\par \medskip
Figure 16: The three-color test is invariant under Reidemeister moves.
\par \medskip
Figure 17: It is possible to ``paint" the trefoil using all three colors.
\par \medskip
Figure 18: The projections related to skein relations.
\par \medskip
Figure 19: When using the Wirtinger presentation one has to distinguish
between these two possibilities.
\par \medskip
Figure 20: $a_k=f(a_i,a_j)$.
\par \medskip
Figure 21: Conditions required to make a color test invariant under
Reidemeister moves.
\par \medskip
Figures 22 and 23: Multiple colorisation possibilities of the knot bearing the
name $\{ (1,6),(3,8),(5,10),(7,12),(9,14),(11,16),(13,2),(15,4) \} $.
\par \medskip
Figure 24: The octahedron test, $M_{ij}$ is a $\pi \over 2$ rotation of $i$
with respect to $j$.
\par \medskip
Figure 25: A cube test is reduced to two tetrehedron tests.
\par \medskip
Figure 26: An octahedron test, where now $M_{ij}$ and $M_{ij'}$ are the
reflections of $i$ with respect to planes passing through $j$ and $j'$.
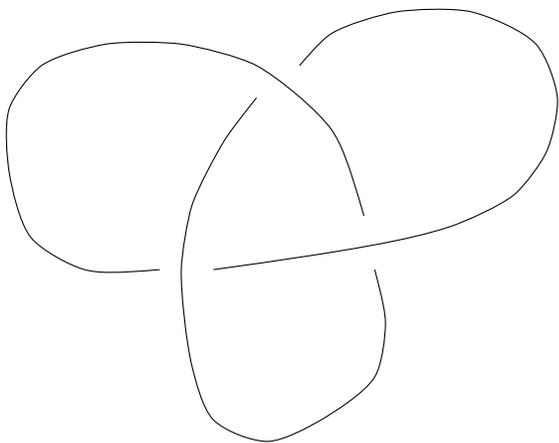
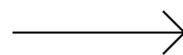
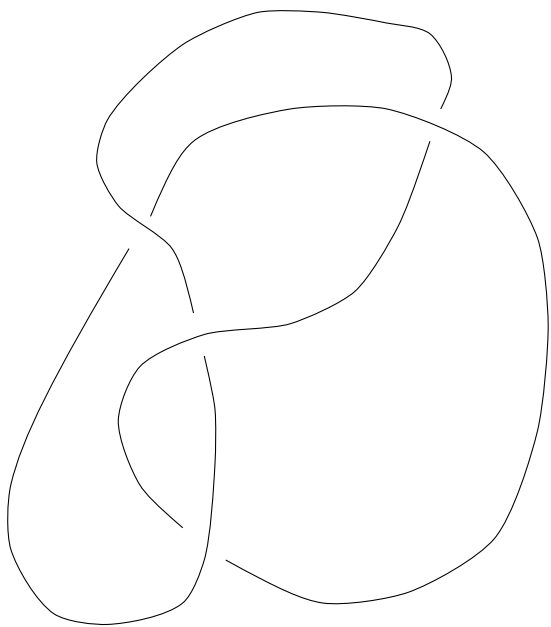
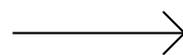
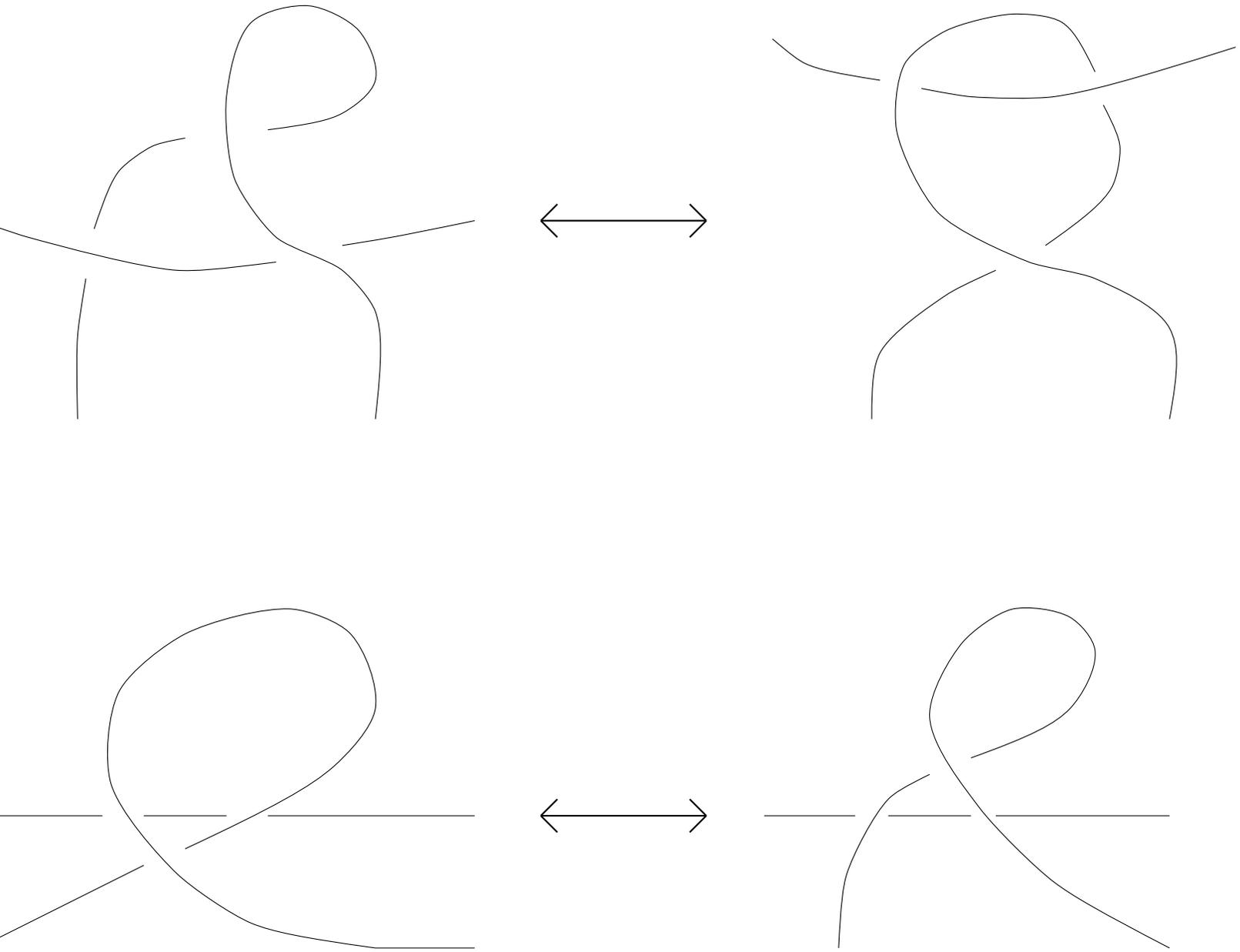
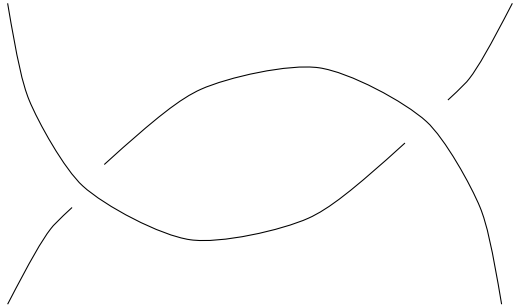
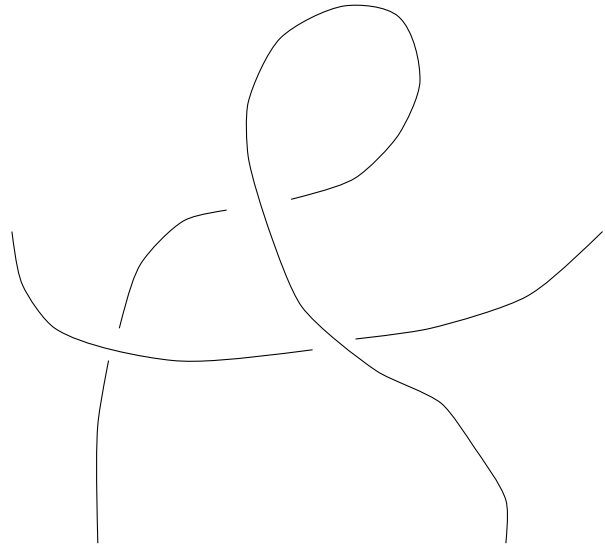
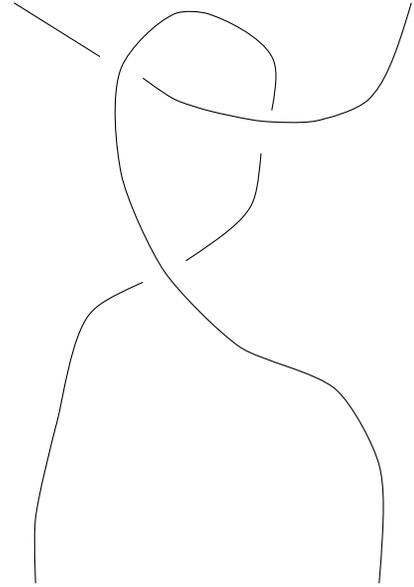
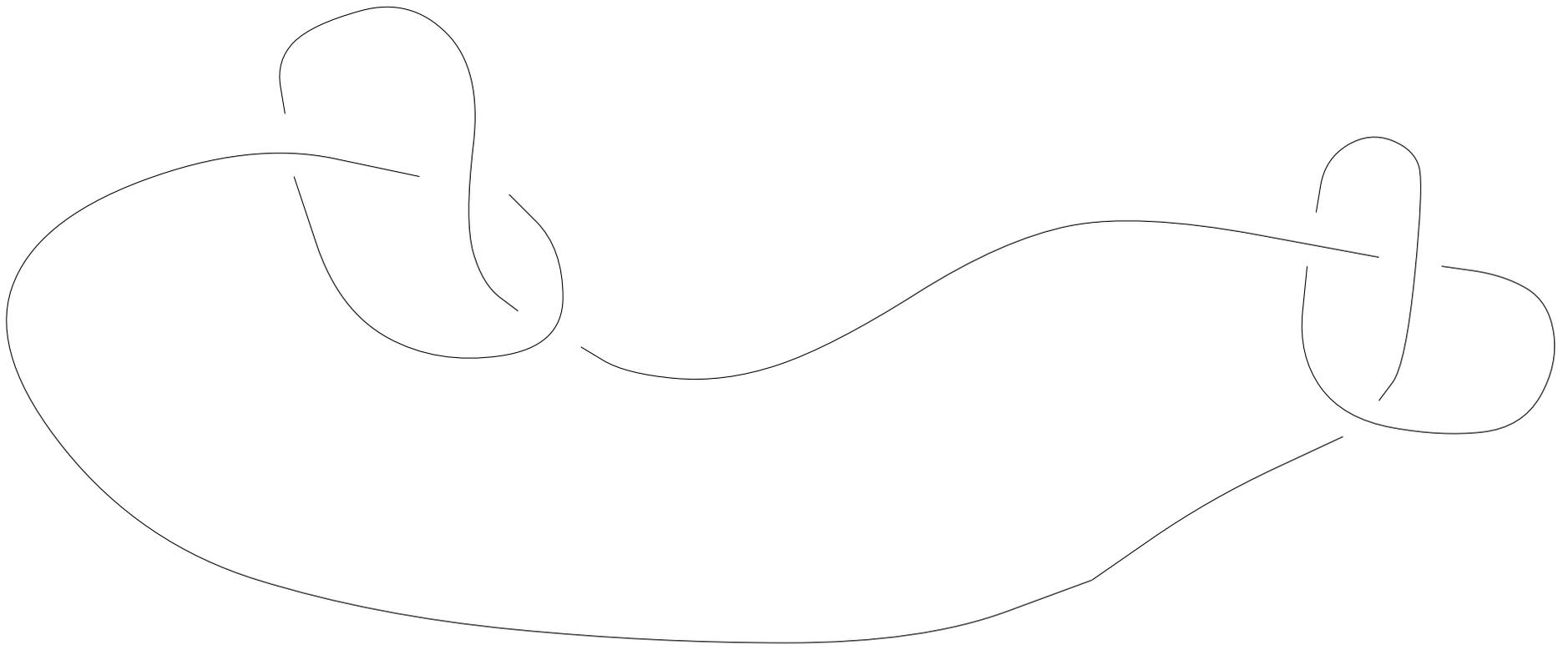
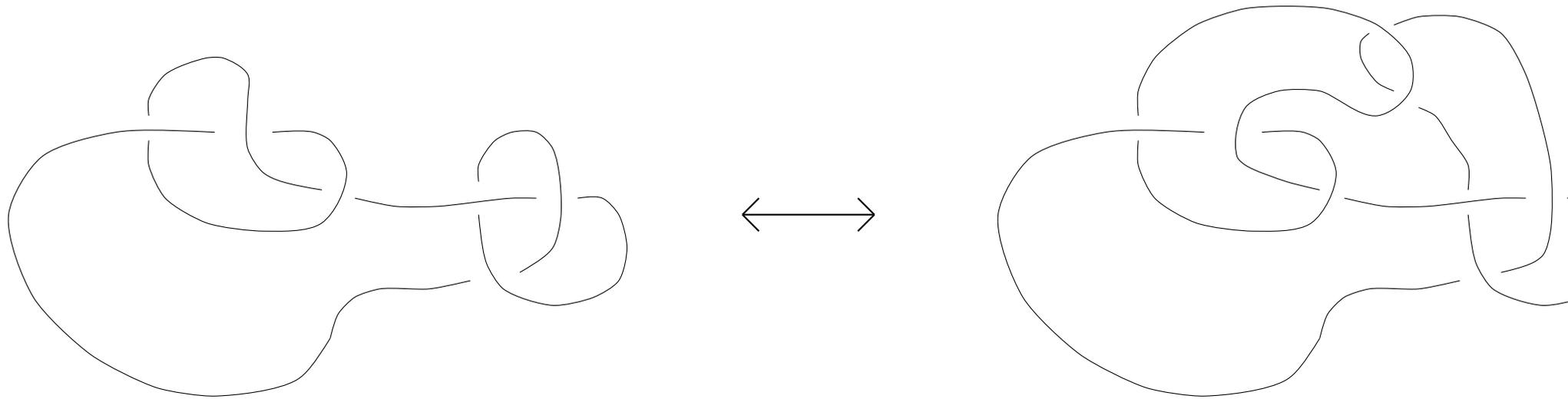
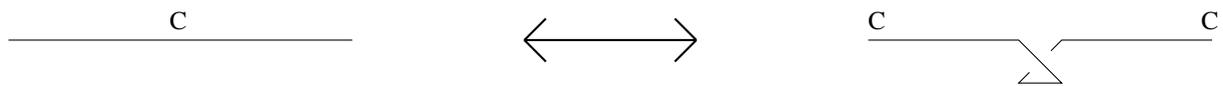
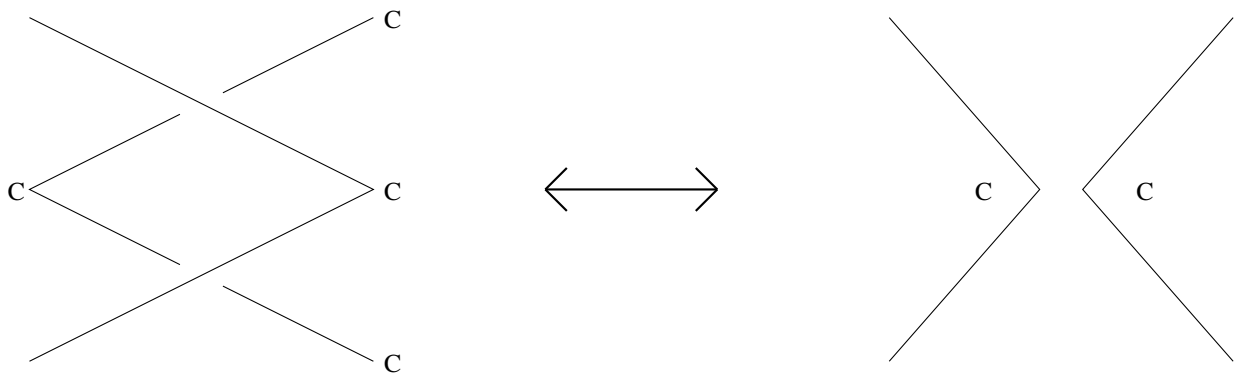
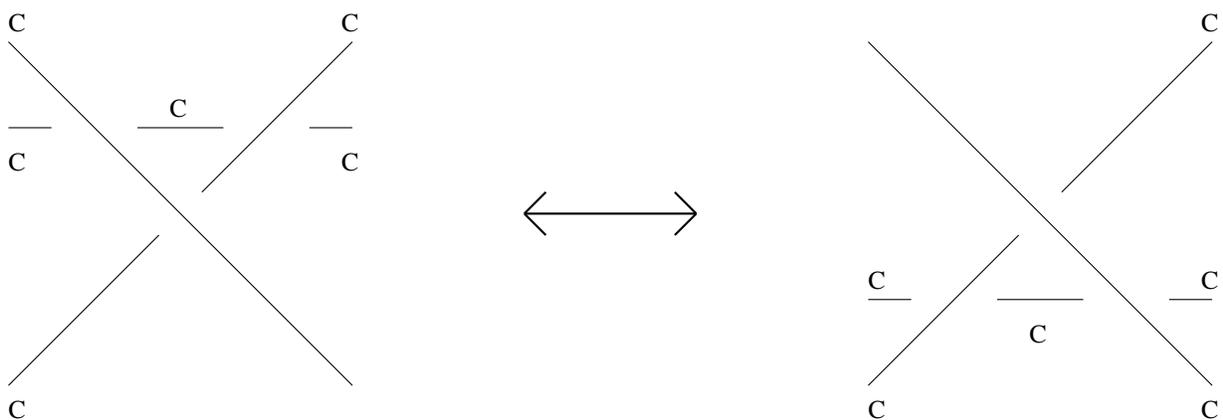
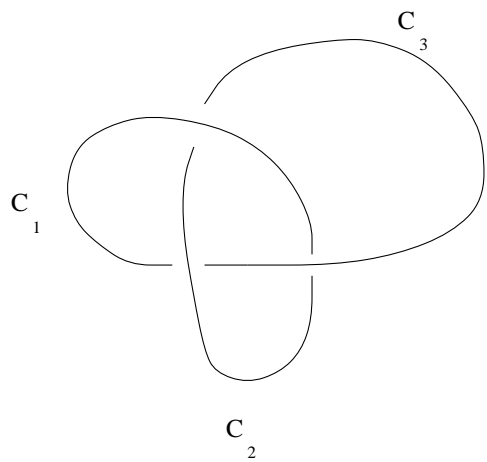
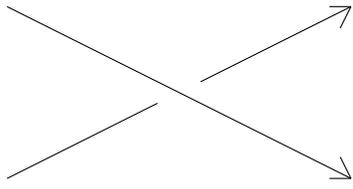
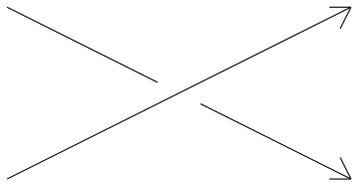
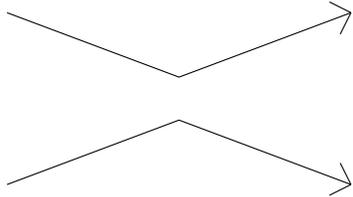
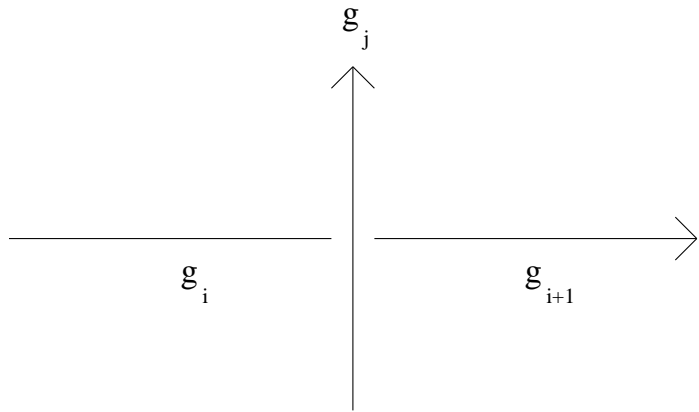
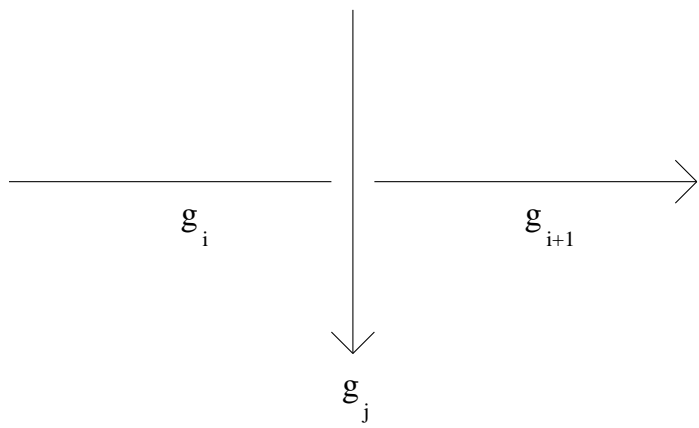
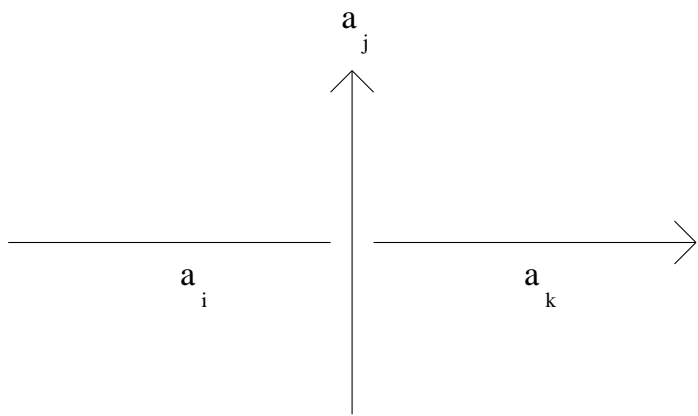
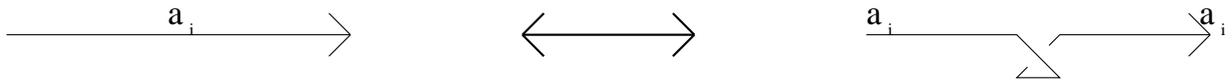
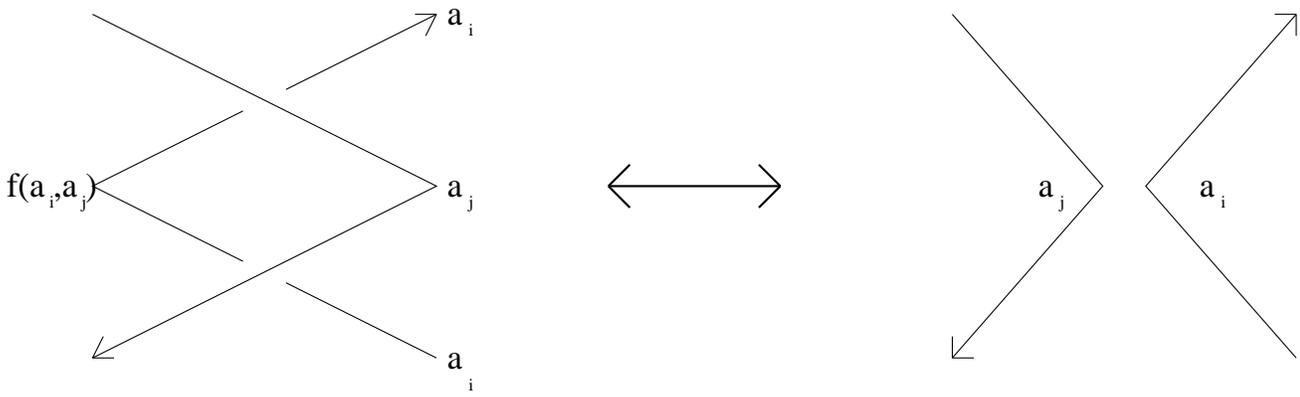
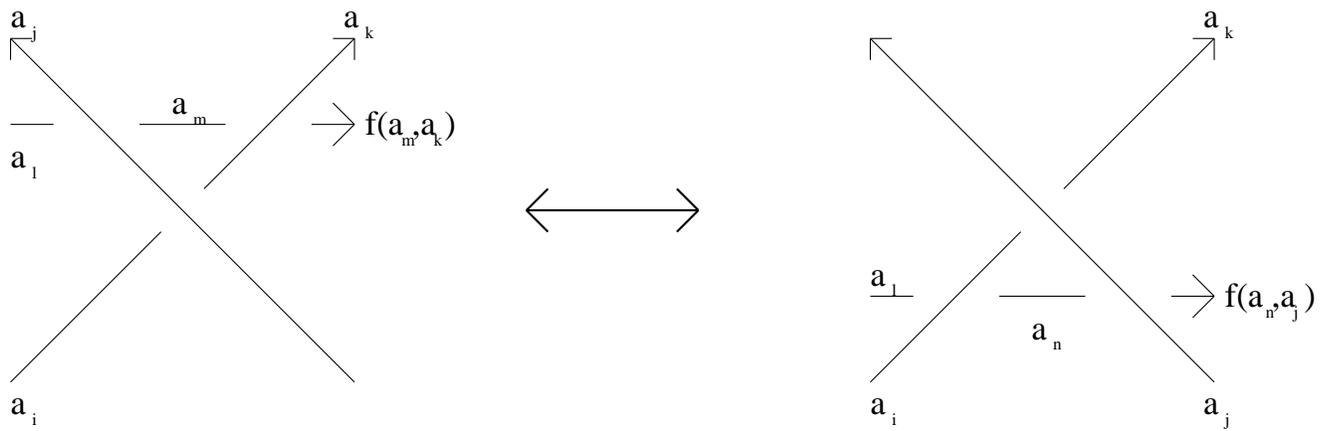
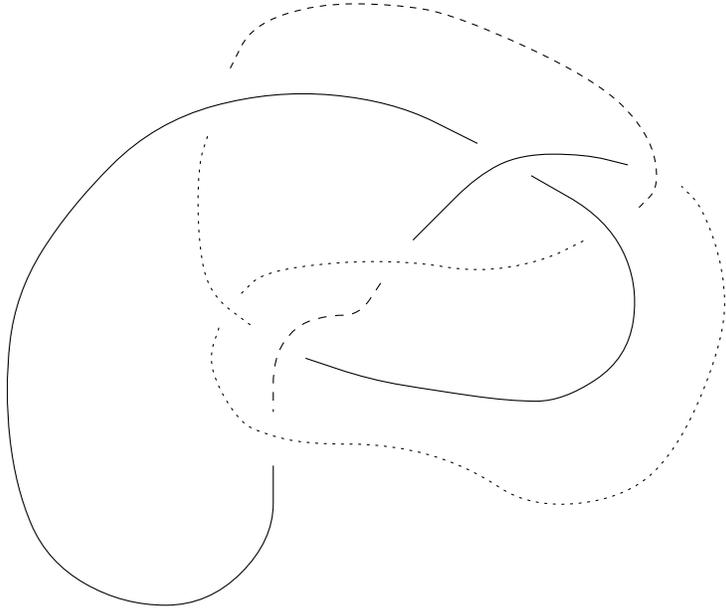
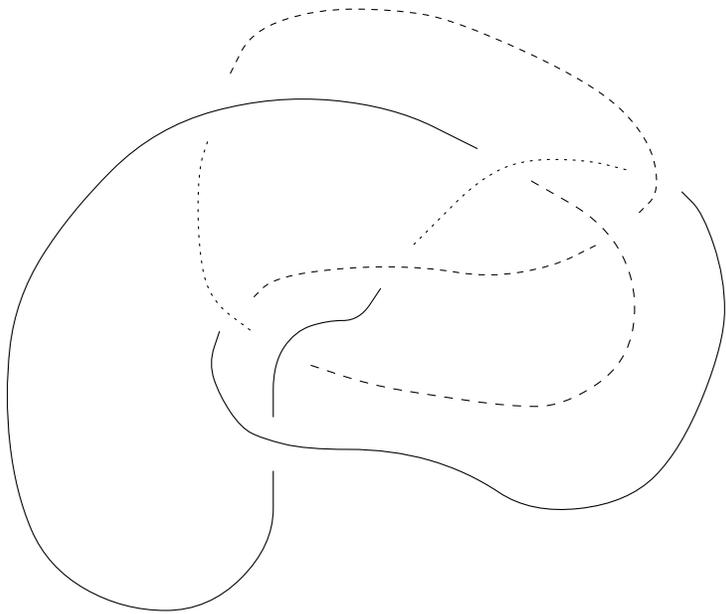
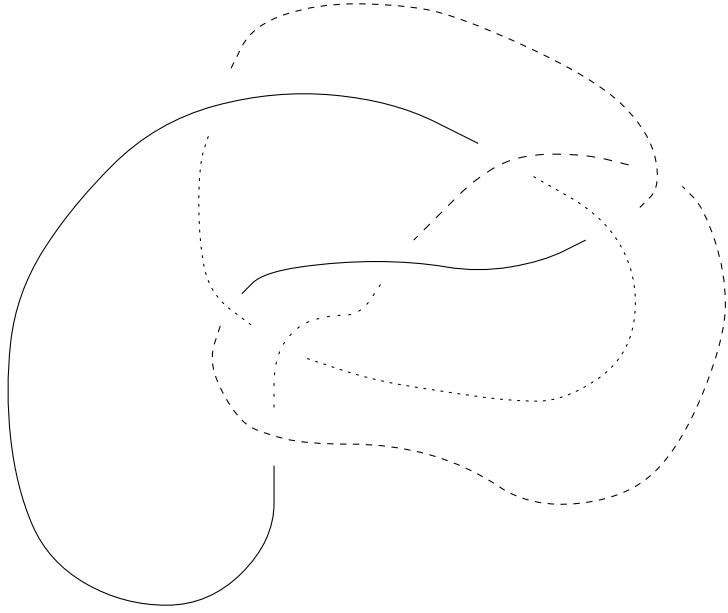
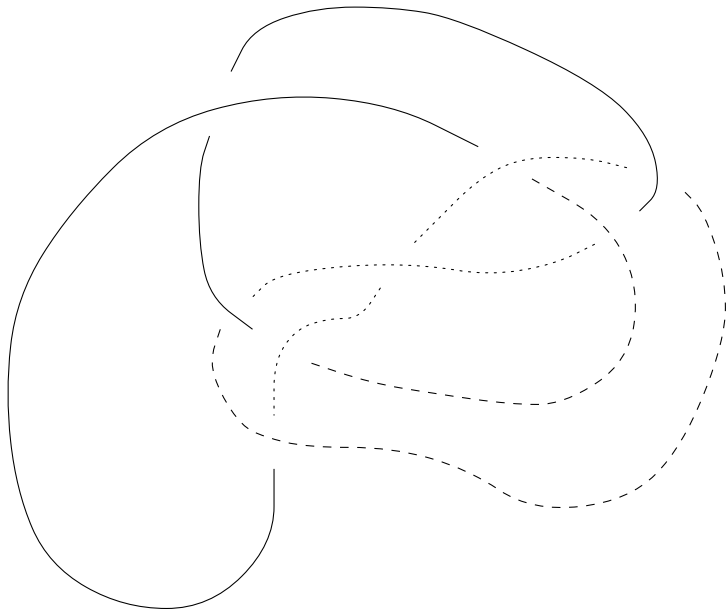
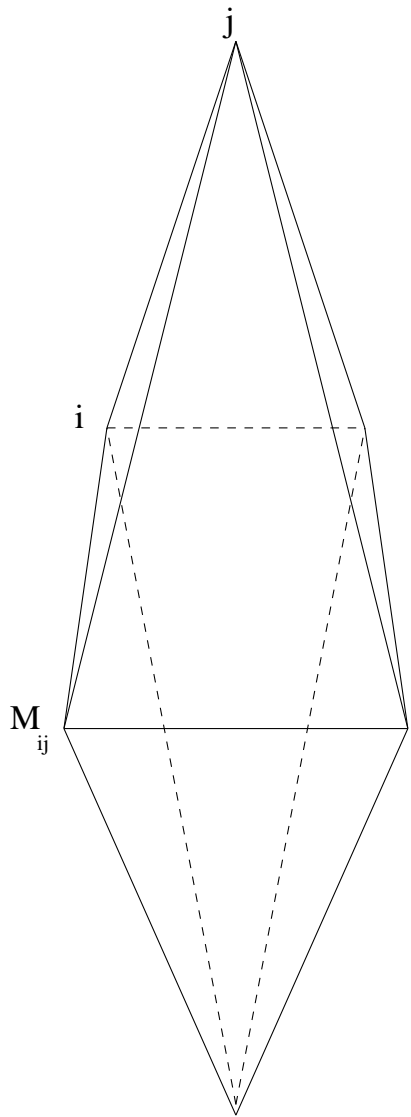
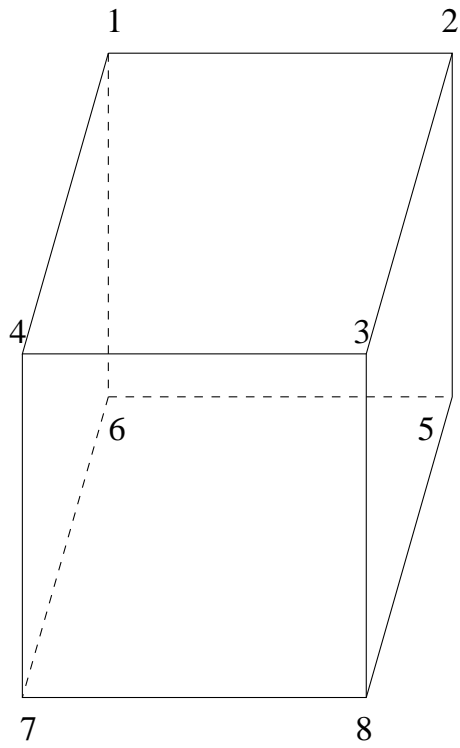
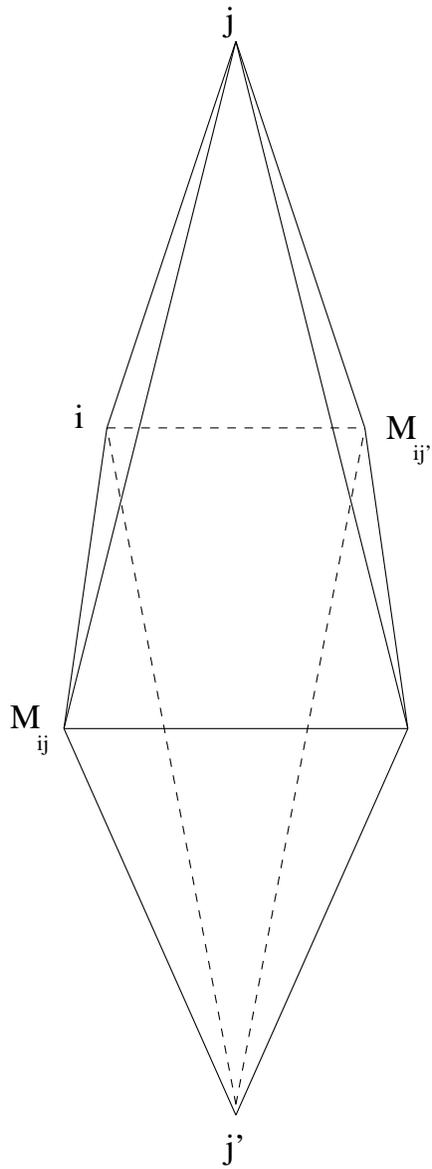
\end